\documentclass[12pt,twoside]{article}
\usepackage{amssymb}
\font\teneufm=eufm10 scaled \magstep1
\font\seveneufm=eufm7 scaled \magstep1
\font\fiveeufm=eufm5  scaled \magstep1
\newfam\eufmfam
\textfont\eufmfam=\teneufm
\scriptfont\eufmfam=\seveneufm
\scriptscriptfont\eufmfam=\fiveeufm
\def\frak#1{{\fam\eufmfam\relax#1}}

\newfam\msbfam
\font\tenmsb=msbm10 scaled \magstep1  \textfont\msbfam=\tenmsb
\font\sevenmsb=msbm7 scaled \magstep1 \scriptfont\msbfam=\sevenmsb
\font\fivemsb=msbm5 scaled \magstep1  \scriptscriptfont\msbfam=\fivemsb
\def\Bbb{\fam\msbfam \tenmsb}

\def\RR{{\Bbb R}}
\def\CC{{\Bbb C}}
\def\QQ{{\Bbb Q}}
\def\NN{{\Bbb N}}
\def\ZZ{{\Bbb Z}}

\def\PP{{\Bbb P}}

\def\ra{\rightarrow}

 \def\HollowBoxx #1#2#3{{\dimen0=#1 \advance\dimen0 by -#2
       \dimen1=#1 \advance\dimen1 by #3
        \vrule height 0pt depth #3 width #2
       \hskip -#3
       \vrule height #1 depth #3 width #3}}
 \def\LeftContraction{\mathord{\kern1.45pt \HollowBoxx{6pt}{3.5pt}{.4pt}}\,}

 \def\HollowBox #1#2#3{{\dimen0=#1 \advance\dimen0 by -#3
       \dimen1=#1 \advance\dimen1 by #3
        \vrule height #1 depth #3 width #3
        \vrule height 0pt depth #3 width #2
        \hskip -#3}}
 \def\RightContraction{\mathord{\, \HollowBox{6pt}{3.1pt}{.4pt}} \kern1.6pt}

\def\qed{{\hfill $\Box$}}
\newtheorem{theorem}{THEOREM}[section]

\newtheorem{lemma}[theorem]{Lemma}

\newtheorem{remark}[theorem]{Remark}
\newtheorem{proposition}[theorem]{Proposition}

\begin{document}

\begin{center}
{\Large \bf Hyperbolic 2-Dimensional Manifolds
\medskip\\
with 3-Dimensional Automorphism Group}
\medskip\\
\normalsize A. V. Isaev
\end{center}

\begin{quotation} \small \sl In this paper we determine all Kobayashi-hyperbolic 2-dimensional complex manifolds for which the group of holomorphic automorphisms has dimension 3. This work concludes a recent series of papers by the author on the classification of hyperbolic $n$-dimensional manifolds, with automorphism group of dimension at least $n^2-1$, where $n\ge 2$.
\end{quotation}

\thispagestyle{empty}

\pagestyle{myheadings}
\markboth{A. V. Isaev}{Hyperbolic Manifolds with 3-Dimensional Automorphism Group}

\setcounter{section}{-1}

\section{Introduction}
\setcounter{equation}{0}

If $M$ is a connected $n$-dimensional Kobayashi-hyperbolic complex manifold, then  the group $\hbox{Aut}(M)$ of holomorphic automorphisms of $M$ is a (real) Lie group in the compact-open topology, of dimension $d(M)$ not exceeding $n^2+2n$, with the maximal value occurring only for manifolds holomorphically equivalent to the unit ball $B^n\subset\CC^n$ \cite{Ko1}, \cite{Ka}. We are interested in describing hyperbolic manifolds with lower (but still sufficiently high) values of $d(M)$. The classification problem for hyperbolic manifolds with high-dimensional automorphism group is a complex-geometric analogue of that for Riemannian manifolds with high-dimensional isometry group, which inspired many results in the 1950's-70's (see \cite{Ko2} for details). The principal underlying property that made the classification in the Riemannian case  possible is that the group of isometries acts properly on the manifold -- see \cite{MS}, \cite{vDvdW} (a topological group $G$ is said to act properly on a manifold $S$ if the map $G\times S\ra S\times S$, $(g,p)\mapsto (gp,p)$ is proper). In the case of hyperbolic manifolds, the action of the group $\hbox{Aut}(M)$ is proper as well (see \cite{Ko1}, \cite{Ka}), and, as in the Riemannian case, this property is critical for our arguments, despite the fact that our techniques are almost entirely different from those utilized for isometry groups.    

In \cite{IKra}, \cite{I1} we completely classified manifolds with $n^2\le d(M)<n^2+2n$ (partial classifications for $d(M)=n^2$ were also obtained in \cite{GIK} and \cite{KV}). Note that for $d(M)=n^2$ the manifold $M$ may not be homogeneous, which makes this case substantially more difficult than the case $d(M)>n^2$, where  homogeneity always takes place \cite{Ka}. A further decrease in $d(M)$ almost immediately leads to unclassifiable situations. Indeed, no reasonable classification exists for $n=2$, $d(M)=2$, in which case $d(M)=n^2-2$ (observe, for example, that the automorphism group of a generic Reinhardt domain in $\CC^2$ is 2-dimensional). While it is possible that there is some classification for $n\ge 3$, $d(M)=n^2-2$, as well as for particular pairs $n$, $d(M)$ with $d(M)< n^2-2$ (see e.g. \cite{GIK} for a study of Reinhardt domains from the point of view of the automorphism group dimension), the case $d(M)=n^2-1$ is probably the only remaining candidate to investigate for the existence of an explicit classification for every $n\ge 2$. It turns out that all hyperbolic manifolds with $n\ge 2$, $d(M)=n^2-1$ indeed can be explicitly described. The case $n\ge 3$ was considered in \cite{I2}. The remaining case $n=2$, $d(M)=3$ is the subject of the present paper.

For brevity we call connected 2-dimensional hyperbolic manifolds with 3-dimensional automorphism group (2,3)-{\it manifolds}.  
For a (2,3)-manifold $M$, we work with the group $G(M):=\hbox{Aut}(M)^0$, the connected identity component of $\hbox{Aut}(M)$. Since the $G(M)$-action on $M$ is proper, for every $p\in M$ its isotropy subgroup $I_p:=\{f\in G(M): f(p)=p\}$ is compact in $G(M)$ and the orbit $O(p):=\{f(p):f\in G(M)\}$ is a connected closed submanifold in $M$. Proposition 2.1 of \cite{I2} gives (see Proposition \ref{dim} below) that for every $p\in M$ the orbit $O(p)$ has (real) codimension 1 or 2 in $M$, and in the latter case $O(p)$ is either a complex curve or a totally real submanifold.

We start by observing that the case when no codimension 1 orbits are present in the manifold can be dealt with as in \cite{I2} and leads to direct products $\Delta\times S$, where $\Delta$ is the unit disk in $\CC$ and $S$ is any hyperbolic Riemann surface with $d(S)=0$ (see Remark \ref{complexhypersurfaces}). Thus, from Section \ref{list} onwards we assume that a codimension 1 orbit is present in the manifold. Clearly, every codimension 1 orbit is either strongly pseudoconvex or Levi-flat.

In Section \ref{list} we give a large number of examples of (2,3)-manifolds. It will be shown in later sections that in fact these examples form a complete classification of (2,3)-manifolds with codimension 1 orbits.

In Section \ref{sectionstronglypseudoconvex} we deal with the case when every orbit is strongly pseudoconvex and classify all (2,3)-manifolds with this property in Theorem \ref{classrealhypersurfaces}. An important ingredient in the proof of Theorem \ref{classrealhypersurfaces} is E. Cartan's classification of 3-dimensional homogeneous strongly pseudoconvex $CR$-manifolds \cite{C}, together with the explicit determination of all covers of the non simply-connected hypersurfaces on Cartan's list \cite{I3}. The explicit realizations of the covers are important for our arguments throughout the paper, especially for those in the proof of Theorem \ref{theoremcodimtwo} in Section \ref{codimtwoorbits}. Another ingredient in the proof of Theorem \ref{classrealhypersurfaces} is an orbit gluing procedure that allows us to join strongly pseudoconvex orbits together to form (2,3)-manifolds.

Studying situations when Levi-flat and codimension 2 orbits can occur is perhaps the most interesting part of the paper. In Section \ref{sectleviflatorbits} we deal with Levi-flat orbits. Every such orbit is foliated by complex manifolds equivalent to $\Delta$ (see Proposition \ref{dim}). We describe Levi-flat orbits together with all possible actions of $G(M)$ in Proposition \ref{propleviflat} and use this description to classify in Theorem \ref{completerealhypersurface} all (2,3)-manifolds for which every orbit has codimension 1 and at least one orbit is Levi-flat. The proof of Theorem \ref{completerealhypersurface} uses the orbit gluing procedure introduced in the proof of Theorem \ref{classrealhypersurfaces}.

Finally, in Section \ref{codimtwoorbits} we allow codimension 2 orbits to be present in the manifold. Every complex curve orbit is equivalent to $\Delta$ (see Proposition \ref{dim}), whereas no a priori description of totally real orbits is available. The properness of the $G(M)$-action implies that there are at most two codimension 2 orbits in $M$ (see \cite{AA}), and in the proof of Theorem \ref{theoremcodimtwo} we investigate how one or two such orbits can be added to the previously obtained manifolds. This is done by studying complex curves invariant under the actions of the isotropy subgroups of points lying in codimension 2 orbits.

In fact, the arguments of the present paper yield not only a classification of (2,3)-manifolds as stated, but a classification of all connected 2-dimensional complex manifolds that admit a proper effective action of a 3-dimensional Lie group by holomorphic transformations. Some manifolds of this kind are not hyperbolic and were excluded during the course of proof (for instance, any 2-dimensional Hopf manifold admits an effective action of $SU_2$, but is clearly not hyperbolic). In addition, we ruled out those hyperbolic manifolds for which the automorphism group has dimension higher than 3 (the automorphism group of every such manifold has a closed 3-dimensional subgroup). Adding the excluded manifolds to our classification is straightforward and leads to a complete list of 2-dimensional complex manifolds with a proper effective action of a 3-dimensional Lie group by holomorphic transformations. We leave details to the reader.

Before proceeding, we would like to thank Stefan Nemirovski for many useful discussions and especially for showing us an elegant realization of a series of manifolds that appear in our classification.

\section{Initial Classification of Orbits}\label{orbits} 
\setcounter{equation}{0}

In this section we list some initial facts about $G(M)$-orbits that follow from the results of \cite{I2}. For $p\in M$ let $L_p:=\{d_pf: f\in I_p\}$ be the linear isotropy subgroup of $p$, where $d_pf$ is the differential of a map $f$ at $p$. The group $L_p$ is a compact subgroup of $GL(T_p(M),\CC)$ isomorphic to $I_p$ by means of the isotropy representation
$$
I_p\ra L_p, \quad f\mapsto d_pf,
$$
where $T_p(M)$ is the tangent space to $M$ at $p$. Proposition 2.1 of \cite{I2} implies the following

\begin{proposition}\label{dim} \sl Let $M$ be a (2,3)-manifold. Fix $p\in M$ and let $V_p:=T_p(O(p))$. Then the following holds:
\vspace{0cm}\\ 

\noindent (i) the orbit $O(p)$ is either a closed real hypersurface, or a closed complex curve, or a closed totally real 2-dimensional submanifold of $M$;
\vspace{0cm}\\

\noindent (ii) if $O(p)$ is a real Levi-flat hypersurface, it is foliated by complex curves holomorphically equivalent to $\Delta$, and there exist coordinates in $T_p(M)$ such that with respect to the orthogonal decomposition $T_p(M)=(V_p\cap iV_p)^{\bot}\oplus(V_p\cap iV_p)$ we have $L_p\subset\{\pm\hbox{id}\}\times L_p'$, where $L_p'$ is a finite subgroup of $U_1$;
\vspace{0cm}\\

\noindent (iii) if $O(p)$ is a complex curve, it is holomorphically equivalent to $\Delta$; furthermore, there exist coordinates $(z,w)$ in $T_p(M)$ in which $V_p=\{z=0\}$ and the identity component $L_p^0$ of $L_p$ is given by either the matrices
\begin{equation}
\left(\begin{array}{ll}
a^{\frac{k_1}{k_2}} & 0\\
0 & a
\end{array}\right),\label{stab1}
\end{equation}
for some  $k_1,k_2\in\ZZ$, $(k_1,k_2)=1$, $k_2\ne 0$, or the matrices
\begin{equation}
\left(\begin{array}{ll}
a & 0\\
0 & 1
\end{array}\right),\label{stab2}
\end{equation}
where $|a|=1$;
\vspace{0cm}\\

\noindent (iv) if $O(p)$ is totally real, then $T_p(M)=V_p\oplus iV_p$, and there are coordinates in $V_p$ such that every transformation from $L_p^0$ has the form: $v_1+iv_2\mapsto Av_1+iAv_2$, $v_1,v_2\in V_p$, where $A\in SO_2(\RR)$.
\end{proposition}

\begin{remark}\label{complexhypersurfaces} \rm Observe that if $O(p)$ is either a complex curve with $L_p^0$ given by (\ref{stab1}) for $k_1\ne 0$ or by (\ref{stab2}), or a totally real submanifold of $M$, then there exists a neighborhood $U$ of $p$ such that for every $q\in U\setminus O(p)$ the values at $q$ of the vector fields on $M$ arising from the action of $G(M)$, span a codimension 1 subspace of $T_q(M)$. Hence in this situation there is a codimension 1 orbit in $M$. Therefore, if no codimension 1 orbits are present in $M$, then every orbit is a complex curve with $L_p^0$ given by (\ref{stab1}) for $k_1=0$. In this case, arguing as in the proof of Proposition 4.1 of \cite{I2} we obtain that $M$ is holomorphically equivalent to a direct product $\Delta\times S$, where $S$ is a hyperbolic Riemann surface with $d(S)=0$.
\end{remark}

From now on we assume that a codimension 1 orbit is present in $M$.

\section{Examples of (2,3)-Manifolds}\label{list}
\setcounter{equation}{0}

In this section we give a large number of examples of (2,3)-manifolds. It will be shown in the forthcoming sections that these examples (upon  excluding equivalent manifolds) give a complete classification of (2,3)-manifolds with codimension 1 orbits.
\vspace{0.5cm}

\noindent {\bf (1)} In this example strongly pseudoconvex and Levi-flat orbits occur.
\vspace{0.2cm}

\noindent {\bf (a)} Fix $b\in\RR$, $b\ne 0,1$, and choose $0\le s<t\le\infty$ with either $s>0$ or $t<\infty$. Define
\begin{equation}
R_{b,s,t}:=\left\{(z,w)\in\CC^2: s\left(\hbox{Re}\,z\right)^b<\hbox{Re}\,w<t\left(\hbox{Re}\,z\right)^b,\,\hbox{Re}\,z>0\right\}.\label{usualrbst}
\end{equation}
The group $G(R_{b,s,t})=\hbox{Aut}(R_{b,s,t})$ consists of all maps
\begin{equation}
\begin{array}{lll}
z & \mapsto & \lambda z+i\beta,\\
w & \mapsto & \lambda^b w+i\gamma,
\end{array}
\label{auttau}
\end{equation}
where $\lambda>0$ and $\beta,\gamma\in\RR$. The $G(R_{b,s,t})$-orbits are the following pairwise $CR$-equivalent strongly pseudoconvex hypersurfaces:
$$
O^{R_b}_{\alpha}:=\left\{(z,w)\in\CC^2: \hbox{Re}\,w=\alpha\left(\hbox{Re}\,z\right)^b,\,\hbox{Re}\,z>0\right\},\quad s<\alpha<t,
$$
and we set 
\begin{equation}
\tau_b:=O^{R_b}_1.\label{tau}
\end{equation}
For every $b\in\RR$ we denote the group of maps of the form (\ref{auttau}) by $G_b$.
\vspace{0.3cm}

\noindent {\bf (b)} If in the definition of $R_{b,s,t}$ we let $-\infty\le s<0<t\le\infty$, where at least one of $s,t$ is finite, we again obtain a hyperbolic domain whose automorphism group coincides with $G_b$, unless $b=1/2$ and $t=-s$ (observe that $R_{1/2,s,-s}$ is equivalent to the unit ball). In such domains, in addition to the strongly pseudoconvex orbits $O_{\alpha}^{R_b}$ for suitable values of $\alpha$ (which are allowed to be negative), there is the following unique Levi-flat orbit:
\begin{equation}
{\cal O}_1:=\left\{(z,w)\in\CC^2: \hbox{Re}\,z>0,\, \hbox{Re}\,w=0\right\}.\label{calo1}
\end{equation}   
\vspace{0.3cm}

\noindent {\bf (c)} For $b>0$, $b\ne 1$, $-\infty<s<0<t<\infty$ define
$$
\hat R_{b,s,t}:=R_{b,s,\infty}\cup\left\{(z,w)\in\CC^2:\hbox{Re}\,w>t\left(-\hbox{Re}\,z\right)^b,\,\hbox{Re}\,z<0\right\}\cup\hat{\cal O}_1,
$$
where
\begin{equation}
\hat{\cal O}_1:=\left\{(z,w)\in\CC^2: \hbox{Re}\,z=0,\, \hbox{Re}\,w>0\right\}.\label{calo1prime}
\end{equation}
The group $G(\hat R_{b,s,t})=\hbox{Aut}(\hat R_{b,s,t})$ coincides with $G_b$, and, in addition to strongly pseudoconvex orbits $CR$-equivalent to $\tau_b$, the Levi-flat hypersurfaces ${\cal O}_1$ and $\hat{\cal O}_1$ are also $G_b$-orbits in $\hat R_{b,s,t}$.  
\vspace{0.5cm}

\noindent{\bf (2)} In this example strongly pseudoconvex orbits and a single Levi-flat orbit arise.   
\vspace{0,3cm}

\noindent {\bf (a)} For $0\le s<t\le\infty$ with either $s>0$ or $t<\infty$ define
\begin{equation}
\begin{array}{lll}
U_{s,t}&:=&\Bigl\{(z,w)\in\CC^2:\hbox{Re}\,w\cdot\ln\left(s\hbox{Re}\,w\right)<\hbox{Re}\,z<\\
&&\hspace{4cm}\hbox{Re}\,w\cdot\ln\left(t\hbox{Re}\,w\right),\,\hbox{Re}\,w>0\Bigr\}.
\end{array}\label{usualubst}
\end{equation}
The group $G(U_{s,t})=\hbox{Aut}(U_{s,t})$ consists of all maps
\begin{equation}
\begin{array}{lll}
z & \mapsto & \lambda z+(\lambda\ln\lambda)w+i\beta,\\
w & \mapsto & \lambda w+i\gamma,\\
\end{array}\label{autxi}
\end{equation}
where $\lambda>0$ and $\beta,\gamma\in\RR$. The $G(U_{s,t})$-orbits are the following pairwise $CR$-equivalent strongly pseudoconvex hypersurfaces:
$$
O^U_{\alpha}:=\left\{(z,w)\in\CC^2:\hbox{Re}\,z=\hbox{Re}\,w\cdot\ln\left(\alpha\hbox{Re}\,w\right),\,\hbox{Re}\,w>0\right\},\quad s<\alpha<t,
$$
and we set
\begin{equation}
\xi:=O^U_1.\label{xi}
\end{equation}
We denote the group of all maps of the form (\ref{autxi}) by ${\frak G}$.
\vspace{0.3cm}

\noindent {\bf (b)} For $-\infty<t<0<s<\infty$ define
$$
\hat U_{s,t}=U_{s,\infty}\cup \left\{(z,w)\in\CC^2:\hbox{Re}\,z>\hbox{Re}\,w\cdot\ln\left(t\hbox{Re}\,w\right),\,\hbox{Re}\,w<0\right\}\cup{\cal O}_1.
$$
The group $G(\hat U_{s,t})=\hbox{Aut}(\hat U_{s,t})$ coincides with ${\frak G}$, and, in addition to strongly pseudoconvex orbits $CR$-equivalent to $\xi$, the Levi-flat hypersurface ${\cal O}_1$ is also a ${\frak G}$-orbit in $\hat U_{s,t}$.
\vspace{0.5cm}

\noindent {\bf (3)} In this example strongly pseudoconvex orbits and a totally real orbit occur.
\vspace{0.3cm}

\noindent {\bf (a)} For $0\le s<t<\infty$ define
\begin{equation}
{\frak S}_{s,t}:=\left\{(z,w)\in\CC^2: s<\left(\hbox{Re}\,z\right)^2+\left(\hbox{Re}\,w\right)^2<t\right\}.\label{frakst}
\end{equation}
The group $G({\frak S}_{s,t})$ consists of all maps of the form
\begin{equation}
\left(
\begin{array}{l}
z\\
w
\end{array}
\right)\mapsto A
\left(
\begin{array}{l}
z\\
w
\end{array}
\right)+i
\left(
\begin{array}{l}
\beta\\
\gamma
\end{array}
\right),\label{autchi}
\end{equation}
where $A\in SO_2(\RR)$ and $\beta,\gamma\in\RR$. The $G({\frak S}_{s,t})$-orbits are the following pairwise $CR$-equivalent strongly pseudoconvex hypersurfaces:
$$
O^{\frak S}_{\alpha}:=\left\{(z,w)\in\CC^2: \left(\hbox{Re}\,z\right)^2+\left(\hbox{Re}\,w\right)^2=\alpha\right\},\quad s<\alpha<t,
$$
and we set
\begin{equation}
\chi:=O^{\frak S}_1.\label{chi}
\end{equation}
We denote the group of all maps of the form (\ref{autchi}) by ${\cal R}_{\chi}$.

\vspace{0.3cm}

\noindent {\bf (b)} For $0<t<\infty$ set
\begin{equation}
{\frak S}_t:=\left\{(z,w)\in\CC^2: \left(\hbox{Re}\,z\right)^2+\left(\hbox{Re}\,w\right)^2<t\right\}.\label{frakt}
\end{equation}
The group $G({\frak S}_t)$ coincides with ${\cal R}_{\chi}$, and, apart from strongly pseudoconvex orbits $CR$-equivalent to $\chi$, its action on ${\frak S}_t$ has the totally real orbit
\begin{equation}
{\cal O}_2:=\left\{(z,w)\in\CC^2:\hbox{Re}\,z=0,\,\hbox{Re}\,w=0\right\}.\label{calo2}
\end{equation}
\vspace{0.5cm}

\noindent {\bf (4)} In this example we explicitly describe all covers of the domains ${\frak S}_{s,t}$ and hypersurface $\chi$ introduced in {\bf (3)} (for more details see \cite{I3}). Only strongly pseudoconvex orbits occur here. 

Let $\Phi_{\chi}^{(\infty)}:\CC^2\ra\CC^2$ be the following map:
$$
\begin{array}{lll}
z & \mapsto & \exp\left(\hbox{Re}\,z\right)\cos\left(\hbox{Im}\,z\right)+i\hbox{Re}\,w,\\
w & \mapsto & \exp\left(\hbox{Re}\,z\right)\sin\left(\hbox{Im}\,z\right)+i\hbox{Im}\,w.
\end{array}
$$
It is easy to see that $\Phi_{\chi}^{(\infty)}$ is an infinitely-sheeted covering map onto $\CC^2\setminus\{\hbox{Re}\,z=0,\,\hbox{Re}\,w=0\}$. Introduce on the domain of $\Phi_{\chi}^{(\infty)}$ the complex structure defined by the condition that the map $\Phi_{\chi}^{(\infty)}$ is holomorphic (the pull-back complex structure under $\Phi_{\chi}^{(\infty)}$), and denote the resulting manifold by $M_{\chi}^{(\infty)}$. Next, for an integer $n\ge 2$, consider the map $\Phi_{\chi}^{(n)}$ from  $\CC^2\setminus\{\hbox{Re}\,z=0,\,\hbox{Re}\,w=0\}$ onto itself defined as follows:
\begin{equation}
\begin{array}{lll}
z & \mapsto & \hbox{Re}\Bigl(\left(\hbox{Re}\,z+i\hbox{Re}\,w\right)^n\Bigr)+i\hbox{Im}\,z,\\
w & \mapsto & \hbox{Im}\Bigl(\left(\hbox{Re}\,z+i\hbox{Re}\,w\right)^n\Bigr) +i\hbox{Im}\,w.
\end{array}\label{phichin}
\end{equation}
Denote by $M_{\chi}^{(n)}$ the domain of $\Phi_{\chi}^{(n)}$ with the pull-back complex structure under $\Phi_{\chi}^{(n)}$.

For $0\le s<t<\infty$, $n\ge 2$ define
\begin{equation}
\begin{array}{lll}
{\frak S}_{s,t}^{(n)}&:=&\left\{(z,w)\in M_{\chi}^{(n)}:s^{1/n}<\left(\hbox{Re}\,z\right)^2+\left(\hbox{Re}\,w\right)^2<t^{1/n}\right\},\\
{\frak S}_{s,t}^{(\infty)}&:=&\left\{(z,w)\in M_{\chi}^{(\infty)}:(\ln s)/2<\hbox{Re}\,z<(\ln t)/2\right\}.
\end{array}\label{frakstinftyn}
\end{equation}
The domains ${\frak S}_{s,t}^{(n)}$ and ${\frak S}_{s,t}^{(\infty)}$ are respectively an $n$- and infinite-sheeted cover of the domain ${\frak S}_{s,t}$. The group $G\left({\frak S}_{s,t}^{(n)}\right)$ for $n\ge 2$ consists of all maps 
\begin{equation}
\begin{array}{lll}
z & \mapsto & \cos\psi\cdot\hbox{Re}\,z+\sin\psi\cdot\hbox{Re}\,w+\\
&&i\Bigl(\cos(n\psi)\cdot\hbox{Im}\,z+\sin(n\psi)\cdot\hbox{Im}\,w+\beta\Bigr),\\
\vspace{0cm}&&\\
w & \mapsto & -\sin\psi\cdot\hbox{Re}\,z+\cos\psi\cdot\hbox{Re}\,w+\\
&&i\Bigl(-\sin(n\psi)\cdot\hbox{Im}\,z+\cos(n\psi)\cdot\hbox{Im}\,w+
\gamma\Bigr),
\end{array}\label{autchin1}
\end{equation}
where $\psi,\beta,\gamma\in\RR$. The $G\left({\frak S}_{s,t}^{(n)}\right)$-orbits are the following pairwise $CR$-equivalent strongly pseudoconvex hypersurfaces:
$$
O^{{\frak S}^{(n)}}_{\alpha}:=\left\{(z,w)\in M_{\chi}^{(n)}:\left(\hbox{Re}\,z\right)^2+\left(\hbox{Re}\,w\right)^2=\alpha\right\},\quad s^{1/n}<\alpha<t^{1/n},
$$
and we set 
\begin{equation}
\chi^{(n)}:=O^{{\frak S}^{(n)}}_1\label{chin}
\end{equation}
(this hypersurface is an $n$-sheeted cover of $\chi$).

The group $G\left({\frak S}_{s,t}^{(\infty)}\right)$ consists of all maps 
\begin{equation}
\begin{array}{ccc}
z & \mapsto & z+i\beta,\\
w & \mapsto & e^{i\beta}w+a,
\end{array}\label{autchiinfty1}
\end{equation}
where $\beta\in\RR$, $a\in\CC$. The $G\left({\frak S}_{s,t}^{(\infty)}\right)$-orbits are the following pairwise $CR$-equivalent strongly pseudoconvex hypersurfaces:
$$
O^{{\frak S}^{(\infty)}}_{\alpha}:=\left\{(z,w)\in M_{\chi}^{(\infty)}:\hbox{Re}\,z=\alpha\right\},\quad (\ln s)/2<\alpha<(\ln t)/2,
$$
and we set
\begin{equation}
\chi^{(\infty)}:=O^{{\frak S}^{(\infty)}}_0\label{chiinfty}
\end{equation}
(this hypersurface is an infinitely-sheeted cover of $\chi$). 
\vspace{0.5cm}

\noindent{\bf (5)} As in the preceding example, only strongly pseudoconvex orbits occur here. 

Fix $b>0$ and for $0<t<\infty$, $e^{-2\pi b}t<s<t$  consider the tube domain
\begin{equation}
V_{b,s,t}:=\left\{(z,w)\in\CC^2:se^{b\phi}<r<te^{b\phi}\right\},\label{usualvbst}
\end{equation}   
where $(r,\phi)$ denote the polar coordinates in the $(\hbox{Re}\,z,\hbox{Re}\,w)$-plane with $\phi$ varying from $-\infty$ to $\infty$ (thus, the boundary of $V_{b,t,s}\cap\RR^2$ consists of two spirals accumulating to the origin and infinity). The group $G(V_{b,s,t})=\hbox{Aut}(V_{b,s,t})$ consists of all maps of the form
\begin{equation}
\left(
\begin{array}{l}
z\\
w
\end{array}
\right)\mapsto e^{b\psi}
\left(
\begin{array}{rr}
\cos\psi & \sin\psi\\
-\sin\psi & \cos\psi
\end{array}
\right)
\left(
\begin{array}{l}
z\\
w
\end{array}
\right)+i
\left(
\begin{array}{l}
\beta\\
\gamma
\end{array}
\right),\label{autrho}
\end{equation}
where $\psi,\beta,\gamma\in\RR$. The $G(V_{b,s,t})$-orbits are the following pairwise $CR$-equivalent strongly pseudoconvex hypersurfaces:
$$
O^{V_b}_{\alpha}:=\left\{(z,w)\in\CC^2:r=\alpha e^{b\phi}\right\},\quad s<\alpha<t,
$$
and we set
\begin{equation}
\rho_b:=O^{V_b}_1.\label{rho}
\end{equation}
\vspace{0.5cm}

\noindent{\bf (6)} In this example strongly pseudoconvex orbits and a totally real orbit arise.
\vspace{0.3cm}

\noindent {\bf (a)} For $1\le s<t<\infty$ define
\begin{equation}
\begin{array}{ll}
E_{s,t}:=&\Bigl\{(\zeta:z:w)\in\CC\PP^2: s |\zeta^2+z^2+w^2|<|\zeta|^2+|z|^2+|w|^2<\\
&\hspace{6cm} t |\zeta^2+z^2+w^2|\Bigr\}.
\end{array}\label{est}
\end{equation} 
The group $G(E_{s,t})=\hbox{Aut}(E_{s,t})$ is given by
\begin{equation}
\left(
\begin{array}{c}
\zeta\\
z\\
w
\end{array}
\right)
\mapsto A
\left(
\begin{array}{c}
\zeta\\
z\\
w
\end{array}
\right),\label{autmu}
\end{equation}
where $A\in SO_3(\RR)$. The orbits of the action of the group $G(E_{s,t})$ on $E_{s,t}$ are the following pairwise $CR$ non-equivalent strongly pseudoconvex hypersurfaces:
\begin{equation}
\begin{array}{ll}
\mu_{\alpha}:=&\left\{(\zeta:z:w)\in\CC\PP^2: |\zeta|^2+|z|^2+|w|^2=\alpha |\zeta^2+z^2+w^2|\right\},\\
&\hspace{8cm} s<\alpha<t.
\end{array}\label{mu}
\end{equation}
We denote the group of all maps of the form (\ref{autmu}) by ${\cal R}_{\mu}$.
\vspace{0.3cm}

\noindent {\bf (b)} For $1<t<\infty$ define
\begin{equation}
E_t:=\Bigl\{(\zeta:z:w)\in\CC\PP^2: |\zeta|^2+|z|^2+|w|^2< t |\zeta^2+z^2+w^2|\Bigr\}.\label{et}
\end{equation}
The group $G(E_t)$ coincides with ${\cal R}_{\mu}$, and its action on $E_t$ has, apart from strongly pseudoconvex orbits, the following totally real orbit:
\begin{equation}
{\cal O}_3:=\RR\PP^2\subset\CC\PP^2.\label{calo3}
\end{equation}
\vspace{0.5cm}

\noindent{\bf (7)} Here we explicitly describe all covers of the domains $E_{s,t}$ and hypersurfaces $\mu_{\alpha}$ introduced in {\bf (6)} (for more details see \cite{I3}). As we will see below, to one of the covers of $E_{1,t}$ a totally real orbit can be attached.
\vspace{0.3cm} 

\noindent {\bf (a)} Let ${\cal Q}_{+}$ be the variety in $\CC^3$ given by
\begin{equation}
z_1^2+z_2^2+z_3^2=1.\label{quadricplus}
\end{equation}

Consider the map $\Phi_{\mu}: \CC^2\setminus\{0\}\ra {\cal Q}_{+}$ defined by the formulas
$$
\begin{array}{lll}
z_1& = &\displaystyle -i(z^2+w^2)+i\frac{z\overline{w}-w\overline{z}}{|z|^2+|w|^2},\\
\vspace{0cm}&&\\
z_2 & = & \displaystyle z^2-w^2-\frac{z\overline{w}+w\overline{z}}{|z|^2+|w|^2},\\
\vspace{0cm}&&\\
z_3 & = & \displaystyle 2zw+\frac{|z|^2-|w|^2}{|z|^2+|w|^2}.
\end{array}
$$
This map was introduced in \cite{R}. It is straightforward to verify that $\Phi_{\mu}$ is a 2-to-1 covering map onto ${\cal Q}_{+}\setminus\RR^3$. We now equip the domain of $\Phi_{\mu}$ with the pull-back complex structure under $\Phi_{\mu}$ and denote the resulting complex manifold by $M_{\mu}^{(4)}$.

For $1\le s<t<\infty$ define
\begin{equation}
\begin{array}{lll}
E_{s,t}^{(2)}&:=&\Bigl\{(z_1,z_2,z_3)\in\CC^3: s<|z_1|^2+|z_2|^2+|z_3|^2<t \Bigr\}\cap {\cal Q}_{+},\\
E_{s,t}^{(4)}&:=&\Bigl\{(z,w)\in M_{\mu}^{(4)}: \sqrt{(s-1)/2}<|z|^2+|w|^2<\\
&&\hspace{7.4cm}\sqrt{(t-1)/2}\Bigr\}.
\end{array}\label{est24}
\end{equation}
These domains are respectively a 2- and 4-sheeted cover of the domain $E_{s,t}$, where $E_{s,t}^{(2)}$ covers $E_{s,t}$ by means of the map $\Psi_{\mu}:\,(z_1,z_2,z_3)\mapsto (z_1:z_2:z_3)$ and $E_{s,t}^{(4)}$ covers $E_{s,t}$ by means of the composition $\Psi_{\mu}\circ\Phi_{\mu}$.   

The group $G\left(E_{s,t}^{(2)}\right)$ consists of all maps
\begin{equation}
\left(
\begin{array}{c}
z_1\\
z_2\\
z_3
\end{array}
\right)
\mapsto A
\left(
\begin{array}{c}
z_1\\
z_2\\
z_3
\end{array}\label{autmu2}
\right),
\end{equation}
where $A\in SO_3(\RR)$. The $G\left(E_{s,t}^{(2)}\right)$-orbits are the following pairwise $CR$ non-equivalent strongly pseudoconvex hypersurfaces:  
\begin{equation}
\mu_{\alpha}^{(2)}:=\left\{(z_1,z_2,z_3)\in\CC^3: |z_1|^2+|z_2|^2+|z_3|^2=\alpha \right\}\cap {\cal Q}_{+},\,s<\alpha<t\label{mu2}
\end{equation}
(note that $\mu_{\alpha}^{(2)}$ is a 2-sheeted cover of $\mu_{\alpha}$). We denote the group of all maps of the form (\ref{autmu2}) by ${\cal R}_{\mu}^{(2)}$. This group is clearly isomorphic to ${\cal R}_{\mu}$. 

The group $G\left(E_{s,t}^{(4)}\right)$ consists of all maps
\begin{equation}
\left(
\begin{array}{c}
z\\
w
\end{array}
\right) \mapsto
A
\left(
\begin{array}{c}
z\\
w
\end{array}
\right),\label{autmu41}
\end{equation}
where $A\in SU_2$. The $G\left(E_{s,t}^{(4)}\right)$-orbits are the following pairwise $CR$ non-equivalent strongly pseudoconvex hypersurfaces:
\begin{equation}
\mu_{\alpha}^{(4)}:=\left\{(z,w)\in M_{\mu}^{(4)}: |z|^2+|w|^2=\sqrt{(\alpha-1)/2}\right\},\quad s<\alpha<t\label{mu4}
\end{equation}
(note that $\mu_{\alpha}^{(4)}$ is a 4-sheeted cover of $\mu_{\alpha}$).
\vspace{0.3cm}

\noindent {\bf (b)} For $1<t<\infty$ define
\begin{equation}
E_t^{(2)}:=\left\{(z_1,z_2,z_3)\in\CC^3: |z_1|^2+|z_2|^2+|z_3|^2<t \right\}\cap {\cal Q}_{+}.\label{et2}
\end{equation}
The group $G\left(E_t^{(2)}\right)$ coincides with ${\cal R}_{\mu}^{(2)}$, and, apart from strongly pseudoconvex orbits, its action on $E_t^{(2)}$ has the totally real orbit
\begin{equation}
{\cal O}_4:=\left\{(x_1,x_2,x_3)\in\RR^3: x_1^2+x_2^2+x_3^2=1\right\}={\cal Q}_{+}\cap\RR^3.\label{calo4}
\end{equation}  
\vspace{0.5cm}

\noindent{\bf (8)} In this example strongly pseudoconvex orbits and a totally real orbit arise.
\vspace{0.3cm}

\noindent {\bf (a)} For $-1\le s<t\le 1$ define
\begin{equation}
\begin{array}{lll}
\Omega_{s,t}&:=&\Bigl\{(z,w)\in\CC^2: s |z^2+w^2-1|<|z|^2+|w|^2-1<\\
&& \hspace{7cm}t |z^2+w^2-1|\Bigr\}.
\end{array}\label{omegast}
\end{equation}
The group $G(\Omega_{s,t})$ consists of all maps
\begin{equation}
\left(
\begin{array}{c}
z\\
w
\end{array}
\right) \mapsto \displaystyle\frac{\left(\begin{array}{cc}
a_{11} & a_{12}\\
a_{21} & a_{22}
\end{array}
\right)\left(
\begin{array}{c}
z\\
w
\end{array}\right)+\left(
\begin{array}{c}
b_1\\
b_2
\end{array}
\right)}{c_1z+c_2w+d},\label{autnu}
\end{equation}
where
\begin{equation}
Q:=\left(\begin{array}{ccc}
a_{11} & a_{12} & b_1\\
a_{21} & a_{22} & b_2\\
c_1& c_2 &d
\end{array}\label{matq}
\right)
\in SO_{2,1}(\RR)^0.
\end{equation}
The orbits of $G(\Omega_{s,t})$ on $\Omega_{s,t}$ are the following pairwise $CR$ non-equivalent strongly pseudoconvex hypersurfaces:
\begin{equation}
\begin{array}{ll}
\nu_{\alpha}:=&\left\{(z,w)\in\CC^2: |z|^2+|w|^2-1=\alpha|z^2+w^2-1|\right\}\setminus\\
&\left\{(x,u)\in\RR^2:x^2+u^2=1\right\},\quad s<\alpha<t.
\end{array}\label{nu}
\end{equation}
We denote the group of all maps of the form (\ref{autnu}) by ${\cal R}_{\nu}$.
\vspace{0.3cm}

\noindent {\bf (b)} For $-1<t\le1$ define
\begin{equation}
\Omega_t:=\left\{(z,w)\in\CC^2: |z|^2+|w|^2-1< t |z^2+w^2-1|\right\}.\label{omegat}
\end{equation}
The group $G(\Omega_t)$ for $t<1$ coincides with ${\cal R}_{\nu}$, and its action on $\Omega_t$, apart from strongly pseudoconvex orbits, has the totally real orbit
\begin{equation}
{\cal O}_5:=\left\{(x,u)\in\RR^2: x^2+u^2<1\right\}\subset\CC^2.\label{calo5}
\end{equation}

We note that $\Omega_1$ is holomorphically equivalent to $\Delta^2$ (see {\bf (11)(c)} below); hence it has a 6-dimensional automorphism group and therefore will be excluded from our considerations.
\vspace{0.5cm}

\noindent{\bf (9)} In this example strongly pseudoconvex orbits and a complex curve orbit occur.
\vspace{0.3cm}

\noindent {\bf (a)} For $1\le s<t\le\infty$ define
\begin{equation}
\begin{array}{lll}
D_{s,t}&:=&\Bigl\{(z,w)\in\CC^2:  s |1+z^2-w^2|<1+|z|^2-|w|^2<\\
&&\hspace{3.5cm} t |1+z^2-w^2|,\,\hbox{Im}\left(z(1+\overline{w})\right)>0\Bigr\},
\end{array}\label{dst}
\end{equation}
where $D_{s,\infty}$ is assumed not to include the complex curve 
\begin{equation}
{\cal O}:=\left\{(z,w)\in{\Bbb C}^2: 1+z^2-w^2=0,\, \hbox{Im}(z(1+\overline{w}))>0\right\}.\label{calo}
\end{equation}
For every matrix $Q\in SO_{2,1}(\RR)^0$ as in (\ref{matq}) consider the map
\begin{equation}
\left(
\begin{array}{c}
z\\
w
\end{array}
\right) \mapsto \displaystyle\frac{\left(\begin{array}{cc}
a_{22} & b_2\\
c_2 & d
\end{array}
\right)\left(
\begin{array}{c}
z\\
w
\end{array}\right)+\left(
\begin{array}{c}
a_{21}\\
c_1
\end{array}
\right)}{a_{12}z+b_1w+a_{11}}.\label{auteta}
\end{equation}
The group $G(D_{s,t})=\hbox{Aut}(D_{s,t})$ consists of all such maps.  
The orbits of $G(D_{s,t})$ on $D_{s,t}$ are the following pairwise $CR$ non-equivalent strongly pseudoconvex hypersurfaces:
\begin{equation}
\begin{array}{ll}
\eta_{\alpha}:=&\Bigl\{(z,w)\in\CC^2: 1+|z|^2-|w|^2=\alpha|1+z^2-w^2|,\\
&\hbox{Im}(z(1+\overline{w}))>0\Bigr\},\quad s<\alpha<t.
\end{array}\label{eta}
\end{equation}
We denote the group of all maps of the form (\ref{auteta}) by ${\cal R}_{\eta}$ (note that ${\cal R}_{\eta}$ is isomorphic to ${\cal R}_{\nu}$).
\vspace{0.3cm}

\noindent {\bf (b)} For $1\le s<\infty$ define
\begin{equation}
\begin{array}{lll}
D_s&:=&\Bigl\{(z,w)\in\CC^2:  1+|z|^2-|w|^2>s |1+z^2-w^2|,\\
&&\hspace{5cm}\hbox{Im}\left(z(1+\overline{w})\right)>0\Bigr\},
\end{array}\label{ds}
\end{equation}
The group $G(D_s)=\hbox{Aut}(D_s)$ coincides with ${\cal R}_{\eta}$. Apart from strongly pseudoconvex orbits, its action on $D_s$ has the complex curve orbit ${\cal O}$.
\vspace{0.5cm}

\noindent{\bf (10)} In this example we explicitly describe all covers of the domains $\Omega_{s,t}$, $D_{s,t}$ and the hypersurfaces $\nu_{\alpha}$, $\eta_{\alpha}$ introduced in {\bf (8)} and {\bf (9)} (for more details see \cite{I3}). Only strongly pseudoconvex orbits occur here. 

Denote by $(z_0:z_1:z_2:z_3)$ homogeneous coordinates in $\CC\PP^3$; we think of the hypersurface $\{z_0=0\}$
as the infinity. Let ${\cal Q}_{-}$ be the variety in $\CC\PP^3$ given by
\begin{equation}
z_1^2+z_2^2-z_3^2=z_0^2.\label{calqminus}
\end{equation}
Next, let $(\zeta:z:w)$ be homogeneous coordinates in $\CC\PP^2$ (where we think of the hypersurface $\{\zeta=0\}$ as the infinity), and let
\begin{equation}
\Sigma:=\left\{(\zeta:z:w)\in\CC\PP^2: |w|<|z|\right\}.\label{Sigma}
\end{equation}
For every integer $n\ge 2$ consider the map $\Phi^{(n)}$ from $\Sigma$ to ${\cal Q}_{-}$ defined as follows:
\begin{equation}
\begin{array}{lll}
z_0&=&\zeta^n,\\
\vspace{0cm}&&\\
z_1 & = &\displaystyle -i(z^n+z^{n-2}w^2)-i\frac{z\overline{w}+w\overline{z}}{|z|^2-|w|^2}\zeta^n,\\
\vspace{0cm}&&\\
z_2 & = & \displaystyle z^n-z^{n-2}w^2+\frac{z\overline{w}-w\overline{z}}{|z|^2-|w|^2}\zeta^n,\\
\vspace{0cm}&&\\
z_3 & = &\displaystyle -2iz^{n-1}w-i\frac{|z|^2+|w|^2}{|z|^2-|w|^2}\zeta^n.
\end{array}\label{bigphin}
\end{equation}
The above maps were introduced in \cite{I3} and are analogous to the map $\Phi_{\mu}$ defined in {\bf (7)}. Further, set 
\begin{equation}
\begin{array}{lll}
{\cal A}_{\nu}^{(n)}&:=&\left\{(z,w)\in\CC^2: 0<|z|^n-|z|^{n-2}|w|^2<1\right\},\\
{\cal A}_{\eta}^{(n)}&:=&\left\{(z,w)\in\CC^2: |z|^n-|z|^{n-2}|w|^2>1\right\}
\end{array}\label{anuetan}
\end{equation}
(both domains lie in the finite part of $\CC\PP^2$ given by $\zeta=1$). Clearly, ${\cal A}_{\nu}^{(n)},\,{\cal A}_{\eta}^{(n)}\subset\Sigma$ for all $n\ge 2$. Let $\Phi_{\nu}^{(n)}$ and $\Phi_{\eta}^{(n)}$ be the restrictions of $\Phi^{(n)}$ to  ${\cal A}_{\nu}^{(n)}$ and ${\cal A}_{\eta}^{(n)}$, respectively. It is straightforward to observe that $\Phi_{\nu}^{(n)}$ and $\Phi_{\eta}^{(n)}$ are $n$-to-1 covering maps onto
\begin{equation}
\begin{array}{ll}
{\cal A}_{\nu}:=\Bigl\{(z_1,z_2,z_3)\in\CC^3:& -1<|z_1|^2+|z_2|^2-|z_3|^2<1,\\
&\hspace{4cm}\hbox{Im}\,z_3<0\Bigr\}\cap {\cal Q}_{-}
\end{array}\label{calanu}
\end{equation}
and
\begin{equation}
\begin{array}{ll}
{\cal A}_{\eta}:=\Bigl\{(z_1,z_2,z_3)\in\CC^3:&|z_1|^2+|z_2|^2-|z_3|^2>1,\\
&\hspace{1.5cm}\hbox{Im}(z_2(\overline{z_1}+\overline{z_3}))>0\Bigr\}\cap {\cal Q}_{-} ,
\end{array}\label{calaeta}
\end{equation}
respectively (both domains lie in the finite part of $\CC\PP^3$ given by $z_0=1$). We now introduce on ${\cal A}_{\nu}^{(n)}$, ${\cal A}_{\eta}^{(n)}$ the pull-back complex structures under the maps $\Phi_{\nu}^{(n)}$, $\Phi_{\eta}^{(n)}$, respectively, and denote the resulting complex manifolds by $M_{\nu}^{(n)}$, $M_{\eta}^{(n)}$.

Further, let $\Lambda:\CC\times\Delta\to\Sigma\cap\{\zeta=1\}$ be the following covering map:
\begin{equation}
\begin{array}{lll}
z & \mapsto & e^z,\\
w & \mapsto & e^z w,
\end{array}\label{maplambda}
\end{equation}
where $z\in\CC$, $w\in\Delta$. Define
$$
\begin{array}{lll}
U_{\nu} & := & \left\{(z,w)\in\CC^2:|w|<1,\,\exp(2\hbox{Re}\,z)(1-|w|^2)<1\right\},\\
U_{\eta} & := & \left\{(z,w)\in\CC^2:|w|<1,\,\exp(2\hbox{Re}\,z)(1-|w|^2)>1\right\}.
\end{array}
$$
Denote by $\Lambda_{\nu}$, $\Lambda_{\eta}$ the restrictions of $\Lambda$ to $U_{\nu}$, $U_{\eta}$, respectively. Clearly, $U_{\nu}$ covers $M_{\nu}^{(2)}$ by means of $\Lambda_{\nu}$, and $U_{\eta}$ covers $M_{\eta}^{(2)}$ by means of $\Lambda_{\eta}$. We now introduce on $U_{\nu}$, $U_{\eta}$ the pull-back complex structures under the maps $\Lambda_{\nu}$, $\Lambda_{\eta}$, respectively, and denote the resulting complex manifolds by $M_{\nu}^{(\infty)}$, $M_{\eta}^{(\infty)}$.

For $-1\le s<t\le 1$, $n\ge 2$ we now define
\begin{equation}
\begin{array}{lll}
\Omega_{s,t}^{(n)}&:=&\Bigl\{(z,w)\in M_{\nu}^{(n)}:\sqrt{(s+1)/2}<|z|^n-|z|^{n-2}|w|^2<\\
&&\hspace{6.4cm}\sqrt{(t+1)/2}\Bigr\},\\
\Omega_{s,t}^{(\infty)}&:=&\Bigl\{(z,w)\in M_{\nu}^{(\infty)}:\sqrt{(s+1)/2}<\\
&&\hspace{1.4cm}\exp\left(2\hbox{Re}\,z\right)(1-|w|^2)<
\sqrt{(t+1)/2}\Bigr\}.
\end{array}\label{omegastinftyn}
\end{equation}
The domain $\Omega_{s,t}^{(n)}$, $n\ge 2$, is an $n$-sheeted cover of the domain $\Omega_{s,t}$ introduced in {\bf (8)} and the domain $\Omega_{s,t}^{(\infty)}$ is its infinitely-sheeted cover. The domain $\Omega_{s,t}^{(n)}$ covers $\Omega_{s,t}$ by means of the composition $\Psi_{\nu}\circ\Phi_{\nu}^{(n)}$, where $\Psi_{\nu}$ is the following 1-to-1 map from ${\cal A}_{\nu}$ to $\CC^2$ : $(z_1,z_2,z_3)\mapsto (z_1/z_3,z_2/z_3)$; the domain $\Omega_{s,t}^{(\infty)}$ covers $\Omega_{s,t}$ by means of the composition $\Psi_{\nu}\circ\Phi_{\nu}^{(2)}\circ\Lambda_{\nu}$.    

The group $G\left(\Omega_{s,t}^{(n)}\right)$ consists of all maps of the form
\begin{equation}
\begin{array}{lll}
z & \mapsto & \displaystyle z\,\sqrt[n]{\left(a+b\, w/z\right)^2}, \\
\vspace{0cm}&&\\ 
w & \mapsto & \displaystyle z\,\frac{\overline{b}+\overline{a}\,w/z}{a+b\,w/z}\sqrt[n]{\left(a+b\,w/z\right)^2},
\end{array}\label{autnun1}
\end{equation}
where $|a|^2-|b|^2=1$. The $G\left(\Omega_{s,t}^{(n)}\right)$-orbits are the following pairwise $CR$ non-equivalent strongly pseudoconvex hypersurfaces:
\begin{equation}
\begin{array}{ll}
\nu_{\alpha}^{(n)}:=&\left\{(z,w)\in M_{\nu}^{(n)}:|z|^n-|z|^{n-2}|w|^2=\sqrt{(\alpha+1)/2}\right\},\\
&\hspace{7.3cm}s<\alpha<t
\end{array}\label{nun}
\end{equation}
(note that $\nu_{\alpha}^{(n)}$ is an $n$-sheeted cover of $\nu_{\alpha}$). We denote the group of all maps of the form (\ref{autnun1}) by ${\cal R}^{(n)}$.

The group $G\left(\Omega_{s,t}^{(\infty)}\right)$ consists of all maps of the form
\begin{equation}
\begin{array}{lll}
z & \mapsto & \displaystyle z+\ln(a+bw),\\
\vspace{0cm}&&\\
w & \mapsto & \displaystyle\frac{\overline{b}+\overline{a}w}{a+bw},
\end{array}\label{autnuinfty1}
\end{equation}
where $|a|^2-|b|^2=1$. The $G\left(\Omega_{s,t}^{(\infty)}\right)$-orbits are the following pairwise $CR$ non-equivalent strongly pseudoconvex hypersurfaces:
\begin{equation}
\begin{array}{ll}
\nu_{\alpha}^{(\infty)}:=&\left\{(z,w)\in M_{\nu}^{(\infty)}:\exp\left(2\hbox{Re}\,z\right)(1-|w|^2)=\sqrt{(\alpha+1)/2}\right\},\\ &\hspace{8.2cm}s<\alpha<t
\end{array}\label{nuinfty}
\end{equation}
(note that $\nu_{\alpha}^{(\infty)}$ is an infinitely-sheeted cover of $\nu_{\alpha}$). We denote the group of all maps of the form (\ref{autnuinfty1}) by ${\cal R}^{(\infty)}$. 

Next, for $1\le s<t\le \infty$, $n\ge 2$ we define
\begin{equation}
\begin{array}{lll}
D_{s,t}^{(2)}& := & \Bigl\{(z_1,z_2,z_3)\in\CC^3:s<|z_1|^2+|z_2|^2-|z_3|^2<t,\\
&&\hspace{4.5cm}\hbox{Im}(z_2(\overline{z_1}+\overline{z_3}))>0\Bigr\}\cap {\cal Q}_{-},\\
D_{s,t}^{(2n)}&:=&\Bigl\{(z,w)\in M_{\eta}^{(n)}:\sqrt{(s+1)/2}<|z|^n-|z|^{n-2}|w|^2<\\
&&\hspace{7cm}\sqrt{(t+1)/2}\Bigr\},\\
D_{s,t}^{(\infty)}&:=&\Bigl\{(z,w)\in M_{\eta}^{(\infty)}:\sqrt{(s+1)/2}<\\
&&\hspace{1.4cm}\exp\left(2\hbox{Re}\,z\right)(1-|w|^2)<
\sqrt{(t+1)/2}\Bigr\}.
\end{array}\label{dstinftyn}
\end{equation}
The domain $D_{s,t}^{(2n)}$, $n\ge 1$, is a $2n$-sheeted cover of the domain $D_{s,t}$ introduced in {\bf (9)} and the domain $D_{s,t}^{(\infty)}$ is its infinitely-sheeted cover. The domain $D_{s,t}^{(2)}$ covers $D_{s,t}$ by means of the map $\Psi_{\eta}$, which is the following 2-to-1 map from ${\cal A}_{\eta}$ to $\CC^2$ : $(z_1,z_2,z_3)\mapsto (z_2/z_1,z_3/z_1)$; the domain $D_{s,t}^{(2n)}$ for $n\ge 2$ covers $D_{s,t}$ by means of the composition $\Psi_{\eta}\circ\Phi_{\eta}^{(n)}$; the domain $D_{s,t}^{(\infty)}$ covers $D_{s,t}$ by means of the composition $\Psi_{\eta}\circ\Phi_{\eta}^{(2)}\circ\Lambda_{\eta}$.  

To obtain an $n$-sheeted cover of $D_{s,t}$ for odd $n\ge 3$, the domain $D_{s,t}^{(4n)}$ must be factored by the action of the cyclic group of four elements generated by the following automorphism of $M_{\eta}^{(2n)}$:
\begin{equation}
\begin{array}{lll}
z & \mapsto & \displaystyle iz^2\overline{z}\sqrt[n]{\frac{1-|w|^2/|z|^2+z^{-2n}\overline{w}/\overline{z}}{\sqrt{|z|^{4n}(1-|w|^2/|z|^2)^2-1}}},\\
\vspace{0cm}&&\\
w & \mapsto & i\displaystyle\frac{1+z^{2n-1}w(1-|w|^2/|z|^2)}{\overline{w}/\overline{z}+z^{2n}(1-|w|^2/|z|^2)}\times\\
\vspace{0cm}&&\\
&&\displaystyle z^2\overline{z}\sqrt[n]{\frac{1-|w|^2/|z|^2+z^{-2n}\overline{w}/\overline{z}}{\sqrt{|z|^{4n}(1-|w|^2/|z|^2)^2-1}}}.
\end{array}\label{autorderfour}
\end{equation}
Let $\Pi^{(n)}$ denote the corresponding factorization map and $M_{\eta}^{'(n)}:=\Pi^{(n)}\left(M_{\eta}^{(2n)}\right)$. Then $D_{s,t}^{(n)}:=\Pi^{(n)}\left(D_{s,t}^{(4n)}\right)$ is an $n$-sheeted cover of $D_{s,t}$.

The group $G\left(D_{s,t}^{(2)}\right)$ consists of all maps of the form (\ref{autmu2}) with $A\in SO_{2,1}(\RR)^0$. We denote this group by ${\cal R}^{(1)}$ (observe that ${\cal R}^{(1)}$ is isomorphic to ${\cal R}_{\eta}$ -- see  (\ref{auteta})). The $G\left(D_{s,t}^{(2)}\right)$-orbits are the following pairwise $CR$ non-equivalent strongly pseudoconvex hypersurfaces:
\begin{equation}
\begin{array}{lll}
\eta_{\alpha}^{(2)}&:=&\left\{(z_1,z_2,z_3)\in{\cal A}_{\eta}: |z_1|^2+|z_2|^2-|z_3|^2=\sqrt{(\alpha+1)/2},\right\},\\
&&\hspace{8.2cm}s<\alpha<t
\end{array}\label{eta2}
\end{equation}
(note that $\eta_{\alpha}^{(2)}$ is a 2-sheeted cover of $\eta_{\alpha}$).
 For $n\ge 2$ the group $G\left(D_{s,t}^{(2n)}\right)$ coincides with ${\cal R}^{(n)}$ (see (\ref{autnun1})), where we think of elements of ${\cal R}^{(n)}$ as maps defined on $D_{s,t}^{(2n)}$ rather than on $\Omega_{s,t}^{(n)}$. The $G\left(D_{s,t}^{(2n)}\right)$-orbits are the following pairwise $CR$ non-equivalent strongly pseudoconvex hypersurfaces:
\begin{equation}
\begin{array}{lll}
\eta_{\alpha}^{(2n)}&:=&\left\{(z,w)\in M_{\eta}^{(n)}:|z|^n-|z|^{n-2}|w|^2=\sqrt{(\alpha+1)/2}\right\},\\
&&\hspace{6cm}s<\alpha<t, \,n\ge 2
\end{array}\label{eta2n}
\end{equation}
(note that $\eta_{\alpha}^{(2n)}$ is a $2n$-sheeted cover of $\eta_{\alpha}$). 

Next, the group $G\left(D_{s,t}^{(n)}\right)$ for odd $n\ge 3$ consists of all lifts from the domain $D_{1,\infty}$ to $D^{(n)}_{1,\infty}=M_{\eta}^{'(n)}$ of all elements of ${\cal R}_{\eta}$ (see (\ref{auteta})). This group is isomorphic to ${\cal R}^{(n)}$. Note, however, that the isotropy subgroup of every point under the action of this group on $D_{s,t}^{(n)}$ consists of two points, whereas the action of ${\cal R}^{(n)}$ on $D_{s,t}^{(2n)}$ is free (observe also that the isotropy subgroup of every point under the action of ${\cal R}_{\eta}$ consists of two points and that the action of ${\cal R}^{(1)}$ on $D_{s,t}^{(2)}$ is free). This difference will be important in the proof of Theorem \ref{classrealhypersurfaces} (see step (II) of the orbit gluing procedure there).

The $G\left(D_{s,t}^{(n)}\right)$-orbits are the following pairwise $CR$ non-equivalent strongly pseudoconvex hypersurfaces:
\begin{equation}
\eta_{\alpha}^{(n)}:=\Pi_n\left(\eta_{\alpha}^{(4n)}\right),\quad s<\alpha<t\label{etan}
\end{equation}
(note that $\eta_{\alpha}^{(n)}$ is an $n$-sheeted cover of $\eta_{\alpha}$). 

Finally, the group $G\left(D_{s,t}^{(\infty)}\right)$ coincides with ${\cal R}^{(\infty)}$ (see (\ref{autnuinfty1})), where we think of the elements of ${\cal R}^{(\infty)}$ as maps defined on $D_{s,t}^{(\infty)}$ rather than on $\Omega_{s,t}^{(\infty)}$. The $G\left(D_{s,t}^{(\infty)}\right)$-orbits are the following pairwise $CR$ non-equivalent hypersurfaces:
\begin{equation}
\begin{array}{ll}
\eta_{\alpha}^{(\infty)}:=&\left\{(z,w)\in M_{\eta}^{(\infty)}:\exp\left(2\hbox{Re}\,z\right)(1-|w|^2)=\sqrt{(\alpha+1)/2}\right\},\\
&\hspace{8.2cm}s<\alpha<t
\end{array}\label{etainfty}
\end{equation}
(here $\eta_{\alpha}^{(\infty)}$ is an infinitely-sheeted cover of $\eta_{\alpha}$).
\vspace{0.5cm}

\noindent{\bf (11)} Here we show how a Levi-flat and complex curve orbit can be attached to some of the domains introduced in {\bf (8)} and {\bf (10)}. 
\vspace{0.3cm}

\noindent{\bf (a)} It is straightforward to show from the explicit from of $\Phi^{(n)}$, for $n\ge 2$ (see (\ref{bigphin})), that the complex structure of $M_{\eta}^{(n)}$ extends to a complex structure on
$$
\tilde{\cal A}_{\eta}^{(n)}:=\left\{(\zeta:z:w)\in\CC\PP^2: |z|^n-|z|^{n-2}|w|^2>|\zeta|^n\right\}.
$$
The set at infinity in $\tilde{\cal A}_{\eta}^{(n)}$ is
\begin{equation}
{\cal O}^{(2n)}:=\left\{(0:z:w)\in\CC\PP^2: |w|<|z|\right\},\label{cal02n}
\end{equation}
and we have $\tilde{\cal A}_{\eta}^{(n)}={\cal A}_{\eta}^{(n)}\cup{\cal O}^{(2n)}$ (see (\ref{anuetan})). Let $\tilde M_{\eta}^{(n)}$ denote $\tilde{\cal A}_{\eta}^{(n)}$ with the extended complex structure.  In the complex structure of $\tilde M_{\eta}^{(n)}$ the set ${\cal O}^{(2n)}$ is a complex curve whose complex structure is identical to that induced from $\CC\PP^2$. The action of the group ${\cal R}^{(n)}$ (see (\ref{autnun1})) extends to an action by holomorphic transformations on $\tilde M_{\eta}^{(n)}$, and ${\cal O}^{(2n)}$ is an orbit of this action. The map $\Phi^{(n)}$ has ramification locus on ${\cal O}^{(2n)}$ and maps it in a 1-to-1 fashion onto the complex curve
\begin{equation}
\begin{array}{lll}
{\cal O}^{(2)}&:=&\Bigl\{(0:z_1:z_2:z_3)\in\CC\PP^3:z_1^2+z_2^2-z_3^2=0,\\
&&\hspace{5cm}\hbox{Im}(z_2(\overline{z_1}+\overline{z_3}))>0
\Bigr\}.
\end{array}\label{calo(2)}    
\end{equation}
Note that ${\cal O}^{(2)}$ is an ${\cal R}^{(1)}$-orbit (clearly, ${\cal R}^{(1)}$ acts on all of ${\cal Q}_{-}$ --  see {\bf (10)}).

For $1\le s<\infty$ and all $n\ge 1$ define
\begin{equation}
D_s^{(2n)}:=D_{s,\infty}^{(2n)}\cup{\cal O}^{(2n)}.\label{ds2n}
\end{equation}
The group $G\Bigl(D_s^{(2n)}\Bigr)$ (with the exception of the case $n=1$, $s=1$) coincides with ${\cal R}^{(n)}$ for all $n$; its orbits in $D_s^{(2n)}$ are the strongly pseudoconvex hypersurfaces $\eta_{\alpha}^{(2n)}$ for $\alpha>s$ (see (\ref{eta2n})) and the complex curve ${\cal O}^{(2n)}$. The map $\Psi_{\eta}$ is a branched covering map on $D_s^{(2)}$, has ramification locus on ${\cal O}^{(2)}$, maps it in a 1-to-1 fashion onto the complex curve ${\cal O}\subset\CC^2$ (see (\ref{calo})), and takes $D_s^{(2)}$ onto $D_s$. Similarly, for $n\ge 2$, the map $\Psi_{\eta}\circ\Phi_{\eta}^{(n)}$ is a branched covering map on $D_s^{(2n)}$, has ramification locus on ${\cal O}^{(2n)}$ and takes $D_s^{(2n)}$ onto $D_s$.

We note that $D_1^{(2)}={\cal A}_{\eta}\cup{\cal O}^{(2)}$ (see (\ref{calaeta})) is holomorphically equivalent to $\Delta^2$ (see {\bf (11)(c)} below), hence it will be excluded from our considerations.
\vspace{0.3cm}   

\noindent{\bf (b)} Fix an odd $n\in\NN$, $n\ge 3$, and let $\Gamma^{(n)}$ be the cyclic group of four elements generated by the obvious extension of automorphism (\ref{autorderfour}) to $\tilde M_{\eta}^{(2n)}=D_1^{(4n)}$. The group $\Gamma^{(n)}$ acts freely properly discontinuously on $M_{\eta}^{(2n)}\subset \tilde M_{\eta}^{(2n)}$ and fixes every point in ${\cal O}^{(4n)}$. It is straightforward to show that the orbifold obtained by factoring $\tilde M_{\eta}^{(2n)}$ by the action of $\Gamma^{(n)}$ can in fact be given the structure of a complex manifold (we denote it by $\tilde M_{\eta}^{'(n)}$) that extends the structure of $M_{\eta}^{'(n)}$ (see {\bf (10)}). The extension of the map $\Pi^{(n)}$ (see {\bf (10)}) is holomorphic on all of $\tilde M_{\eta}^{(2n)}$, has ramification locus on ${\cal O}^{(4n)}$  and maps ${\cal O}^{(4n)}$ onto a complex curve ${\cal O}^{(n)}\subset \tilde M_{\eta}^{'(n)}$ in a 1-to-1 fashion (note that $\tilde M_{\eta}^{'(n)}=M_{\eta}^{'(n)}\cup{\cal O}^{(n)}$). The covering map from $M_{\eta}^{'(n)}$ onto $D_{1,\infty}$ (see (\ref{dst})) extends to a branched covering map from $\tilde M_{\eta}^{'(n)}$ onto $D_1$ (see (\ref{ds})) with ramification locus ${\cal O}^{(n)}$, and takes ${\cal O}^{(n)}$ onto ${\cal O}$ (see (\ref{calo})) in a 1-to-1 fashion.

For $1\le s<\infty$ define
\begin{equation}
D_s^{(n)}:=\Pi^{(n)}\left(D_s^{(4n)}\right).\label{dsn}
\end{equation}
The group $G\left(D_s^{(n)}\right)$ is isomorphic to ${\cal R}^{(n)}$ and consists of the extensions from $D_{s,\infty}^{(n)}=D_s^{(n)}\setminus{\cal O}^{(n)}$ to $D_s^{(n)}$ of all elements of the group $G\left(D_{s,\infty}^{(n)}\right)$. The $G\left(D_s^{(n)}\right)$-orbits are the strongly pseudoconvex hypersurfaces $\eta_{\alpha}^{(n)}$ with $\alpha>s$ (see (\ref{etan})) and the complex curve ${\cal O}^{(n)}$.
\vspace{0.2cm}

\noindent{\bf (c)} Define
\begin{equation}
M^{(1)}:=\Psi_{\nu}^{-1}(\Omega_1)\cup D_1^{(2)}\cup {\cal O}_0^{(1)},\label{m1}
\end{equation}
where
\begin{equation}
\begin{array}{lll}
{\cal O}_0^{(1)}&:=&\Bigl\{(z_1,z_2,z_3)\in\CC^3\setminus{\Bbb R}^3:|iz_1+z_2|=|iz_3-1|,\\
&&\hspace{2cm}|iz_1-z_2|=|iz_3+1|,\,\hbox{Im}\,z_3<0\Bigr\}\cap {\cal Q}_{-}
\end{array}\label{calo61}
\end{equation}
(see {\bf (10)} for the definition of $\Psi_{\nu}$ and (\ref{calqminus}) for the definition of ${\cal Q}_{-}$). Clearly, $M^{(1)}$ is invariant under the action of the group ${\cal R}^{(1)}$ (defined in {\bf (10)}) on ${\cal Q}_{-}$. We will now describe the orbits of the ${\cal R}^{(1)}$-action on $M^{(1)}$. The hypersurfaces $\eta_{\alpha}^{(2)}$ for $\alpha>1$ (see (\ref{eta2})) and
$$
\nu_{\alpha}^{(1)}:=\left\{(z_1,z_2,z_3)\in\CC^3: |z_1|^2+|z_2|^2-|z_3|^2=\alpha\right\}\cap{\cal Q}_{-}
$$
for $-1<\alpha<1$ are strongly pseudoconvex orbits (note that $\nu_{\alpha}^{(1)}$ is equivalent to $\nu_{\alpha}$ (see (\ref{nu})) by means of the map $\Psi_{\nu}$); the hypersurface ${\cal O}_0^{(1)}$ is the unique Levi-flat orbit; the surfaces
\begin{equation}
{\cal O}_6 :=\Bigl\{(z_1,z_2,z_3)\in i{\Bbb R}^3:\hbox{Im}\,z_3<0\Bigr\}\cap {\cal Q}_{-}\label{calo6}
\end{equation}
and ${\cal O}^{(2)}$ are codimension 2 totally real and complex curve orbits, respectively (observe that $\Psi_{\nu}^{-1}(\Omega_1)={\cal A}_{\nu}\cup{\cal O}_6$ (see (\ref{calanu})) with ${\cal O}_6=\Psi_{\nu}^{-1}({\cal O}_5)$ (see (\ref{calo5}))).

The manifold $M^{(1)}$ can be mapped onto $\Delta\times\CC\PP^1\subset\CC\PP^1\times\CC\PP^1$ by the inverse to a variant of the Segre map. \footnote{We are grateful to Stefan Nemirovski for showing us this realization of $M^{(1)}$.} Let $\left[\left(Z_0:Z_1\right),\left(W_0:W_1\right)\right]$ denote two pairs of homogeneous coordinates in $\CC\PP^1\times\CC\PP^1$, where the infinity in $\CC\PP^1$ is given by the vanishing of the coordinate that carries index 0. Consider the following map ${\cal S}$ from $\CC\PP^1\times\CC\PP^1$ to $\CC\PP^3$:
$$
\begin{array}{lll}
z_0&=&i\left(Z_0W_0-Z_1W_1\right),\\
z_1&=&Z_0W_1+Z_1W_0,\\
z_2&=&i\left(Z_0W_1-Z_1W_0\right),\\
z_3&=&Z_0W_0+Z_1W_1.
\end{array}
$$
It is straightforward to see that this map takes $\Delta\times\CC\PP^1$ biholomorphically onto $M^{(1)}$. Under the inverse map ${\cal S}^{-1}$ the action of ${\cal R}^{(1)}$ on $M^{(1)}$ is transformed into the following action of $SU_{1,1}/\{\pm\hbox{id}\}\simeq {\cal R}^{(1)}$ on $\Delta\times\CC\PP^1$: the element $g\{\pm\hbox{id}\}\in SU_{1,1}/\{\pm\hbox{id}\}$ acts on the vector $\left(Z_0:Z_1\right)$ by applying the matrix $g$ to the vector and on the vector $\left(W_0:W_1\right)$ by applying the matrix $\overline{g}$ to it. The map ${\cal S}^{-1}$ takes the orbit ${\cal O}_0^{(1)}$ into the $SU_{1,1}/\{\pm\hbox{id}\}$-orbit $\Delta\times\partial\Delta$, the orbit ${\cal O}_6$ into
$$
\left\{\left[\left(1:Z\right),\left(1:\overline{Z}\right)\right],\,|Z|<1\right\},
$$
and the orbit ${\cal O}^{(2)}$ into
$$
\left\{\left[\left(1:Z\right),\left(1:1/Z\right)\right],\,0<|Z|<1\right\}\cup\left\{\left[\left(1:0\right),\left(0:1\right)\right]\right\}.
$$
The domain $\Psi_{\nu}^{-1}(\Omega_1)$ is mapped by ${\cal S}^{-1}$ onto $\Delta\times\Delta$ and $D_1^{(2)}$ onto
$$
\Delta\times\Bigl(\left\{(1:W),\,|W|>1\right\}\cup\left\{(0:1)\right\}\Bigr)
$$
(hence each of $\Omega_1$, $D_1^{(2)}$ is equivalent to $\Delta^2$). For more general examples of this kind arising from actions of non-compact forms of complex reductive groups see \cite{AG}, \cite{FH}.

It is clear from the above description of $M^{(1)}$ that in order to obtain a hyperbolic ${\cal R}^{(1)}$-invariant submanifold of $M^{(1)}$ containing the Levi-flat orbit ${\cal O}_0^{(1)}$, one must remove from $M^{(1)}$ an ${\cal R}^{(1)}$-invariant neighborhood of either ${\cal O}_6$ or ${\cal O}^{(2)}$. Namely, each of the domains
\begin{equation}
\begin{array}{lll}
{\frak D}_s^{(1)}&:=&\Psi_{\nu}^{-1}(\Omega_{s,1})\cup D_1^{(2)}\cup {\cal O}_0^{(1)},  \quad-1<s<1,\\
\hat{\frak D}_t^{(1)}&:=&\Psi_{\nu}^{-1}(\Omega_1)\cup D_{1,t}^{(2)}\cup {\cal O}_0^{(1)},  \quad\,\,\, 1<t<\infty,\\
{\frak D}_{s,t}^{(1)}&:=&\Psi_{\nu}^{-1}(\Omega_{s,1})\cup D_{1,t}^{(2)}\cup {\cal O}_0^{(1)}, \quad -1\le s<1<t\le\infty,\\
&&\hbox{where $s=-1$ and $t=\infty$ do not hold simultaneously,}
\end{array}\label{frakdst1}
\end{equation}
is a (2,3)-manifold of this kind (see (\ref{omegast}), (\ref{omegat}), (\ref{dstinftyn})). Observe here that
$$
\begin{array}{ll}
\Psi_{\nu}^{-1}(\Omega_{s,1})=\Bigl\{(z_1,z_2,z_3)\in\CC^3:& s<|z_1|^2+|z_2|^2-|z_3|^2<1,\\
&\hspace{4cm}\hbox{Im}\,z_3<0\Bigr\}\cap {\cal Q}_{-}.
\end{array}
$$
Each of the groups $G\left({\frak D}_s^{(1)}\right)=\hbox{Aut}\left({\frak D}_s^{(1)}\right)$, $G\left(\hat{\frak D}_t^{(1)}\right)=\hbox{Aut}\left(\hat{\frak D}_t^{(1)}\right)$, $G\left({\frak D}_{s,t}^{(1)}\right)=\hbox{Aut}\left({\frak D}_{s,t}^{(1)}\right)$ coincides with ${\cal R}^{(1)}$.    
\vspace{0.5cm}

\noindent{\bf (d)} We now consider covers of
\begin{equation}
{\frak D}_{-1,\infty}^{(1)}:=\Psi_{\nu}^{-1}(\Omega_{-1,1})\cup D_{1,\infty}^{(2)}\cup {\cal O}_0^{(1)}.\label{frakdst1spec}
\end{equation}
For $n\ge 2$ the domain $\Sigma\setminus{\cal O}^{(2n)}$ (see (\ref{Sigma})) is an $n$-sheeted cover of ${\frak D}_{-1,\infty}^{(1)}$ with covering map $\Phi^{(n)}$. We equip $\Sigma\setminus{\cal O}^{(2n)}$ with the pull-back complex structure under $\Phi^{(n)}$. This complex structure extends the structure of each of $M_{\nu}^{(n)}$, $M_{\eta}^{(n)}$ (see {\bf (10)}) and can be extended to a complex structure on all of $\Sigma$. Let $M^{(n)}$ be the domain $\Sigma$ with this extended complex structure. The map $\Phi^{(n)}$ takes $M^{(n)}$ onto
$$
\Psi_{\nu}^{-1}(\Omega_{-1,1})\cup D_1^{(2)}\cup {\cal O}_0^{(1)},
$$
has ramification locus on ${\cal O}^{(2n)}$ and maps it in a 1-to-1 fashion onto ${\cal O}^{(2)}$.

Clearly, the group ${\cal R}^{(n)}$ acts on $M^{(n)}$. We will now describe the orbits of this action. The hypersurfaces  $\nu_{\alpha}^{(n)}$ for $-1<\alpha<1$ (see (\ref{nun})) and $\eta_{\alpha}^{(2n)}$ for $\alpha>1$ (see (\ref{eta2n})) are strongly pseudoconvex orbits; the hypersurface 
\begin{equation}
{\cal O}_0^{(n)}:=\left\{(z,w)\in M^{(n)}: |z|^n-|z|^{n-2}|w|^2=1\right\}\label{calo6n}
\end{equation}
is the unique Levi-flat orbit ($CR$-equivalent to $\Delta\times\partial\Delta$ for every $n$ and covering ${\cal O}_0^{(1)}$ by means of the $n$-to-1 map $\Phi^{(n)}$); the complex curve ${\cal O}^{(2n)}$ is the unique codimension 2 orbit.

We now introduce the following domains in $M^{(n)}$:      
$$
\begin{array}{lll}
{\frak D}_s^{(n)}&:=&\Biggl\{(\zeta:z:w)\in M^{(n)}:|z|^n-|z|^{n-2}|w|^2>\sqrt{(s+1)/2}|\zeta|^n\Biggr\}=\\
&&\Omega_{s,1}^{(n)}\cup D_1^{(2n)}\cup {\cal O}_0^{(n)},\, -1<s<1,\\
{\frak D}_{s,t}^{(n)}&:=&\Biggl\{(\zeta:z:w)\in M^{(n)}:\sqrt{(s+1)/2}|\zeta|^n<|z|^n-|z|^{n-2}|w|^2<\\
&&\sqrt{(t+1)/2}|\zeta|^n\Biggr\}=\Omega_{s,1}^{(n)}\cup D_{1,t}^{(2n)}\cup {\cal O}_0^{(n)},\,-1\le s<1<t\le \infty,\\
&&\hbox{where $s=-1$ and $t=\infty$ do not hold simultaneously}
\end{array}
$$
(see (\ref{omegastinftyn}), (\ref{dstinftyn})). Each of these domains is a (2,3)-manifold. Each of the groups $G\left({\frak D}_s^{(n)}\right)=\hbox{Aut}\left({\frak D}_s^{(n)}\right)$, $G\left({\frak D}_{s,t}^{(n)}\right)=\hbox{Aut}\left({\frak D}_{s,t}^{(n)}\right)$ coincides with ${\cal R}^{(n)}$.

Next, $\CC\times\Delta$ is an infinitely-sheeted cover of ${\frak D}_{-1,\infty}^{(1)}$ with covering map $\Lambda$ (see (\ref{maplambda})). We equip the domain $\CC\times\Delta$ with the pull-back complex structure under $\Lambda$ and denote the resulting manifold by $M^{(\infty)}$. Clearly, the complex structure of $M^{(\infty)}$ extends the structure of each of $M_{\nu}^{(\infty)}$, $M_{\eta}^{(\infty)}$ (see {\bf (10)}). The group ${\cal R}^{(\infty)}$ (see (\ref{autnuinfty1})) acts on $M^{(\infty)}$, and the orbits of this action are the strongly pseudoconvex hypersurfaces  $\nu_{\alpha}^{(\infty)}$ for $-1<\alpha<1$ (see (\ref{nuinfty})) and $\eta_{\alpha}^{(\infty)}$ for $\alpha>1$ (see (\ref{etainfty})), as well as the Levi-flat hypersurface 
\begin{equation}
{\cal O}_0^{(\infty)}:=\left\{(z,w)\in M^{(\infty)}: \exp(2\hbox{Re}\,z)(1-|w|^2)=1\right\}\label{calo6infty}
\end{equation}
(note that ${\cal O}_0^{(\infty)}$ is $CR$-equivalent to ${\cal O}_1$ -- see (\ref{calo1})).

We now introduce the following domains in $M^{(\infty)}$:      
$$
\begin{array}{lll}
{\frak D}_s^{(\infty)}&:=&\Biggl\{(z,w)\in M^{(\infty)}:\exp\left(2\hbox{Re}\,z\right)(1-|w|^2)>\sqrt{(s+1)/2}\Biggr\}=\\
&&\Omega_{s,1}^{(\infty)}\cup D_{1,\infty}^{(\infty)}\cup {\cal O}_0^{(\infty)}, \, -1<s<1,\\
{\frak D}_{s,t}^{(\infty)}&:=&\Biggl\{(z,w)\in M^{(\infty)}:\sqrt{(s+1)/2}<\exp\left(2\hbox{Re}\,z\right)(1-|w|^2)<\\
&&\sqrt{(t+1)/2}\Biggr\}=\Omega_{s,1}^{(\infty)}\cup D_{1,t}^{(\infty)}\cup {\cal O}_0^{(\infty)}, \, -1\le s<1<t\le\infty,\\
&&\hbox{where $s=-1$ and $t=\infty$ do not hold simultaneously}
\end{array}
$$
(see (\ref{omegastinftyn}), (\ref{dstinftyn})). Each of these domains is a (2,3)-manifold. The groups $G\left({\frak D}_s^{(\infty)}\right)=\hbox{Aut}\left({\frak D}_s^{(\infty)}\right)$, $G\left({\frak D}_{s,t}^{(\infty)}\right)=\hbox{Aut}\left({\frak D}_{s,t}^{(\infty)}\right)$ all coincide with ${\cal R}^{(\infty)}$.

\section{Strongly Pseudoconvex Orbits}\label{sectionstronglypseudoconvex}
\setcounter{equation}{0}

In this section we give a complete classification of (2,3)-manifolds $M$ for which every $G(M)$-orbit is a strongly pseudoconvex real hypersurface in $M$. In the formulation below we use the notation introduced in the previous section.

\begin{theorem}\label{classrealhypersurfaces}\sl Let $M$ be a (2,3)-manifold. Assume that the $G(M)$-orbit of every point in $M$ is a strongly pseudoconvex real hypersurface. Then $M$ is holomorphically equivalent to one of the following pairwise non-equivalent manifolds:
$$
\begin{array}{ll}

\hbox{(i)} & \hbox{$R_{b,s,t}$, $b\in\RR$, $|b|\ge1$, $b\ne 1$, with either $s=0$, $t=1$, or}\\
&\hbox{$s=1$, $1<t\le\infty$};\\
\hbox{(ii)} & \hbox{$U_{s,t}$, with either $s=0$, $t=1$, or $s=1$, $1<t\le\infty$};\\

\hbox{(iii)} & \hbox{${\frak S}_{s,t}$, with either $s=0$, $t=1$, or $s=1$, $1<t<\infty$};\\

\hbox{(iv)} & \hbox{${\frak S}_{s,t}^{(\infty)}$, with either $s=0$, $t=1$, or $s=1$, $1<t<\infty$};\\

\hbox{(v)} & \hbox{${\frak S}_{s,t}^{(n)}$, $n\ge 2$, with either $s=0$, $t=1$, or $s=1$, $1<t<\infty$};\\

\hbox{(vi)} & \hbox{$V_{b,s,1}$, $b>0$, with $e^{-2\pi b}<s<1$};\\

\hbox{(vii)} & \hbox{$E_{s,t}$, with $1\le s<t<\infty$};\\

\hbox{(viii)} & \hbox{$E_{s,t}^{(4)}$, with $1\le s<t<\infty$};\\

\hbox{(ix)} & \hbox{$E_{s,t}^{(2)}$, with $1\le s<t<\infty$};\\

\hbox{(x)} & \hbox{$\Omega_{s,t}$, with $-1\le s<t\le 1$};\\

\hbox{(xi)} & \hbox{$\Omega_{s,t}^{(\infty)}$, with $-1\le s<t\le 1$};\\

\hbox{(xii)} & \hbox{$\Omega_{s,t}^{(n)}$, $n\ge 2$, with $-1\le s<t\le 1$};\\

\hbox{(xiii)} & \hbox{$D_{s,t}$, with $1\le s<t\le\infty$};\\

\hbox{(xiv)} & \hbox{$D_{s,t}^{(\infty)}$, with $1\le s<t\le\infty$};\\

\hbox{(xv)} & \hbox{$D_{s,t}^{(n)}$, $n\ge 2$, with $1\le s<t\le\infty$}.
\end{array}
$$
\end{theorem}

\noindent {\bf Proof:} In \cite{C} E. Cartan classified all homogeneous 3-dimensional strongly pseudoconvex $CR$-manifolds. Since the $G(M)$-orbit of every point in $M$ is such a manifold, every $G(M)$-orbit is $CR$-equivalent to a manifold on Cartan's list. We reproduce Cartan's classification below together with the corresponding groups of $CR$-automorphisms. Note that all possible covers of the hypersurfaces $\chi$, $\mu_{\alpha}$, $\nu_{\alpha}$ and $\eta_{\alpha}$ appear below as explicitly realized in \cite{I3}.   
$$
\begin{array}{ll}
\hbox{(a)} & S^3;\\
\hbox{(b)}& {\cal L}_m:=S^3/\ZZ_m,\, m\in\NN, m\ge 2;\\
\hbox{(c)} & \sigma:=\left\{(z,w)\in\CC^2:\hbox{Re}\,w=|z|^2\right\};\\
\hbox{(d)} & \varepsilon_b:=\left\{(z,w)\in\CC^2: \hbox{Re}\,z=|w|^b,\, w\ne 0\right\}, \, b>0;\\
\hbox{(e)} & \omega:=\left\{(z,w)\in\CC^2:\hbox{Re}\,z=\exp\left(\hbox{Re}\,w\right)\right\};\\
\hbox{(f)} & \delta:=\left\{(z,w)\in\CC^2: |w|=\exp\left(|z|^2\right)\right\};\\
\hbox{(g)} & \hbox{$\tau_b$, $b\in\RR$, $|b|\ge 1$, $b\ne 1$ (see (\ref{tau}))};\\
\hbox{(h)} & \hbox{$\xi$ (see (\ref{xi}))};\\
\hbox{(j)} & \hbox{$\chi$ (see (\ref{chi}))};\\
\hbox{(j')} & \hbox{$\chi^{(\infty)}$ (see (\ref{chiinfty}))};\\
\hbox{(j'')} & \hbox{$\chi^{(n)}$, $n\ge 2$ (see (\ref{chin}))};
\end{array}
$$
$$
\begin{array}{ll}
\hbox{(k)} & \hbox{$\rho_b$, $b>0$ (see (\ref{rho}))};\\
\hbox{(l)} & \hbox{$\mu_{\alpha}$, $\alpha>1$ (see (\ref{mu}))};\\
\hbox{(l')} & \hbox{$\mu_{\alpha}^{(4)}$, $\alpha>1$ (see (\ref{mu4}))};\\\hbox{(l'')} & \hbox{$\mu_{\alpha}^{(2)}$, $\alpha>1$ (see (\ref{mu2}))};\\
\hbox{(m)} & \hbox{$\nu_{\alpha}$, $-1<\alpha<1$, (see (\ref{nu}))};\\
\hbox{(m')} & \hbox{$\nu_{\alpha}^{(\infty)}$, $-1<\alpha<1$, (see (\ref{nuinfty}))};\\
\hbox{(m'')} & \hbox{$\nu_{\alpha}^{(n)}$, $-1<\alpha<1$, $n\ge 2$ (see (\ref{nun}))};\\
\hbox{(n)} & \hbox{$\eta_{\alpha}$, $\alpha>1$, (see (\ref{eta}))};\\
\hbox{(n')} & \hbox{$\eta_{\alpha}^{(\infty)}$, $\alpha>1$, (see (\ref{etainfty}))};\\
\hbox{(n'')} & \hbox{$\eta_{\alpha}^{(n)}$, $\alpha>1$, $n\ge 2$ (see (\ref{eta2}), (\ref{eta2n}), (\ref{etan}))}.
\end{array}
$$

The above hypersurfaces are pairwise $CR$ non-equivalent. The corresponding groups of $CR$-automorphisms are as follows:
\vspace{0.5cm}

$\underline{\hbox{(a)}\,\hbox{Aut}_{CR}(S^3)}:$  maps of the form (\ref{autnu}), where the matrix $Q$ defined in (\ref{matq}) belongs to $SU_{2,1}$;
\vspace{0.5cm}

$\underline{\hbox{(b)}\,\hbox{Aut}_{CR}({\cal L}_m),\, m\ge 2}:$
$$
\left[\left(
\begin{array}{c}
z\\
w
\end{array}
\right)\right] \mapsto
\left[U
\left(
\begin{array}{c}
z\\
w
\end{array}
\right)\right],
$$
where $U\in U_2$, and $[(z,w)]\in{\cal L}_m$ denotes the equivalence class of $(z,w)\in S^3$ under the action of $\ZZ_m$ embedded in $U_2$ as a subgroup of scalar matrices;
\vspace{0.5cm}

$\underline{\hbox{(c)}\,\hbox{Aut}_{CR}(\sigma)}:$
\begin{equation}
\begin{array}{lll}
z & \mapsto & \lambda e^{i\psi}z+a,\\
w & \mapsto & \lambda^2w+2\lambda e^{i\psi}\overline{a}z+|a|^2+i\gamma,
\end{array}\label{groupofsigma}
\end{equation}
where $\lambda>0$, $\psi,\gamma\in\RR$, $a\in\CC$;
\vspace{0.5cm}

$\underline{\hbox{(d)}\,\hbox{Aut}_{CR}(\varepsilon_b)}:$
\begin{equation}
\begin{array}{lll}
z&\mapsto&\displaystyle \frac{\lambda z+i\beta}{i\mu z+\kappa},\\
\vspace{0mm}&&\\
w &\mapsto&\displaystyle
\frac{e^{i\psi}}{(i\mu z+\kappa)^{2/b}}w,
\end{array}\label{groupvarepsilon}
\end{equation}
where $\lambda,\beta,\mu,\kappa,\psi\in\RR$, $\lambda\kappa+\mu\beta=1$;
\vspace{0.5cm}

$\underline{\hbox{(e)}\,\hbox{Aut}_{CR}(\omega)}:$
\begin{equation}
\begin{array}{lll}
z&\mapsto&\displaystyle \frac{\lambda z+i\beta}{i\mu z+\kappa},\\
\vspace{0mm}&&\\
w &\mapsto&\displaystyle w-2\ln(i\mu z+\kappa)+i\gamma,
\end{array}\label{groupomega}
\end{equation}
where $\lambda,\beta,\mu,\kappa,\gamma\in\RR$, $\lambda\kappa+\mu\beta=1$;
\vspace{0.5cm}

$\underline{\hbox{(f)}\,\hbox{Aut}_{CR}(\delta)}:$
\begin{equation}
\begin{array}{lll}
z & \mapsto & e^{i\psi}z+a,\\
w & \mapsto &e^{i\theta}\exp\Bigl(2e^{i\psi}\overline{a}z+|a|^2\Bigr)w,
\end{array}\label{autdelta}
\end{equation}
where $\psi,\theta\in\RR$, $a\in\CC$;
\vspace{0.5cm}

$\underline{\hbox{(g)}\,\hbox{Aut}_{CR}(\tau_b)}:$ the group $G_b$ (see (\ref{auttau}));
\vspace{0.5cm}

$\underline{\hbox{(h)}\,\hbox{Aut}_{CR}(\xi)}:$ the group ${\frak G}$ (see (\ref{autxi}));
\vspace{0.5cm}

$\underline{\hbox{(j)}\,\hbox{Aut}_{CR}(\chi)}$ is generated by ${\cal R}_{\chi}$ (see (\ref{autchi})) and the map
\begin{equation}
\begin{array}{rrr}
z & \mapsto & z,\\
w & \mapsto & -w;
\end{array}\label{extramap}
\end{equation}
\vspace{0.5cm}

$\underline{\hbox{(j')}\,\hbox{Aut}_{CR}\left(\chi^{(\infty)}\right)}$ is generated by maps (\ref{autchiinfty1}) and the map
$$
\begin{array}{ccc}
z & \mapsto &\overline{z},\\
w & \mapsto & \overline{w};
\end{array}
$$
\vspace{0.5cm}

$\underline{\hbox{(j'')}\,\hbox{Aut}_{CR}\left(\chi^{(n)}\right),\, n\ge 2,}$ is generated by maps (\ref{autchin1}) and map (\ref{extramap});
\vspace{0.5cm}

$\underline{\hbox{(k)}\,\hbox{Aut}_{CR}(\rho_b)}:$ see (\ref{autrho}); 
\vspace{0.5cm}

$\underline{\hbox{(l)}\,\hbox{Aut}_{CR}(\mu_{\alpha})}:$ the group ${\cal R}_{\mu}$ (see (\ref{autmu})); 
\vspace{0.5cm}

$\underline{\hbox{(l')}\,\hbox{Aut}_{CR}\left(\mu_{\alpha}^{(4)}\right)}$ is generated by maps (\ref{autmu41}) and the map
$$
\begin{array}{lll}
z & \mapsto & \displaystyle i\frac{z(|z|^2+|w|^2)-\overline{w}}{\sqrt{1+(|z|^2+|w|^2)^2}},\\
\vspace{0cm}&&\\
w & \mapsto & \displaystyle i\frac{w(|z|^2+|w|^2)+\overline{z}}{\sqrt{1+(|z|^2+|w|^2)^2}};
\end{array}
$$
\vspace{0.5cm}

$\underline{\hbox{(l'')}\,\hbox{Aut}_{CR}\left(\mu_{\alpha}^{(2)}\right)}$ is generated by ${\cal R}_{\mu}^{(2)}$ (see (\ref{autmu2})) and the map
\begin{equation}
\begin{array}{rrr}
z_1 & \mapsto & -z_1,\\
z_2 & \mapsto & -z_2,\\
z_3 & \mapsto & -z_3;
\end{array}\label{mapminus}
\end{equation}
\vspace{0.5cm}

$\underline{\hbox{(m)}\,\hbox{Aut}_{CR}(\nu_{\alpha})}$ is generated by ${\cal R}_{\nu}$ (see (\ref{autnu})) and map (\ref{extramap});
\vspace{0.5cm}

$\underline{\hbox{(m')}\,\hbox{Aut}_{CR}\left(\nu_{\alpha}^{(\infty)}\right)}$ is generated by ${\cal R}^{(\infty)}$ (see (\ref{autnuinfty1})) and the map 
$$
\begin{array}{lll}
z & \mapsto & \displaystyle \overline{z} +\ln\left(-\frac{1+e^{2z}w(1-|w|^2)}{\sqrt{1-\exp\left(4\hbox{Re}\,z\right)(1-|w|^2)^2}}\right),\\
\vspace{0cm}&&\\
w & \mapsto & \displaystyle -\frac{\overline{w}+e^{2z}(1-|w|^2)}{1+e^{2z}w(1-|w|^2)};
\end{array}
$$
\vspace{0.5cm}

$\underline{\hbox{(m'')}\,\hbox{Aut}_{CR}\left(\nu_{\alpha}^{(n)}\right)\, n\ge 2}$ is generated by ${\cal R}^{(n)}$ (see (\ref{autnun1})) and the map
$$
\begin{array}{lll}
z & \mapsto & \displaystyle\overline{z}\sqrt[n]{\frac{\Bigl(1+z^{n-1}w(1-|w|^2/|z|^2)\Bigr)^2}{1-|z|^{2n}(1-|w|^2/|z|^2)^2}},\\
\vspace{0cm}&&\\
w & \mapsto & -\displaystyle\frac{\overline{w}/\overline{z}+z^n(1-|w|^2/|z|^2)}{1+z^{n-1}w(1-|w|^2/|z|^2)}\times\\
\vspace{0cm}&&\\
&&\displaystyle\hspace{1.2cm}\overline{z}\sqrt[n]{\frac{\Bigl(1+z^{n-1}w(1-|w|^2/|z|^2)\Bigr)^2}{1-|z|^{2n}(1-|w|^2/|z|^2)^2}};
\end{array}
$$
\vspace{0.5cm}

$\underline{\hbox{(n)}\,\hbox{Aut}_{CR}(\eta_{\alpha})}:$ the group ${\cal R}_{\eta}$ (see (\ref{auteta}));
\vspace{0.5cm}

$\underline{\hbox{(n')}\,\hbox{Aut}_{CR}\left(\eta_{\alpha}^{(\infty)}\right)}$ is generated by ${\cal R}^{(\infty)}$ (see (\ref{autnuinfty1})) and the map
$$
\begin{array}{lll}
z & \mapsto & \displaystyle 2z+\overline{z} +\ln\left(i\frac{1-|w|^2+e^{-2z}\overline{w}}{\sqrt{\exp\left(4\hbox{Re}\,z\right)(1-|w|^2)^2-1}}\right),\\
\vspace{0cm}&&\\
w & \mapsto & \displaystyle \frac{1+e^{2z}w(1-|w|^2)}{\overline{w}+e^{2z}(1-|w|^2)};
\end{array}
$$

\vspace{0.5cm}

$\underline{\hbox{(n'')}\,\hbox{Aut}_{CR}\left(\eta_{\alpha}^{(2)}\right)\, }$ is generated by ${\cal R}^{(1)}$ (see {\bf (10)}) and map (\ref{mapminus});
\vspace{0.5cm}

$\underline{\hbox{(n'')}\,\hbox{Aut}_{CR}\left(\eta_{\alpha}^{(2n)}\right)\, n\ge 2}$ is generated by ${\cal R}^{(n)}$ (see (\ref{autnun1})) and the map
$$
\begin{array}{lll}
z & \mapsto & \displaystyle z^2\overline{z}\sqrt[n]{\frac{\Bigl(1-|w|^2/|z|^2+z^{-n}\overline{w}/\overline{z}\Bigr)^2}{|z|^{2n}(1-|w|^2/|z|^2)^2-1}},\\
\vspace{0cm}&&\\
w & \mapsto & \displaystyle\frac{1+z^{n-1}w(1-|w|^2/|z|^2)}{\overline{w}/\overline{z}+z^{n}(1-|w|^2/|z|^2)}\times\\
\vspace{0cm}&&\\
&&\displaystyle z^2\overline{z}\sqrt[n]{\frac{\Bigl(1-|w|^2/|z|^2+z^{-n}\overline{w}/\overline{z}\Bigr)^2}{|z|^{2n}(1-|w|^2/|z|^2)^2-1}};
\end{array}
$$
\vspace{0.5cm}

$\underline{\hbox{(n'')}\,\hbox{Aut}_{CR}\left(\eta_{\alpha}^{(n)}\right),\, \hbox{$n\ge 3$ is odd}}:$ this group is isomorphic to ${\cal R}^{(n)}$ and consists of all lifts from the domain $D_{1,\infty}$ (see (\ref{dst})) to $M_{\eta}^{'(n)}$ (see {\bf (10)}) of maps from ${\cal R}_{\eta}$ (see (\ref{auteta})). 
\vspace{0.5cm}

We will now show that the presence of an orbit of a particular kind in $M$ determines the group $G(M)$ as a Lie group. Fix $p\in M$ and suppose that $O(p)$ is $CR$-equivalent to ${\frak m}$, where ${\frak m}$ is one of the hypersurfaces listed above in (a)--(n''). In this case we say that ${\frak m}$ is the {\it model}\, for $O(p)$.    
Since $G(M)$ acts properly and effectively on $O(p)$, the $CR$-equivalence induces an isomorphism between $G(M)$ and a closed connected 3-dimensional subgroup $R_{\frak m}$ of the Lie group $\hbox{Aut}_{CR}({\frak m})$, that acts transitively on ${\frak m}$ (note that the Lie group topology of $\hbox{Aut}_{CR}({\frak m})$ coincides with the compact-open topology -- see e.g. \cite{Sch}). The subgroup $R_{\frak m}$ a priori depends on the choice of $CR$-equivalence between $O(p)$ and ${\frak m}$, but, as we will see below, this dependence is insignificant.

We will now list all possible groups $R_{\frak m}$ for each model in (a)-(n''). In the following lemma ${\cal P}$ denotes the right half-plane $\left\{z\in\CC:\hbox{Re}\,z>0\right\}$.

\begin{lemma}\label{groupsofmodels}\sl We have
$$
\begin{array}{ll}
\hbox{(A)} & \hbox{$R_{\frak m}=\hbox{Aut}_{CR}({\frak m})^0$, if ${\frak m}$ is one of (g)-(n'')};\\
\hbox{(B)} & \hbox{$R_{S^3}$ is conjugate in $\hbox{Aut}_{CR}(S^3)$ to $SU_2$};\\
\hbox{(C)} & R_{{\cal L}_m}=SU_2/(SU_2\cap\ZZ_m),\, m\ge 2;\\
\hbox{(D)} & \hbox{$R_{\sigma}$ is the Heisenberg group, that is, it consists of all elements}\\
&\hbox{of $\hbox{Aut}_{CR}(\sigma)$ with $\lambda=1$, $\psi=0$ in formula (\ref{groupofsigma})};\\
\hbox{(E)} & \hbox{$R_{\varepsilon_b}$ either is the subgroup of $\hbox{Aut}_{CR}(\varepsilon_b)$ corresponding to a subgroup}\\
&\hbox{of $\hbox{Aut}({\cal P})$, conjugate in $\hbox{Aut}({\cal P})$ to the subgroup ${\cal T}$ given by}
\end{array}
$$
\begin{equation}
z\mapsto\lambda z+i\beta,\label{groupt}
\end{equation}
$$
\begin{array}{ll}
&\hbox{where $\lambda>0$, $\beta\in\RR$, or, for $b\in\QQ$, is the subgroup ${\frak V}_b$ given by $\psi=0$}\\
&\hbox{in formula (\ref{groupvarepsilon})};\\
\hbox{(F)} &\hbox{$R_{\omega}$ either is the subgroup of $\hbox{Aut}_{CR}(\omega)$ corresponding to a subgroup}\\
&\hbox{of $\hbox{Aut}({\cal P})$ conjugate in $\hbox{Aut}({\cal P})$ to the subgroup ${\cal T}$ specified in (E),}\\
&\hbox{or is the subgroup ${\frak V}_{\infty}$ given by $\gamma=0$ in formula (\ref{groupomega})};\\
\hbox{(G)} & \hbox{$R_{\delta}$ coincides with the subgroup of $\hbox{Aut}_{CR}(\delta)$ given by $\psi=0$ in}\\
&\hbox{formula (\ref{autdelta})}.
\end{array}
$$
\end{lemma}

\noindent {\bf Proof:} Case (A) is clear since in (g)-(n'') we have $\hbox{dim}\,\hbox{Aut}_{CR}({\frak m})=d(M)=3$. Further, in case (B) the orbit $O(p)$ is compact and, since $I_p$ is compact as well, it follows that $G(M)$ is compact. Thus $R_{S^3}$ is compact, and hence it is conjugate to a subgroup of $U_2$, which is a maximal compact subgroup in $\hbox{Aut}_{CR}(S^3)$. Since $R_{S^3}$ is 3-dimensional, it is in fact conjugate to $SU_2$, as required. In case (C) the group $R_{{\cal L}_m}$ is of codimension 1 in $\hbox{Aut}_{CR}\left({\cal L}_m\right)=U_2/\ZZ_m$, hence $R_{{\cal L}_m}=SU_2/(SU_2\cap\ZZ_m)$.

Further, (D) and (G) follow since $R_{\sigma}$ and $R_{\delta}$ act transitively on $\sigma$ and $\delta$, respectively. In cases (E) and (F) note that every codimension 1 subgroup of $\hbox{Aut}({\cal P})$ is conjugate in $\hbox{Aut}({\cal P})$ to the subgroup ${\cal T}$ defined in (\ref{groupt}). The only codimension 1 subgroups of $\hbox{Aut}_{CR}(\varepsilon_b)$ and $\hbox{Aut}_{CR}(\omega)$ that do not arise from codimension 1 subgroups of $\hbox{Aut}({\cal P})$ are ${\frak V}_b$ and ${\frak V}_{\infty}$, respectively. Observe, however, that ${\frak V}_b$ is not closed in $\hbox{Aut}_{CR}(\varepsilon_b)$ unless $b\in\QQ$.   

The lemma is proved.\qed
\vspace{0.5cm}

Lemma \ref{groupsofmodels} implies, in particular, that if for some point $p\in M$ the model for $O(p)$ is $S^3$, then $M$ admits an effective action of $SU_2$. Therefore, $M$ is holomorphically equivalent to one of the manifolds listed in \cite{IKru2}. However, none of the (2,3)-manifolds on the list has a spherical orbit. Hence we have ruled out case (a). 

We now observe -- directly from the explicit forms of the $CR$-automorphism groups of the models listed above -- that for each ${\frak m}$ every element of $\hbox{Aut}_{CR}({\frak m})$ extends to a holomorphic automorphism of a certain complex manifold $M_{\frak m}$ containing ${\frak m}$, such that every $R_{\frak m}$-orbit $O$ in $M_{{\frak m}}$ is strongly pseudoconvex and exactly one of the following holds: (a) $O$ is $CR$-equivalent to ${\frak m}$ (cases (b)--(k)); (b) $O$ belongs to the same family to which ${\frak m}$ belongs and the $R_{\frak m}$-orbits are pairwise $CR$ non-equivalent (cases (l)--(n'')). The manifolds $M_{\frak m}$ are as follows:
\vspace{0.5cm}

\noindent (b) $M_{{\cal L}_m}=\CC^2\setminus\{0\}/\ZZ_m$, $m\ge 2$;
\vspace{0.5cm}

\noindent (c) $M_{\sigma}=\CC^2$;
\vspace{0.5cm}

\noindent (d) $M_{\varepsilon_b}=\left\{(z,w)\in\CC^2: \hbox{Re}\,z> 0,\,w\ne 0\right\}$;
\vspace{0.5cm}

\noindent (e) $M_{\omega}=\left\{(z,w)\in\CC^2: \hbox{Re}\,z> 0\right\}$;
\vspace{0.5cm}

\noindent (f) $M_{\delta}=\left\{(z,w)\in\CC^2:w\ne 0\right\}$;
\vspace{0.5cm}

\noindent (g) $M_{\tau_b}=\left\{(z,w)\in\CC^2:\hbox{Re}\,z>0,\,\hbox{Re}\,w> 0\right\}$;
\vspace{0.5cm}

\noindent (h) $M_{\xi}=\left\{(z,w)\in\CC^2:\hbox{Re}\,w>0\right\}$;
\vspace{0.5cm}

\noindent (j) $M_{\chi}=\CC^2\setminus\left\{(z,w)\in\CC^2:\hbox{Re}\,z=0,\,\hbox{Re}\,w=0\right\}$;
\vspace{0.5cm}

\noindent (j') $M_{\chi^{(\infty)}}=M_{\chi}^{(\infty)}$ (see {\bf (4)});
\vspace{0.5cm}

\noindent (j'') $M_{\chi^{(n)}}=M_{\chi}^{(n)}$, $n\ge 2$ (see {\bf (4)});
\vspace{0.5cm}

\noindent (k) $M_{\rho_b}=\CC^2\setminus\left\{(z,w)\in\CC^2:\hbox{Re}\,z=0,\,\hbox{Re}\,w=0\right\}$;
\vspace{0.5cm}

\noindent (l) $M_{\mu_{\alpha}}=\displaystyle\bigcup_{\alpha>1}\mu_{\alpha}=\CC\PP^2\setminus\RR\PP^2$; 
\vspace{0.5cm}

\noindent (l') $M_{\mu_{\alpha}^{(4)}}=\displaystyle\bigcup_{\alpha>1}\mu_{\alpha}^{(4)}=M_{\mu}^{(4)}$ (see {\bf (7)});
\vspace{0.5cm}

\noindent (l'') $M_{\mu_{\alpha}^{(2)}}=\displaystyle\bigcup_{\alpha>1}\mu_{\alpha}^{(2)}={\cal Q}_{+}\setminus\RR^3$ (see (\ref{quadricplus})); 
\vspace{0.5cm}

\noindent (m)
$M_{\nu_{\alpha}}=\displaystyle\bigcup_{-1<\alpha<1}\nu_{\alpha}=\Omega_{-1,1}$ (see (\ref{omegast}));
\vspace{0.5cm}

\noindent (m') $M_{\nu_{\alpha}^{(\infty)}}=\displaystyle\bigcup_{-1<\alpha<1}\nu_{\alpha}^{(\infty)}=M_{\nu}^{(\infty)}=\Omega^{(\infty)}_{-1,1}$ (see (\ref{omegastinftyn})); 
\vspace{0.5cm}

\noindent (m'') $M_{\nu_{\alpha}^{(n)}}=\displaystyle\bigcup_{-1<\alpha<1}\nu_{\alpha}^{(n)}=M_{\nu}^{(n)}=\Omega^{(n)}_{-1,1}$, $n\ge 2$ (see (\ref{omegastinftyn}));
\vspace{0.5cm}

\noindent (n) $M_{\eta_{\alpha}}=\displaystyle\bigcup_{\alpha>1}\eta_{\alpha}=D_{1,\infty}$ (see (\ref{dst}));
\vspace{0.5cm}

\noindent (n') $M_{\eta_{\alpha}^{(\infty)}}=\displaystyle\bigcup_{\alpha>1}\eta_{\alpha}^{(\infty)}=M_{\eta}^{(\infty)}=D^{(\infty)}_{1,\infty}$ (see (\ref{dstinftyn})); 
\vspace{0.5cm}

\noindent (n'') $M_{\eta_{\alpha}^{(2)}}=\displaystyle\bigcup_{\alpha>1}\eta_{\alpha}^{(2)}={\cal A}_{\eta}$ (see (\ref{calaeta})); 
\vspace{0.5cm}

\noindent (n'') $M_{\eta_{\alpha}^{(2n)}}=\displaystyle\bigcup_{\alpha>1}\eta_{\alpha}^{(2n)}=M_{\eta}^{(n)}=D^{(2n)}_{1,\infty}$, $n\ge 2$ (see (\ref{dstinftyn}));
\vspace{0.5cm}

\noindent (n'') $M_{\eta_{\alpha}^{(n)}}=\displaystyle\bigcup_{\alpha>1}\eta_{\alpha}^{(n)}=M_{\eta}^{'(n)}=D^{(n)}_{1,\infty}$, $n\ge 3$ is odd (see {\bf (10)}).
\vspace{0.5cm}

In each of cases (b)--(k) every two $R_{\frak m}$-orbits are $CR$-equivalent (and equivalent to ${\frak m}$) by means of an automorphism of $M_{\frak m}$ of one of the simple forms specified below:
\begin{equation}
\begin{array}{ll}
\hbox{(b)}& [(z,w)]  \mapsto  [(az,aw)],\, a>0;\\
\vspace{0cm}&\\ 
\hbox{(c)}& z  \mapsto  z,\,w  \mapsto  w+a,\,a\in\RR;\\
\vspace{0cm}&\\     
\hbox{(d)} & z  \mapsto  az,\, w  \mapsto  w,\,a>0;\\
\vspace{0cm}&\\    
\hbox{(e)}& \hbox{as in (d)};\\
\vspace{0cm}&\\ 
\hbox{(f)} & z  \mapsto  z,\, w  \mapsto  aw,\,a>0;\\
\vspace{0cm}&\\ 
\hbox{(g)} & \hbox{as in (f)};\\
\vspace{0cm}&\\ 
\hbox{(h)} & z  \mapsto  az,\, w  \mapsto  aw,\,a>0;\\
\vspace{0cm}&\\ 
\hbox{(j)} & \hbox{as in (h)};\\
\vspace{0cm}&\\ 
\hbox{(j')} & z  \mapsto  z+a,\, w  \mapsto  e^aw,\,a\in\RR;\\
\vspace{0cm}&\\ 
\hbox{(j'')}& z \mapsto a\hbox{Re}\,z+ia^n\hbox{Im}\,z,\, w \mapsto a\hbox{Re}\,w+ia^n\hbox{Im}\,w,\,a>0;\\
\vspace{0cm}&\\ 
\hbox{(k)}& \hbox{as in (h)}.
\end{array}\label{simplelist}
\end{equation}
 
We will now show how strongly pseudoconvex orbits can be glued together to form (2,3)-manifolds. The procedure comprises the following steps.
\vspace{0cm}\\

\noindent (I). Start with a real hypersurface orbit $O(p)$ with model ${\frak m}$ and consider a real-analytic $CR$-isomorphism $f:O(p)\ra{\frak m}$. Clearly, $f$ satisfies
\begin{equation}
f(gq)=\varphi(g)f(q),\label{equivar}
\end{equation} 
for all $g\in G(M)$ and $q\in O(p)$, where $\varphi:G(M)\ra R_{\frak m}$ is a Lie group isomorphism.
\vspace{0cm}\\

\noindent (II). Observe that $f$ can be extended to a biholomorphic map from a $G(M)$-invariant connected neighborhood of $O(p)$ in $M$ onto an $R_{\frak m}$-invariant neighborhood of ${\frak m}$ in the corresponding manifold $M_{{\frak m}}$. If $G(M)$ is compact (in which case ${\frak m}$ is one of ${\cal L}_m$ with $m\ge 2$, $\mu_{\alpha}$, $\mu_{\alpha}^{(2)}$, $\mu_{\alpha}^{(4)}$ with $\alpha>1$), then every neighborhood of $O(p)$ contains a $G(M)$-invariant neighborhood. In this case, we extend $f$ biholomorphically to some neighborhood of $O(p)$ (this can be done due to the real-analyticity of $f$) and choose a $G(M)$-invariant neighborhood in it. 

We now assume that $G(M)$ is non-compact. In this case it will be more convenient for us to extend the inverse map ${\frak F}:=f^{-1}$. First of all, extend ${\frak F}$ to some neighborhood $U$ of ${\frak m}$ in $M_{{\frak m}}$ to a biholomorphic map onto a neighborhood $W$ of $O(p)$ in $M$. It can be seen from the explicit form of the $R_{\frak m}$-action on $M_{{\frak m}}$ that $U$ can be chosen to satisfy the following condition that we call Condition (${}*{}$): for every two points $s_1,s_2\in U$ and every element $h\in R_{\frak m}$ such that $hs_1=s_2$ there exists a curve $\gamma\subset U$ joining $s_1$ with a point in ${\frak m}$ for which $h\gamma\subset U$ (clearly, $h\gamma$ is a curve joining $s_2$ with a point in ${\frak m}$).  

To extend ${\frak F}$ to a $R_{\frak m}$-invariant neighborhood of ${\frak m}$, fix $s\in U$ and $s_0\in {\cal O}(s)$, where ${\cal O}(y)$ denotes the $R_{\frak m}$-orbit of a point $y\in M_{{\frak m}}$.  Choose $h_0\in R_{\frak m}$ such that $s_0=h_0s$ and define ${\frak F}(s_0):=\varphi^{-1}(h_0){\frak F}(s)$. We will now show that ${\frak F}$ is well-defined. Suppose that for some $s_1,s_2\in U$ and $h_1,h_2\in R_{\frak m}$ we have $s_0=h_1s_1=h_2s_2$. To show that $\varphi^{-1}(h_1){\frak F}(s_1)=\varphi^{-1}(h_2){\frak F}(s_2)$ we set $h:=h_2^{-1}h_1$ and, according to Condition (${}*{}$), find a curve $\gamma\subset U$ that joins $s_1$ with a point in ${\frak m}$ and such that $h\gamma\subset U$. 

Clearly, for all $q\in{\frak m}$ we have 
\begin{equation}
{\frak F}(hq)=\varphi^{-1}(h){\frak F}(q).\label{equivar1}
\end{equation}
Consider the open set $h^{-1}U\cap U$ and let $U_h$ be its connected component containing ${\frak m}$. For $q\in U_h$ identity (\ref{equivar1}) holds. It now follows from the existence of a curve $\gamma$ as above that $s_1\in U_h$. Thus, (\ref{equivar1}) holds for $q=s_1$, and we have shown that ${\frak F}$ is well-defined at $s_0$. The same argument gives that for $s_0\in U$ our definition agrees with the original value ${\frak F}(s_0)$. Thus, we have extended ${\frak F}$ to $U':=\cup_{s\in U}{\cal O}(s)$. The extended map is locally biholomorphic, satisfies (\ref{equivar1}), and maps $U'$ onto a $G(M)$-invariant neighborhood $W'$ of $O(p)$ in $M$. We will now show that the extended map is 1-to-1 on an $R_{{\frak m}}$-invariant neighborhood of ${\frak m}$ contained in $U'$. 

Suppose that for some $s_0,s_0'\in U'$, $s_0\ne s_0'$, we have ${\frak F}(s_0)={\frak F}(s_0')$. This can only occur if $s_0$ and $s_0'$ lie in the same $R_{\frak m}$-orbit, and therefore there exist a point $s\in U$ and elements $h,h'\in R_{\frak m}$ such that $s_0=hs$, $s_0'=h's$. Then  $h^{'-1}h\not\in J_s$ and $\varphi^{-1}(h^{'-1}h)\in I_{{\frak F}(s)}$, where $J_s$ denotes the isotropy subgroup of $s$ under the $R_{\frak m}$-action. At the same time, we have $\varphi^{-1}(J_s)\subset I_{{\frak F}(s)}$. Thus, $I_{{\frak F}(s)}$ contains more points than $J_s$. Observe also that if ${\frak m}'$ is the model for $O\left({\frak F}(s)\right)$, then $R_{{\frak m}'}$ is isomorphic to $R_{\frak m}$. 

Assume first that ${\cal O}(s)$ is non-spherical. It follows from the explicit forms of the models and the corresponding groups (see Lemma \ref{groupsofmodels}) that if for two locally $CR$-equivalent non-spherical models ${\frak m}_1$, ${\frak m}_2$ the groups $R_{{\frak m}_1}$ and $R_{{\frak m}_2}$ are isomorphic and the isotropy subgroup of a point in ${\frak m}_1$ contains more points than that of a point in ${\frak m}_2$, then ${\frak m}_1=\eta_{\alpha}^{(n)}$ and ${\frak m}_2=\eta_{\alpha}^{(2n)}$ for some $\alpha$ and odd $n\ge 1$ (here we set $\eta_{\alpha}^{(1)}:=\eta_{\alpha}$). Hence ${\cal O}(s)=\eta_{\alpha}^{(2n)}$, and the model for $O({\frak F}(s))$ is $\eta_{\alpha}^{(n)}$ for some $\alpha$ and odd $n\ge 1$; consequently, ${\frak m}=\eta_{\beta}^{(2n)}$ for some $\beta$. If there is a neighborhood of $p$ not containing a point $q$ such that the model for $O(q)$ is some $\eta_{\gamma}^{(n)}$, then ${\frak F}$ is biholomorphic on an $R_{{\frak m}}$-invariant open subset of $U'$. Suppose now that in every neighborhood of $p$ (that we assume to be contained in $W$) there is a point $q$ such that ${\cal O}(f(q))=\eta_{\gamma}^{(2n)}$ and the model for $O(q)$ is $\eta_{\gamma}^{(n)}$ for some $\gamma$. Note that $I_q$ consists of two elements and $J_{f(q)}$ is trivial. Choose a sequence of such points $\{q_j\}\subset W$ converging to $p$. Let $g_j$ be the non-trivial element of $I_{q_j}$. Since the action of $G(M)$ on $M$ is proper and $I_p$ is trivial, the sequence $\{g_j\}$ converges to the identity in $G(M)$. At the same time, the sequence $\{f(q_j)\}$ converges to $f(p)$ and therefore $\varphi(g_j)f(q_j)$ lies in $U$ for large $j$. For large $j$ we have ${\frak F}(\varphi(g_j)f(q_j))=q_j$. Since $\varphi(g_j)$ is a non-trivial element in $R_{{\frak m}}$, the point $\varphi(g_j)f(q_j)$ does not coincide with $f(q_j)$. Thus, we have found two distinct points in $U$ (namely, $f(q_j)$ and $\varphi(g_j)f(q_j)$ for large $j$) mapped by ${\frak F}$ into the same point in $W$, which contradicts the fact that ${\frak F}$ is 1-to-1 on $U$.

Assume now that ${\cal O}(s)$ is spherical. It follows from the explicit forms of the models and the corresponding groups (see Lemma \ref{groupsofmodels}) that if for two spherical models ${\frak m}_1$, ${\frak m}_2$ the groups $R_{{\frak m}_1}$ and $R_{{\frak m}_2}$ are isomorphic and the isotropy subgroup of a point in ${\frak m}_1$ contains more points than that of a point in ${\frak m}_2$, then ${\frak m}_1=\varepsilon_{n/k_1}$, $R_{{\frak m}_1}={\frak V}_{n/k_1}$, and ${\frak m}_2=\varepsilon_{n/k_2}$, $R_{{\frak m}_2}={\frak V}_{n/k_2}$, where $n,k_1,k_2\in\NN$, $(n,k_1)=1$, $(n,k_2)=1$, $k_1>k_2$. Hence ${\frak m}=\varepsilon_{n/k_2}$, ${\cal O}(s)$ is equivalent to $\varepsilon_{n/k_2}$ by means of a map of the form (d) on list (\ref{simplelist}), and the model for $O({\frak F}(s))$ is $\varepsilon_{n/k_1}$ for some $n,k_1,k_2$ as above. If there is a neighborhood of $p$ not containing a point $q$ such that the model for $O(q)$ is some $\varepsilon_{n/k}$, with $k\in\NN$, $(n,k)=1$, $k>k_2$, then ${\frak F}$ is biholomorphic on an $R_{{\frak m}}$-invariant open subset of $U'$. Suppose now that in every neighborhood of $p$ (that we assume to be contained in $W$) there is a point $q$ such that the model for $O(q)$ is $\varepsilon_{n/k}$, with $k\in\NN$, $(n,k)=1$, $k>k_2$. Choose a sequence of such points $\{q_j\}\subset W$ converging to $p$. Since the action of $G(M)$ on $M$ is proper, the isotropy subgroups $I_{q_j}$ converge to $I_p$. Every subgroup $I_{q_j}$ contains more points than $J_{f(q_{j_k})}$, and therefore for a subsequence $\{j_k\}$ of the sequence of indices $\{j\}$ there is a sequence of elements $\{g_{j_k}\}$ with $g_{j_k}\in I_{q_{j_k}}$, $\varphi(g_{j_k})\not\in J_{f(q_{j_k})}$ and such that $\{g_{j_k}\}$ converges to an element of $I_p$. Arguing now as in the non-spherical case, we obtain a contradiction with the fact that ${\frak F}$ is 1-to-1 on $U$.

Hence we have shown that $f$ can be extended to a biholomorphic map satisfying (\ref{equivar}) between a $G(M)$-invariant neighborhood of $O(p)$ in $M$ and a $R_{{\frak m}}$-invariant neighborhood of ${\frak m}$ in $M_{{\frak m}}$.    
\vspace{0cm}\\  

\noindent (III). Consider a maximal $G(M)$-invariant domain $D\subset M$ from which there exists a biholomorphic map $f$ onto an $R_{\frak m}$-invariant domain in $M_{\frak m}$ satisfying (\ref{equivar}) for all $g\in G(M)$ and $q\in D$. The existence of such a domain is guaranteed by the previous step. Assume that $D\ne M$ and consider $x\in\partial D$. Let ${\frak m}_1$ be the model for $O(x)$ and let $f_1:O(x)\ra{\frak m}_1$ be a real-analytic $CR$-isomorphism satisfying (\ref{equivar}) for all $g\in G(M)$, $q\in O(x)$ and some isomorphism $\varphi_1: G(M)\ra R_{{\frak m}_1}$ in place of $\varphi$. 
As in (II), extend $f_1$ to a biholomorphic map from a connected $G(M)$-invariant neighborhood $V$ of $O(x)$ onto an $R_{{\frak m}_1}$-invariant neighborhood of ${\frak m}_1$ in $M_{{\frak m}_1}$. The extended map satisfies (\ref{equivar}) for all $g\in G(M)$, $q\in V$ and $\varphi_1$ in place of $\varphi$. Consider $s\in V\cap D$. The maps $f$ and $f_1$ take $O(s)$ onto an $R_{\frak m}$-orbit ${\frak m}'$ in $M_{\frak m}$ and an $R_{{\frak m}_1}$-orbit ${\frak m}_1'$ in $M_{{\frak m}_1}$, respectively. Then $F:=f_1\circ f^{-1}$ establishes a $CR$-isomorphism between ${\frak m}'$ and ${\frak m}_1'$. Therefore, ${\frak m}_1$ lies in $M_{\frak m}$, that is, we have $M_{\frak m}=M_{{\frak m}_1}$. Moreover, $F$ is either an element of $\hbox{Aut}_{CR}({\frak m}')$ (if ${\frak m}'={\frak m}_1'$), or is a composition of an element of $\hbox{Aut}_{CR}({\frak m}')$ and a non-trivial map from list (\ref{simplelist}) that takes ${\frak m}'$ onto ${\frak m}_1'$ (if ${\frak m}'\ne{\frak m}_1'$); the latter is only possible in cases (b)--(k). It now follows from the explicit forms of $CR$-automorphisms of the models and the maps on list (\ref{simplelist}) that $F$ extends to a holomorphic automorphism of $M_{\frak m}$.   
\vspace{0cm}\\   

\noindent (IV). Since $O(x)$ is strongly pseudoconvex and closed in $M$, for $V$ sufficiently small we have $V=V_1\cup V_2\cup O(x)$, where $V_j$ are open connected non-intersecting sets. Furthermore, if $V$ is sufficiently small, then each $V_j$ is either a subset of $D$ or disjoint from it. Suppose first that there is only one $j$ for which $V_j\subset D$. In this case $V\cap D$ is connected and $V\setminus(D\cup O(x))\ne\emptyset$. Setting now
\begin{equation}
\tilde f:=\Biggl\{\begin{array}{l}
f\hspace{1.7cm}\hbox{on $D$},\\
F^{-1}\circ f_1\hspace{0.45cm}\hbox{on $V$},
\end{array}\label{ext}
\end{equation}
we obtain a biholomorphic extension of $f$ to $D\cup V$. By construction, $\tilde f$ satisfies (\ref{equivar}) for $g\in G(M)$ and $q\in D\cup V$. Since $D\cup V$ is strictly larger than $D$, we obtain a contradiction with the maximality of $D$. Thus, 
in this case $D=M$, and hence $M$ is holomorphically equivalent to an $R_{\frak m}$-invariant domain in $M_{\frak m}$ (all such domains will explicitly appear below).

Suppose now that $V_j\subset D$ for $j=1,2$. Applying formula (\ref{ext}) to suitable $f_1$ and $F$, we can extend $f|_{V_1}$ and $f|_{V_2}$ to biholomorphic maps $\hat f_1$, $\hat f_2$, respectively, from a neighborhood  of $O(x)$ into $M_{\frak m}$; each of these maps satisfies (\ref{equivar}). Let $\hat{\frak m}_j:=\hat f_j(O(x))$, $j=1,2$. Then $\partial D=\hat{\frak m}_1\cup\hat{\frak m}_2$, $\hat{\frak m}_1\ne\hat{\frak m}_2$, and $M\setminus O(x)$ is holomorphically equivalent to $D$. The map $\hat F:=\hat f_2\circ\hat f_1^{-1}$ is a $CR$-isomorphism from $\hat{\frak m}_1$ onto $\hat{\frak m}_2$ (hence ${\frak m}$ is one of the hypersurfaces occurring in cases (b)--(k)), and  $M$ is holomorphically equivalent to the manifold $M_{\hat F}$ obtained from $\overline{D}$ by identifying $\hat{\frak m}_1$ with $\hat{\frak m}_2$ by means of $\hat F$. Since the action of $R_{\frak m}$ on $M\setminus O(x)$ extends to an action on $M$, the map $\hat F$ is $R_{\frak m}$-equivariant. In each of cases (b)--(f) this implies that $\hat F$ extends to a holomorphic automorphism of $M_{\frak m}$ of a simple form (similar to the corresponding form on list (\ref{simplelist})). Let $\Gamma$ denote the group of automorphisms of $M_{\frak m}$ generated by $\hat F$. It follows from the explicit forms of $\hat F$ and $M_{\frak m}$ in each of cases (b)--(f) that $\Gamma$ acts freely properly discontinuously on $M_{\frak m}$ and that $M_{\frak m}$ covers $M$, with $\Gamma$ being the group of deck transformations of the covering map. Observe, however, that for every ${\frak m}$ the manifold $M_{\frak m}$ is not hyperbolic. Next, in cases (g)--(k) the $R_{\frak m}$-equivariance of $\hat F$ implies that $\hat F=\hbox{id}$, which is impossible. These contradictions show that exactly one of $V_j$, $j=1,2$, is a subset of $D$, and hence $M$ is holomorphically equivalent to an $R_{\frak m}$-invariant domain in $M_{\frak m}$.
\vspace{0.5cm}

All hyperbolic $R_{\frak m}$-invariant domains in each of cases (b)-(n'') are described as follows:

\begin{equation}
\begin{array}{ll}
\hspace{-4.3cm}\hbox{(b)} & \hspace{-3.8cm}\hbox{$S_{m,s,t}:=\left\{(z,w)\in\CC^2:s<|z|^2+|w|^2<t\right\}/\ZZ_m$},\\
&\hspace
{-2.2cm}\hbox{$0\le s<t<\infty$};
\end{array}\label{smst}
\end{equation}

\noindent (c) $\left\{(z,w)\in\CC^2: s+|z|^2<\hbox{Re}\,w<t+|z|^2\right\}$, $-\infty<s<t\le\infty$;
\vspace{0.5cm}

\begin{equation}
\begin{array}{ll}
\hspace{-1.65cm}\hbox{(d)} & \hspace{-1.3cm}\hbox{${\frak R}_{b,s,t}:=\left\{(z,w)\in\CC^2: s|w|^b<\hbox{Re}\,z<t|w|^b,\,w\ne 0\right\}$,}\\
&\hspace{-0.3cm}\hbox{$0\le s<t\le\infty$, where $s=0$ and $t=\infty$ do not hold}\\
&\hspace{-0.3cm}\hbox{simultaneously};
\end{array}\label{frakrbst}
\end{equation}

$$
\begin{array}{ll}
\hspace{-2.3cm}\hbox{(e)} & \hspace{-1.9cm} \hbox{${\frak R}_{s,t}:=\left\{(z,w)\in\CC^2: s\exp\left(\hbox{Re}\,w\right)<\hbox{Re}\,z<t\exp\left(\hbox{Re}\,w\right)\right\}$,}\\ 
&\hspace{-0.6cm}\hbox{$0\le s<t\le\infty$, where $s=0$ and $t=\infty$ do not hold}\\
&\hspace{-0.6cm}\hbox{simultaneously};
\end{array}
$$ 

\noindent (f) $\left\{(z,w)\in\CC^2: s\exp\left({|z|^2}\right)<|w|<t\exp\left({|z|^2}\right)\right\}$, $0<s<t\le\infty$;
\vspace{0.5cm}

\noindent (g) $R_{b,s,t}$, $0\le s<t\le \infty$, where $s=0$ and $t=\infty$ do not hold simultaneously\\
\indent (see (\ref{usualrbst}));
\vspace{0.5cm}

\noindent (h) $U_{s,t}$, $0\le s<t\le\infty$, where $s=0$ and $t=\infty$ do not hold simultaneously\\ 
\indent (see (\ref{usualubst}));
\vspace{0.5cm}

\noindent (j) ${\frak S}_{s,t}$, $0\le s<t<\infty$ (see (\ref{frakst}));
\vspace{0.5cm}

\noindent (j') ${\frak S}_{s,t}^{(\infty)}$, $0\le s<t<\infty$ (see (\ref{frakstinftyn}));
\vspace{0.5cm}

\noindent (j'') ${\frak S}_{s,t}^{(n)}$, $0\le s<t<\infty$, $n\ge 2$ (see (\ref{frakstinftyn}));
\vspace{0.5cm}

\noindent (k) $V_{b,s,t}$, $0<t<\infty$, $e^{-2\pi b}t<s<t$ (see (\ref{usualvbst}));
\vspace{0.5cm}

\noindent (l) $E_{s,t}$, $1\le s<t<\infty$ (see (\ref{est}));
\vspace{0.5cm}

\noindent (l') $E_{s,t}^{(4)}$, $1\le s<t<\infty$ (see (\ref{est24}));
\vspace{0.5cm}

\noindent (l'') $E_{s,t}^{(2)}$, $1\le s<t<\infty$ (see (\ref{est24}));
\vspace{0.5cm}

\noindent (m) $\Omega_{s,t}$, $-1\le s<t\le 1$ (see (\ref{omegast}));
\vspace{0.5cm}

\noindent (m') $\Omega_{s,t}^{(\infty)}$, $-1\le s<t\le 1$ (see (\ref{omegastinftyn}));
\vspace{0.5cm}

\noindent (m'') $\Omega_{s,t}^{(n)}$, $-1\le s<t\le 1$, $n\ge 2$ (see (\ref{omegastinftyn}));
\vspace{0.5cm}

\noindent (n) $D_{s,t}$, $1\le s<t\le\infty$ (see (\ref{dst}));
\vspace{0.5cm}

\noindent (n') $D_{s,t}^{(\infty)}$, $1\le s<t\le\infty$ (see (\ref{dstinftyn}));
\vspace{0.5cm}

\noindent (n'') $D_{s,t}^{(2)}$, $1\le s<t\le\infty$ (see (\ref{dstinftyn}));
\vspace{0.5cm}

\noindent (n'') $D_{s,t}^{(2n)}$, $1\le s<t\le\infty$, $n\ge 2$ (see (\ref{dstinftyn}));
\vspace{0.5cm}

\noindent (n'') $D_{s,t}^{(n)}$, $1\le s<t\le\infty$, $n\ge 3$, $n$ is odd (see {\bf (10)}).
\vspace{0.5cm}

\noindent This concludes our orbit gluing procedure. Note that in each of cases (d) and (e) we have two non-isomorphic possibilities for $R_{\frak m}$. Each of the possibilities leads to the same set of $R_{\frak m}$-invariant domains.   

We now observe that the automorphism groups of all $R_{\frak m}$-invariant domains that appear in cases (b)--(f) have dimension at least 4. Finally, excluding equivalent domains leads to list (i)--(xv) as stated in the theorem.

The proof of the theorem is complete.\qed

\section{Levi-Flat orbits}\label{sectleviflatorbits}
\setcounter{equation}{0}

In this section we give a classification of (2,3)-manifolds $M$ for which every $G(M)$-orbit has codimension 1 and at least one orbit is Levi-flat. We start by classifying all possible Levi-flat orbits up to $CR$-diffeomorphisms together with group actions. 

Observe, first of all, that the Levi-flat hypersurface ${\cal O}_1$ (see (\ref{calo1})) is an orbit of the action of each of the groups $G_b$ (see (\ref{auttau})) for $b\in\RR$ (including $b=0$) and ${\frak G}$ (see (\ref{autxi})) on $\CC^2$. Recall next that the Levi-flat hypersurface ${\cal O}_0^{(n)}$ (see (\ref{calo61}), (\ref{calo6n})) is an orbit of the action of ${\cal R}^{(n)}$ (see {\bf (10)}, (\ref{autnun1})) on ${\cal Q}_{-}$ for $n=1$ (see (\ref{calqminus})) and on $M^{(n)}$ for $n\in\NN$, $n\ge 2$ (see {\bf (11)(d)}). Furthermore, the Levi-flat hypersurface ${\cal O}_0^{(\infty)}$ (see (\ref{calo6infty})) is an orbit of the action of ${\cal R}^{(\infty)}$ (see (\ref{autnuinfty1})) on $M^{(\infty)}$ (see {\bf (11)(d)}). Note also that the Levi-flat hypersurface
\begin{equation}
{\cal O}_1':=\left\{(z,w)\in\CC^2:\hbox{Re}\,z>0,\,|w|=1\right\}\label{calo1primeprime}
\end{equation}
is an orbit of the action on $\CC^2$ of the group $G_0'$ of all maps     
\begin{equation}
\begin{array}{lll}
z & \mapsto & \lambda z+i\beta,\\
w & \mapsto & e^{i\psi}w,
\end{array}\label{groupcircle}
\end{equation}
where $\lambda>0$, $\beta,\psi\in\RR$. The hypersurface ${\cal O}_0^{(\infty)}$ is $CR$-equivalent to ${\cal O}_1$, and the hypersurface ${\cal O}_0^{(n)}$ is $CR$-equivalent to ${\cal O}_1'$ for every $n\in\NN$.

We will now prove the following proposition. Note that it applies to  (2,3)-manifolds possibly containing codimension 2 orbits. 

\begin{proposition}\label{propleviflat}\sl Let $M$ be a (2,3)-manifold. Assume that for a point $p\in M$ its orbit $O(p)$ is Levi-flat. Then one of the following holds:
\smallskip\\

\noindent (i) $O(p)$ is equivalent to ${\cal O}_1$ by means of a real-analytic $CR$-map that transforms $G(M)|_{O(p)}$ into either the group $G_b|_{{\cal O}_1}$ for some $b\in\RR$ or the group ${\frak G}|_{{\cal O}_1}$;

\noindent (ii) $O(p)$ is equivalent to ${\cal O}_1'$ by means of a real-analytic $CR$-map that transforms $G(M)|_{O(p)}$ into the group $G_0'|_{{\cal O}_1'}$;

\noindent (iii) $O(p)$ is equivalent to ${\cal O}_0^{(j)}$ for some $j\in\{1,2,\dots,\infty\}$ by means of a real-analytic $CR$-map that transforms $G(M)|_{O(p)}$ into the group ${\cal R}^{(j)}|_{{\cal O}_0^{(j)}}$.
\end{proposition}

\noindent {\bf Proof:} Recall that the hypersurface $O(p)$ is foliated by complex manifolds equivalent to $\Delta$ (see (ii) of Proposition \ref{dim}). For convenience, we realize $\Delta$ as the right half-plane ${\cal P}:=\left\{z\in\CC:\hbox{Re}\,z>0\right\}$. Denote by ${\frak g}(M)$ the Lie algebra of vector fields on $M$ arising from the action of $G(M)$. The algebra ${\frak g}(M)$ is isomorphic to the Lie algebra of $G(M)$.  We identify every vector field from ${\frak g}(M)$ with its restriction to $O(p)$. For $q\in O(p)$ we consider the leaf $M_q$ of the foliation passing through $q$ and the subspace ${\frak l}_q\subset{\frak g}(M)$ of all vector fields tangent to $M_q$ at $q$. Since vector fields in ${\frak l}_q$ remain tangent to $M_q$ at each point in $M_q$, the subspace ${\frak l}_q$ is in fact a Lie subalgebra of ${\frak g}(M)$. It follows from the definition of ${\frak l}_q$ that $\hbox{dim}\,{\frak l}_q=2$.

Denote by $H_q$ the (possibly non-closed) connected subgroup of $G(M)$ with Lie algebra ${\frak l}_q$. It is straightforward to verify that the group $H_q$ acts on $M_q$ by holomorphic transformations. If some element $g\in H_q$ acts trivially on $M_q$, then $g\in I_q$. If for every non-identical element of $L_q$ its projection to $L_q'$ is non-identical (see (ii) of Proposition \ref{dim}), then every non-identical element of $I_q$ acts non-trivially on $M_q$ and thus $g=\hbox{id}$; if $L_q$ contains a non-identical element with an identical projection to $L_q'$ and $g\ne\hbox{id}$, then $g=g_q$, where $g_q$ denotes the element of $I_q$ corresponding to the non-trivial element in $\ZZ_2$ (see (ii) of Proposition \ref{dim}). Thus, $\hbox{dim}\,H_q=2$, and either $H_q$ or $H_q/\ZZ_2$ acts effectively on $M_q$ (the former case occurs if $g_q\not\in H_q$, the latter if $g_q\in H_q$). As we noted in the proof of Lemma \ref{groupsofmodels}, every 2-dimensional (a priori not necessarily closed) subgroup of $\hbox{Aut}({\cal P})$ is conjugate in $\hbox{Aut}({\cal P})$ to the subgroup ${\cal T}$ (see (\ref{groupt})). The Lie algebra of this subgroup is isomorphic to the 2-dimensional Lie algebra ${\frak h}$ given by two generators $X$ and $Y$ satisfying $[X,Y]=X$. Therefore, ${\frak l}_q$ is isomorphic to ${\frak h}$ for every $q\in O(p)$.

It is straightforward to determine all 3-dimensional Lie algebras containing a subalgebra isomorphic to ${\frak h}$. Every such algebra has generators $X,Y,Z$ that satisfy one of the following sets of relations:
\begin{equation}
\begin{array}{lllll}
\hbox{(R1)} &[X,Y]=X,& [Z,X]=0, & [Z,Y]=b Z, & b\in\RR,\\
\hbox{(R2)} &[X,Y]=X,& [Z,X]=0,& [Z,Y]=X+Z, & \\
\hbox{(R3)} & [X,Y]=X,& [Z,X]=Y, & [Z,Y]=-Z.& 
\end{array}\label{algebraforms}
\end{equation}

Suppose first that ${\frak g}(M)$ is given by relations (R1). In this case ${\frak g}(M)$ is isomorphic to the Lie algebra of the simply-connected Lie group $G_b$ (see (\ref{auttau})). Indeed, the Lie algebra of $G_b$ is isomorphic to the Lie algebra of vector fields on $\CC^2$ with the generators 
$$
\begin{array}{lll}
X_1 & := & \displaystyle i\,\partial/\partial z,\\
Y_1 &: = & z\,\displaystyle \partial/\partial z+bw\,\partial/\partial w,\\
Z_1 & := & \displaystyle i\, \partial/\partial w,
\end{array}
$$
that clearly satisfy (R1).

Assume first that $b\ne 0$. In this case the center of $G_b$ is trivial, and hence $G_b$ is the only (up to isomorphism) connected Lie group whose Lie algebra is given by relations (R1). Therefore, $G(M)$ is isomorphic to $G_b$. Assume further that $b\ne 1$. In this case, it is straightforward to observe that every subalgebra of ${\frak g}(M)$ isomorphic to ${\frak h}$ is generated either by $X_1$ and $Y_1+\nu Z_1$, or by $Z_1$ and $\nu X_1+Y_1$ for some $\nu\in\RR$. The connected subgroup of $G_b$ with Lie algebra generated by $X_1$ and $Y_1+\nu Z_1$ is conjugate in $G_b$ to the closed subgroup $H^1_b$ given by $\gamma=0$ in (\ref{auttau}); similarly, the connected subgroup of $G_b$ with Lie algebra generated by $Z_1$ and $\nu X_1+Y_1$ is conjugate  to the closed subgroup $H^2_b$ given by $\beta=0$ in (\ref{auttau}). Moreover, the conjugating element can be chosen to belong to the subgroup ${\cal W}^1$ of maps of the form
\begin{equation}
\begin{array}{lll}
z & \mapsto & z\\
w & \mapsto & w+i\gamma,\quad \gamma\in\RR,
\end{array}\label{formmm1}
\end{equation}
in the first case, and to the subgroup ${\cal W}^2$ of maps of the form
$$
\begin{array}{lll}
z & \mapsto & z+i\beta,\quad \beta\in\RR,\\
w & \mapsto & w,
\end{array}
$$
in the second case. These subgroups are one-parameter subgroups of $G_b$ arising from $Z_1$ and $X_1$, respectively.

Thus, upon identifying $G(M)$ with $G_b$, the subgroup $H_q$ for every $q\in O(p)$ is conjugate to either $H^1_b$ or $H^2_b$ by an an element of either ${\cal W}^1$ or ${\cal W}^2$, respectively. In particular, $H_q$ is isomorphic to ${\cal T}$ and hence does not have subgroups isomorphic to $\ZZ_2$. Therefore, $H_q$ acts effectively on $M_q$. Since the subgroups $H_q$ are conjugate to each other, it follows that either $H_q$ is conjugate to $H^1_b$ for every $q$, or $H_q$ is conjugate to $H^2_b$ for every $q$. Suppose first that the former holds. Then for every $q\in O(p)$ every element of $G(M)$ can be written as $gh$, where $g\in{\cal W}^1$, $h\in H_q$. Hence for every $q_1,q_2\in O(p)$ there exists $g\in{\cal W}^1$ such that $gM_{q_1}=M_{q_2}$. Furthermore, since the normalizer of $H^1_b$ in $G_b$ coincides with $H^1_b$, such an element $g$ is unique. Let $q_0\in O(p)$ be a point for which $H_{q_0}=H_b^1$, and let $f:M_{q_0}\ra {\cal P}$ be a holomorphic equivalence that transforms $H_{q_0}|_{M_{q_0}}$ into the group ${\cal T}$. Let $\hat X_1$ and $\hat Y_1$ be the vector fields on $O(p)$ corresponding to $X_1,Y_1$. Under the map $f$ the vector fields $\hat X_1|_{M_{q_0}}$ and $\hat Y_1|_{M_{q_0}}$ (which are tangent to $M_{q_0}$) transform into some vector fields $X_1^*$ and $Y_1^*$ on ${\cal P}$ such that $[X_1^*,Y_1^*]=X_1^*$. Clearly, $X_1^*$ and $Y_1^*$ generate the algebra of vector fields on ${\cal P}$ arising from the action of ${\cal T}$. It is straightforward to verify that one can find an element of ${\cal T}$ that transforms $Y_1^*$ into $z\,\partial/\partial z=Y_1|_{\cal P}$ and $X_1^*$ into one of $\pm i\,\partial/\partial z=\pm X_1|_{\cal P}$, and therefore we can assume that $f$ is chosen so that it transforms $\hat X_1|_{M_{q_0}}$ and $\hat Y_1|_{M_{q_0}}$ into $\pm X_1|_{\cal P}$, $Y_1|_{\cal P}$, respectively. 

For every $q\in O(p)$ we now find the unique element $g\in{\cal W}^1$ such that $gM_{q_0}=M_q$ and define $F(q):=\Bigl(f(g^{-1}(q)),i\gamma\Bigr)\in\CC^2$, with $\gamma$ corresponding to $g$ as in formula (\ref{formmm1}). Clearly, $F$ is a real-analytic $CR$-isomorphism between $O(p)$ and ${\cal O}_1$ that transforms $\hat Z_1$ into $i\,\partial/\partial w|_{{\cal O}_1}=Z_1|_{{\cal O}_1}$, where $\hat Z_1$ is the vector field on $O(p)$ corresponding to $Z_1$ (recall that ${\cal W}^1=\{\exp(sZ_1), s\in\RR\}$).

Denote by $\tilde X$, $\tilde Y$ the vector fields on ${\cal O}_1$ into which $F$ transforms $\hat X_1$, $\hat Y_1$, respectively. Since $F$ is real-analytic, it extends to a biholomorphic map from a neighborhood of $O(p)$ in $M$ onto a neighborhood of ${\cal O}_1$ in $\CC^2$. Clearly, $\hat X_1$, $\hat Y_1$, extend to holomorphic vector fields on all of $M$ and hence $\tilde X$, $\tilde Y$, extend to holomorphic vector fields defined in a neighborhood of ${\cal O}_1$. Since the restrictions of $\tilde X$ and $\tilde Y$ to ${\cal P}\times\{0\}\subset {\cal O}_1$ are $\pm\, i\partial/\partial z$ and $z\,\partial/\partial z$, respectively, these vector fields have the forms
\begin{equation}
\tilde X=(\pm i+\rho(z,w))\partial/\partial z+\sigma(z,w)\partial/\partial w,\label{vectorspecform1},
\end{equation}
and
\begin{equation}
\tilde Y=(z+\mu(z,w))\partial/\partial z+\tau(z,w)\partial/\partial w,
\label{vectorspecform2}
\end{equation}
where $\rho,\sigma,\mu,\tau$ are functions holomorphic near ${\cal O}_1$ and such that
\begin{equation}
\rho(z,0)\equiv\sigma(z,0)\equiv\mu(z,0)\equiv\tau(z,0)\equiv 0.\label{zo}
\end{equation}
Since $[\hat Z_1,\hat X_1]=0$ and $[\hat Z_1,\hat Y_1]=b\hat Z_1$ on $O(p)$, on a neighborhood of ${\cal O}_1$ we obtain
\begin{equation}
[Z_1,\tilde X]=0,\,[Z_1,\tilde Y] =bZ_1.\label{commut1}
\end{equation}
Conditions (\ref{zo}), (\ref{commut1}) imply: $\rho\equiv 0$, $\sigma\equiv 0$, $\mu\equiv 0$, $\tau=bw$. Thus, $\tilde X=\pm X_1$, $\tilde Y=Y_1$, and hence $F$ transforms $G(M)|_{O(p)}$ into $G_b|_{{\cal O}_1}$. 

The case when $H_q$ is conjugate to $H^2_b$ for every $q\in O(p)$ is treated similarly; arguing as above we construct a real-analytic $CR$-isomorphism between $O(p)$ and $\hat{\cal O}_1$ (see (\ref{calo1prime})) that transforms $G(M)|_{O(p)}$ into $G_b|_{\hat{\cal O}_1}$. Further, interchanging the variables turns $\hat{\cal O}_1$ into ${\cal O}_1$ and $G_b$ into $G_{1/b}$. 

Suppose now that $b=1$. In this case, in addition to the subalgebras arising for $b\ne 1$, a subalgebra of ${\frak g}(M)$ isomorphic to ${\frak h}$ can also be generated by $X_1+\eta Z_1$ and $Y_1+\nu Z_1$ for some $\eta,\nu\in\RR$, $\eta\ne 0$. The connected subgroup of $G_1$ corresponding to this subalgebra is conjugate in $G_1$ to the closed subgroup $H_{1,\eta}$ of all maps of the form (\ref{auttau}) with $b=1$, $\gamma=\beta\eta$. Moreover, the conjugating element can be chosen to belong to the subgroup ${\cal W}^1$ (see (\ref{formmm1})). Thus, upon identifying $G(M)$ with $G_1$, the subgroup $H_q$ for every $q\in O(p)$ is conjugate to either $H^1_1$ or $H^2_1$, or $H_{1,\eta}$ for some $\eta\ne 0$ (all these subgroups are closed). In particular, $H_q$ is isomorphic to ${\cal T}$ and hence acts effectively on $M_q$. Since the subgroups $H_q$ are conjugate to each other, it follows that either $H_q$ is conjugate to $H^1_1$ for every $q$, or $H_q$ is conjugate to $H^2_1$ for every $q$, or $H_q$ is conjugate to $H_{1,\eta}$ for every $q$ and a fixed $\eta$. The first two cases are treated as for $b\ne 1$. Suppose that $H_q$ is conjugate to $H_{1,\eta}$ for every $q\in O(p)$. It can be shown, as before, that for every $q_1,q_2\in O(p)$ there exists a unique $g\in {\cal W}^1$ such that $gM_{q_1}=M_{q_2}$. Fix $q_0\in O(p)$ with the property $H_{q_0}=H_{1,\eta}$, and let $f:M_{q_0}\ra{\cal P}$, with be a holomorphic equivalence that transforms $H_{q_0}|_{M_{q_0}}$ into the group ${\cal T}$ and such that $\hat X_1+\eta\hat Z_1|_{M_{q_0}}$ and $\hat Y_1|_{M_{q_0}}$ (which are tangent to $M_{q_0}$) are transformed into the vector fields $\pm X_1|_{{\cal P}}$ and $Y_1|_{{\cal P}}$, respectively. For every $q\in O(p)$ we now find the unique map $g\in{\cal W}^1$ such that $gM_{q_0}=M_q$ and define $F(q):=\Bigl(f(g^{-1}(q)),i\gamma\Bigr)$, with $\gamma$ corresponding to $g$ as in formula (\ref{formmm1}). Analogously to the case $b\ne 1$ we obtain: $\tilde X=\pm X_1-\eta Z_1$, $\tilde Y=Y_1$. Hence $F$ transforms $G(M)|_{O(p)}$ into $G_1|_{{\cal O}_1}$. 

Suppose now that $b=0$. In this case there are exactly two (up to isomorphism) connected Lie groups with Lie algebra ${\frak g}(M)$: $G_0$ and $G_0'$ (see (\ref{groupcircle})). It is straightforward to see that every subalgebra of ${\frak g}(M)$ isomorphic to ${\frak h}$ is generated by $X_1$ and $Y_1+\nu Z_1$ for some $\nu\in\RR$. Clearly, the connected subgroup of $G_0$ with Lie algebra generated by $X_1$ and $Y_1+\nu Z_1$ coincides with the closed normal subgroup $H_{0,\nu}$ given by $\lambda=e^t$, $\gamma=\nu t$, $t\in\RR$ (see (\ref{auttau})).
It then follows that if $G(M)$ is isomorphic to $G_0$, there exists $\nu\in\RR$, such that, identifying $G(M)$ and $G_0$, we have $H_q=H_{0,\nu}$ for every $q\in O(p)$.
Further, let us realize the Lie algebra of $G_0'$ as the Lie algebra generated by the following vector fields on $\CC^2$: $X_1$, $Y_1$, $Z_1':=iw\, \partial/\partial w$, which clearly satisfy (R1) of (\ref{algebraforms}). The connected subgroup of $G_0'$ with Lie algebra generated by $X_1$ and $Y_1+\nu Z_1'$ coincides with the closed normal subgroup $H_{0,\nu}'$ of $G_0'$ given by $\lambda=e^t$, $\psi=\nu t$, $t\in\RR$ (see (\ref{groupcircle})). It then follows that if $G(M)$ is isomorphic to $G_0'$, there exists $\nu\in\RR$, such that, identifying $G(M)$ and $G_0'$, we have $H_q=H_{0,\nu}'$ for every $q\in O(p)$. 

Thus, if $b=0$, every subgroup $H_q$ is normal, closed, isomorphic to ${\cal T}$ (hence acts effectively on $M_q$). In particular, all these subgroups coincide for $q\in O(p)$. Denote by $H$ the coinciding subgroups $H_q$. The group $H$ acts properly on $O(p)$, and the orbits of this action are the leaves $M_q$ of the foliation on $O(p)$. Further, we have $G(M)=H\times L$, where $L$ is either the subgroup 
${\cal W}^1$ (see (\ref{formmm1})), or the subgroup ${\cal W}^{'1}$ given by $\lambda=1$, $\beta=0$ in formula (\ref{groupcircle}), and hence is isomorphic to either $\RR$ or $S^1$. For every $q\in O(p)$ let $S_q:=\left\{g\in L: gM_q=M_q\right\}$. Since $M_q$ is closed, $S_q$ is a closed subgroup of $L$. Clearly, for every $g\in S_q$ there is $h\in H$ such that $hg\in I_q$. The elements $g$ and $h$ lie in the projections of $I_q$ to $L$ and $H$, respectively. Since $H$ is isomorphic to ${\cal T}$, it does not have non-trivial finite subgroups, hence the projection of $I_q$ to $H$ is trivial, and therefore $S_q=I_q$. Since all isotropy subgroups are contained in the Abelian subgroup $L$ and are conjugate to each other in $G(M)$, they are in fact identical. The effectiveness of the action of $G(M)$ on $M$ now implies that all isotropy subgroups are trivial and hence every $S_q$ is trivial as well. 

Thus, we have shown that for every $q_1,q_2\in O(p)$ there is a unique $g\in L$, such that $gM_{q_1}=M_{q_2}$. Suppose first that $L={\cal W}^1$. Fix $q_0\in O(p)$, and let $f:M_{q_0}\ra {\cal P}$ be a holomorphic equivalence that transforms $H|_{M_{q_0}}$ into the group ${\cal T}$ and $\hat X_1|_{M_{q_0}}$, $\hat Y_1+\nu \hat Z_1|_{M_{q_0}}$ into $\pm X_1|_{\cal P}$, $Y_1|_{\cal P}$, respectively. For every $q\in O(p)$ find the unique map $g\in {\cal W}^1$ such that $gM_{q_0}=M_q$ and define $F(q):=\Bigl(f(g^{-1}(q)),i\gamma\Bigr)$, with $\gamma$ corresponding to $g$ as in formula (\ref{formmm1}). It can now be shown as in the case $b\ne 0$ that $F$ transforms $G(M)|_{O(p)}$ into $G_0|_{{\cal O}_1}$.

Suppose now that $L={\cal W}^{'1}$. Fix $q_0\in O(p)$, and let $f:M_{q_0}\ra {\cal P}$ be a holomorphic equivalence that transforms $H|_{M_{q_0}}$ into the group ${\cal T}$ and $\hat X_1|_{M_{q_0}}$, $\hat Y_1+\nu \hat Z_1'|_{M_{q_0}}$ into $\pm X_1|_{\cal P}$, $Y_1|_{\cal P}$, respectively, where $\hat Z_1'$ denotes the vector field on $O(p)$ corresponding to $Z_1'$. For every $q\in O(p)$ find the unique map $g\in {\cal W}^{'1}$ such that $gM_{q_0}=M_q$ and define $F(q):=\Bigl(f(g^{-1}(q)),e^{i\psi}\Bigr)$, with $e^{i\psi}$ corresponding to $g$ as in formula (\ref{groupcircle}). Clearly, $F$ is a real-analytic $CR$-isomorphism between $O(p)$ and ${\cal O}_1'$ (see (\ref{calo1primeprime}))  that transforms $\hat Z_1'$ into $iw\,\partial/\partial w|_{{\cal O}_1'}=Z_1'|_{{\cal O}_1'}$.

As before, denote by $\tilde X$, $\tilde Y$ the vector fields on ${\cal O}_1$ into which $F$ transforms $\hat X_1$, $\hat Y_1$, respectively. These vector fields extend to holomorphic vector fields defined in a neighborhood of ${\cal O}_1'$. Since the restrictions of $\tilde X$ and $\tilde Y+\nu Z_1'$ to ${\cal P}\times\{1\}\subset {\cal O}_1'$ are $\pm X_1|_{{\cal P}}$ and $Y_1|_{{\cal P}}$, these vector fields have the forms that appear in the right-hand sides of formulas (\ref{vectorspecform1}), (\ref{vectorspecform2}), respectively, where $\rho,\sigma,\mu,\tau$ are functions holomorphic near ${\cal O}_1'$ and such that
\begin{equation}
\rho(z,1)\equiv\sigma(z,1)\equiv\mu(z,1)\equiv\tau(z,1)\equiv 0.\label{z1}
\end{equation}
Since $[\hat Z_1',\hat X_1]=[\hat Z_1',\hat Y_1+\nu \hat Z_1']=0$ on $O(p)$, on a neighborhood of ${\cal O}_1'$ we obtain
\begin{equation}
[Z_1',\tilde X]=[Z_1',\tilde Y+\nu Z_1']=0.\label{commut2}
\end{equation}
Conditions (\ref{z1}), (\ref{commut2}) imply: $\rho\equiv\sigma\equiv\mu\equiv\tau\equiv0$. Thus, $\tilde X=\pm X_1$, $\tilde Y=Y_1-\nu Z_1'$, and hence $F$ transforms $G(M)|_{O(p)}$ into $G_0'|_{{\cal O}_1'}$.

Suppose next that ${\frak g}(M)$ is given by relations (R2) (see (\ref{algebraforms})). In this case ${\frak g}(M)$ is isomorphic to the Lie algebra of the simply-connected Lie group ${\frak G}$ (see (\ref{autxi})). Indeed, the Lie algebra of ${\frak G}$ is isomorphic to the Lie algebra of holomorphic vector fields on $\CC^2$ with the following generators: $X_1$, $Y_2:=(z+w)\,\partial/\partial z+ w\,\displaystyle \partial/\partial w$, $Z_1$, which clearly satisfy (R2). It is straightforward to observe that the center of ${\frak G}$ is trivial, and hence ${\frak G}$ is the only (up to isomorphism) connected Lie group whose Lie algebra is given by relations (R2). Therefore, $G(M)$ is isomorphic to ${\frak G}$. In this case every subalgebra of ${\frak g}(M)$ isomorphic to ${\frak h}$ is generated by $X_1$ and $Y_2+\nu Z_1$ for some $\nu\in\RR$. The connected subgroup of ${\frak G}$ with Lie algebra generated by $X_1$ and $Y_2+\nu Z_1$ is conjugate in ${\frak G}$ to the closed subgroup $Q$ given by $\gamma=0$ (see (\ref{autxi})). Moreover, the conjugating element can be chosen to belong to ${\cal W}^1$ (see (\ref{formmm1})). 

Thus -- upon identification of $G(M)$ and ${\frak G}$ -- the subgroup $H_q$ for every $q\in O(p)$ is conjugate to $Q$ by an element of ${\cal W}^1$. Further, since the normalizer of $Q$ in ${\frak G}$ coincides with $Q$, we proceed as in the case of the group $G_b$ for $b\ne 0$ and obtain that there exists a real-analytic $CR$-isomorphism $F$ between $O(p)$ and ${\cal O}_1$ that transforms $\hat Z_1$ into $Z_1|_{{\cal O}_1}$ and the corresponding vector fields $\hat X_1$, $\hat Y_2$ on a neighborhood of $O(p)$ in $M$ into holomorphic vector fields $\tilde X$, $\tilde Y$ of the forms appearing in the right-hand sides of formulas (\ref{vectorspecform1}), (\ref{vectorspecform2}), respectively, where $\rho,\sigma,\mu,\tau$ are functions holomorphic near ${\cal O}_1$ and satisfying (\ref{zo}). Since $[\hat Z_1,\hat X_1]=0$ and $[\hat Z_1,\hat Y_2]=\hat X_1+\hat Z_1$ on $O(p)$, on a neighborhood of ${\cal O}_1$ we obtain
\begin{equation}
[Z_1,\tilde X]=0,\,[Z_1,\tilde Y] =\tilde X+Z_1.\label{commut5}
\end{equation}
Conditions (\ref{zo}), (\ref{commut5}) imply: $\rho\equiv 0$, $\sigma\equiv 0$, $\mu\equiv\pm w$, $\tau=w$, respectively. Thus, we have either $\tilde X= X_1$, $\tilde Y=Y_2$, or $\tilde X=-X_1$, $\tilde Y=(z-w)\,\partial/\partial z+w\,\partial/\partial w$. Hence either $F$ or $S\circ F$ transforms $G(M)|_{O(p)}$ into ${\frak G}|_{{\cal O}_1}$, where $S$ is the map given by formula (\ref{extramap}).

Suppose finally that ${\frak g}(M)$ is given by relations (R3) (see (\ref{algebraforms})). In this case ${\frak g}(M)$ is isomorphic to the algebra ${\frak {so}}_{2,1}(\RR)$. All connected Lie groups with such Lie algebra are described as follows: any simply-connected group is isomorphic to the group ${\frak V}_{\infty}$, and any non simply-connected group is isomorphic to ${\frak V}_n$ with $n\in\NN$, where ${\frak V}_{\infty}$ and ${\frak V}_n$ are the Lie groups defined in Lemma  \ref{groupsofmodels}. Clearly, the set ${\cal C}:=\left\{(z,w)\in\CC^2:\hbox{Re}\,z>0\right\}$ is ${\frak V}_j$-invariant for $j\in\{1,2,\dots,\infty\}$.

Consider in ${\frak V}_j$ three one-parameter subgroups of transformations of ${\cal C}$
$$
\begin{array}{lllllll}
\hbox{for $j=\infty$:}&&&&&&\\
&z&\mapsto &\displaystyle z-\frac{i}{2}\beta, & w & \mapsto & w,\\
&z&\mapsto &\lambda z, & w & \mapsto & w+\ln_0 \lambda,\\
&z&\mapsto &\displaystyle\frac{z}{i\mu z+1},& w & \mapsto & \displaystyle w-2\ln_0(i\mu z+1),
\end{array}
$$
$$
\begin{array}{lllllll}
\hbox{for $j=n\in\NN$:}&&&&&&\\
&z&\mapsto &\displaystyle z-\frac{i}{2}\beta, & w & \mapsto & w,\\
&z&\mapsto &\lambda z, & w & \mapsto & \lambda^{1/n}w,\\
&z&\mapsto &\displaystyle\frac{z}{i\mu z+1},& w & \mapsto & \displaystyle \frac{1}{(i\mu z+1)^{2/n}}w.
\end{array}
$$
where $\lambda>0$, $\beta,\mu\in\RR$, $t^{2/n}=\exp(2/n\,\ln_0 t)$ for $t\in\CC\setminus(-\infty,0]$, and $\ln_0$ is the branch of the logarithm in $\CC\setminus(-\infty,0]$ defined by the condition $\ln_0 1=0$. The vector fields corresponding to these subgroups generate the Lie algebras of ${\frak V}_j$ for $j\in\{1,2,\dots,\infty\}$ and are as follows:
$$
\begin{array}{llll}
\hbox{for $j=\infty$:}&&&\\
&X_3&:=&\displaystyle-\frac{i}{2}\,\partial/\partial z,\\
\vspace{0cm}&&&\\
&Y_3&:=&z\,\partial/\partial z+\partial/\partial w,\\
\vspace{0cm}&&&\\
&Z_3&:=&-iz^2\,\partial/\partial z-2iz\,\partial/\partial w,
\end{array}
$$
$$
\begin{array}{llll}
\hbox{for $j=n\in\NN$:}&&&\\
&X_3&:=&\displaystyle-\frac{i}{2}\,\partial/\partial z,\\
\vspace{0cm}&&&\\
&Y_3&:=&\displaystyle z\,\partial/\partial z+\frac{w}{n}\,\partial/\partial w,\\
\vspace{0cm}&&&\\
&Z_3&:=&\displaystyle -iz^2\,\partial/\partial z-\frac{2izw}{n}\,\partial/\partial w.
\end{array}
$$
One can verify that these vector fields indeed satisfy relations (R3).

Next, it is straightforward to show that any subalgebra of ${\frak g}(M)$ isomorphic to ${\frak h}$ is generated by either $X_3+\eta Y_3-\eta^2/2\, Z_3$, $Y_3-\eta Z_3$, with $\eta\in\RR$, or by $Y_3$, $Z_3$. For every $j\in\{1,2,\dots,\infty\}$ the connected subgroups of ${\frak V}_j$ corresponding to the subalgebras generated by $X_3+\eta Y_3-\eta^2/2\, Z_3$, $Y_3-\eta Z_3$, with $\eta\in\RR$, or by $Y_3$, $Z_3$ are isomorphic to ${\cal T}$ (hence $H_q$ acts effectively on $M_q$ for every $q$) and are all closed and conjugate to each other in ${\frak V}_j$ by elements of the one-parameter subgroup of ${\frak V}_j$ arising from $X_3-1/2\,Z_3$. We denote this subgroup by ${\cal W}_j$ and describe it in more detail. Let first $j=\infty$. Transform ${\cal C}$ into $\left\{(z,w):|z|<1\right\}$ by means of the map
$$
\begin{array}{lll}
z & \mapsto & \displaystyle\frac{z-1}{z+1},\\
\vspace{0cm}&&\\
w & \mapsto & w-2\ln_0\left((z+1)/2\right).
\end{array}
$$
Then ${\cal W}_{\infty}$ is the subgroup of ${\frak V}_{\infty}$ that transforms into the group of maps
\begin{equation}
\begin{array}{lll}
z & \mapsto & e^{it}z,\\
w & \mapsto & w+it,
\end{array}\label{winfty}
\end{equation}
where $t\in\RR$. Let now $j=n$. Transform ${\cal C}$ into $\left\{(z,w):|z|<1 \right\}$ by means of the map
$$
\begin{array}{lll}
z & \mapsto & \displaystyle\frac{z-1}{z+1},\\
\vspace{0cm}&&\\
w & \mapsto & \displaystyle\left(\frac{2}{z+1}\right)^{2/n}w.
\end{array}
$$
Then ${\cal W}_n$ is the subgroup of ${\frak V}_n$ that transforms into the group of maps
\begin{equation}
\begin{array}{lll}
z & \mapsto & e^{it}z,\\
w & \mapsto & e^{ it/n}w,
\end{array}\label{wn}
\end{equation}
where $0\le t<2\pi n$.

Observe that -- upon identifying $G(M)$ with ${\frak V}_j$ for a particular value of $j$ -- for every $q\in O(p)$ every element of $G(M)$ can be written in the form $gh$ with $g\in{\cal W}_j$, $h\in H_q$, and ${\cal W}_j\cap H_q=\{\hbox{id}\}$. Since no element in a sufficiently small neighborhood of the identity in ${\cal W}_j$ lies in the normalizer of $H_q$ in $G(M)$, for every $q\in O(p)$ there exists a tubular neighborhood ${\cal U}$ of $M_q$ in $O(p)$ with the following property: for every curve $\gamma\subset{\cal U}$ transversal to the leaves $M_{q'}$ for $q'\in{\cal U}$ and every $q_1,q_2\in\gamma$, $q_1\ne q_2$, we have $H_{q_1}\ne H_{q_2}$. Further, for every two points $q_1,q_2\in O(p)$ there exists $g\in{\cal W}_j$ such that $gM_{q_1}=M_{q_2}$.  If for some $q\in O(p)$ there is a non-trivial $g\in{\cal W}_j$ such that $gM_q=M_q$, then  $gH_qg^{-1}=H_q$ and hence $g$ has the form
$$
\begin{array}{llll}
\hbox{for $j=\infty$:}&&&\\
& z & \mapsto & z,\\
& w & \mapsto & w+2\pi i k_0, \quad k_0\in\ZZ\setminus\{0\},\\
\vspace{1cm}&&\\
\hbox{for $j=n\in\NN$:}&&\\
& z & \mapsto & z,\\
& w & \mapsto & e^{2\pi i k_0/n}w, \quad k_0\in\NN,\, k_0 \le n-1.
\end{array}
$$
It then follows that $g$ lies in the centralizer of $H_{q'}$ for every $q'\in O(p)$. Let $h\in H_q$ be such that $hg\in I_q$. Every element of $I_q$ has finite order (see (ii) of Proposition \ref{dim}), which implies that each of $h$ and $g$ is of finite order. At the same time, if $G(M)={\frak V}_{\infty}$, then $g$ is clearly of infinite order; hence $gM_q\ne M_q$ for every $q\in O(p)$ and every non-trivial $g\in {\cal W}_{\infty}$. Assume now that $G(M)={\frak V}_n$ for some $n\in\NN$. Since every non-trivial element of ${\cal T}$ has infinite order, we obtain $h=\hbox{id}$ and thus $g\in I_q$. This argument can be applied to any point in $M_q$, and thus we obtain that $g$ fixes every point in $M_q$. Hence $g=g_q$, where, as before, $g_q$ denotes the element of $I_q$ corresponding to the non-trivial element in $\ZZ_2$ (see (ii) of Proposition \ref{dim}). Then if a point $q_1\not\in M_q$ is sufficiently close to $q$, the point $q_2:=gq_1$ is also
close to $q$, and we can assume that $q_1,q_2\in{\cal U}$. It follows from the explicit form of the action of the linear isotropy subgroup $L_p$ on $T_p(M)$ that $q_1\ne q_2$ and that $q_1$, $q_2$ lie on a curve transversal to every leaf in ${\cal U}$; hence $H_{q_1}\ne H_{q_2}$. At the same time, we have $H_{q_2}=g_qH_{q_1}g_{q}^{-1}=H_{q_1}$. This contradiction shows that $gM_q\ne M_q$ for every $q\in O(p)$ and every non-trivial $g\in{\cal W}_n$.

Suppose that $G(M)={\frak V}_{\infty}$. Fix $q_0\in O(p)$ for which $H_{q_0}=H_0$, where $H_0$ is the subgroup of $G(M)$ with Lie algebra generated by $X_3$, $Y_3$, and let $f:M_{q_0}\ra {\cal P}$ be a holomorphic equivalence that transforms $H_{q_0}|_{M_{q_0}}$ into the group ${\cal T}$ and such that $\hat X_3|_{M_{q_0}}$ and $\hat Y_3|_{M_{q_0}}$ are transformed into the vector fields $\pm X_1|_{\cal P}$ and $Y_1|_{\cal P}$, respectively, where $\hat X_3$, $\hat Y_3$ are the vector fields on $O(p)$ corresponding to $X_3$, $Y_3$. For every $q\in O(p)$ we now find the unique map $g\in {\cal W}_{\infty}$ such that $gM_{q_0}=M_q$ and define $F(q):=\Bigl(f(g^{-1}(q)),it\Bigr)\in\CC^2$, where $t$ is the parameter value corresponding to $g$ (see (\ref{winfty})). Clearly, $F$ is a real-analytic $CR$-isomorphism between $O(p)$ and ${\cal O}_1$ that transforms $\hat X_3-1/2\,\hat Z_3$ into $Z_1|_{{\cal O}_1}$, where $\hat Z_3$ is the vector field on $O(p)$ corresponding to $Z_3$ (recall that ${\cal W}_{\infty}=\left\{\exp\left(s\left(X_3-1/2\,Z_3\right)\right),\,s\in\RR\right\}$), and transforms $\hat X_3$, $\hat Y_3$ on a neighborhood of $O(p)$ into holomorphic vector fields $\tilde X$, $\tilde Y$ of the forms appearing in the right-hand sides of formulas (\ref{vectorspecform1}), (\ref{vectorspecform2}), respectively, where $\rho,\sigma,\mu,\tau$ are holomorphic in a neighborhood of ${\cal O}_1$ and satisfy (\ref{zo}). Since $[\hat X_3-1/2\, \hat Z_3,\hat X_3]=-1/2\, \hat Y_3$ and $[\hat X_3-1/2\, \hat Z_3,\hat Y_3]=\hat X_3+1/2\, \hat Z_3$ on $O(p)$, on a neighborhood of ${\cal O}_1$ we obtain
\begin{equation}
[Z_1,\tilde X]=-\displaystyle\frac{1}{2}\tilde Y,\, [Z_1,\tilde Y]=2\tilde X-Z_1.\label{commut8}
\end{equation}
Conditions (\ref{zo}), (\ref{commut8}) uniquely determine the functions $\rho,\sigma,\mu,\tau$ as follows:
$$
\begin{array}{ll}
\rho&=\pm\displaystyle\frac{i}{4}\Bigl((z+2)e^{\pm w}-(z-2)e^{\mp w}\Bigr)\mp i,\quad
\sigma=-\displaystyle\frac{i}{4}\Bigl(e^w+e^{-w}-2\Bigr),\\
\vspace{0cm}\\
\mu&=\displaystyle\frac{1}{2}\Bigl((z+2)e^{\pm w}+(z-2)e^{\mp w}\Bigr)-z,\quad
\tau=-\displaystyle\frac{1}{2}\Bigl(e^w-e^{-w}\Bigr).
\end{array}
$$

Thus, we have shown that if $M$ is a (2,3)-manifold such that $G(M)$ is the universal cover of $SO_{2,1}(\RR)^0$ and $O(p)$ is a Levi-flat $G(M)$-orbit in $M$, then there exists a $CR$-isomorphism from $O(p)$ onto ${\cal O}_1$ that transforms near $O(p)$ vector fields from the Lie algebra ${\frak g}(M)$ into vector fields near ${\cal O}_1$ from the Lie algebra ${\frak a}^{(\infty)}$ generated by $Z_1$ and 
$$
\begin{array}{ll}
&\displaystyle i\Bigl((z+2)e^w-(z-2)e^{-w}\Bigr)\partial/\partial z-
\displaystyle i\Bigl(e^w+e^{-w}\Bigr)\partial/\partial w,\\
\vspace{0.1cm}&\\
&\displaystyle \Bigl((z+2)e^w+(z-2)e^{-w}\Bigr)\partial/\partial z-\displaystyle \Bigl(e^w-e^{-w}\Bigr)\partial/\partial w.
\end{array}
$$
The $CR$-isomorphism is either the map $F$ constructed above or the map $S\circ F$, where $S$ is given by (\ref{extramap}).

Let $N$ be any of the (2,3)-manifolds ${\frak D}_s^{(\infty)}$, ${\frak D}_{s,t}^{(\infty)}$ introduced in {\bf (11)(d)}. The group $G(N)$ coincides with ${\cal R}^{(\infty)}$ (see (\ref{autnuinfty1})) and hence is isomorphic to the universal cover of $SO_{2,1}(\RR)^0$. Furthermore, ${\cal O}_0^{(\infty)}$ (see (\ref{calo6infty})) is a Levi-flat $G(N)$-orbit in $N$. Hence, as we have shown above, there exists a $CR$-isomorphism from ${\cal O}_0^{(\infty)}$ onto ${\cal O}_1$ that transforms ${\frak g}(N)$ near ${\cal O}_0^{(\infty)}$  into ${\frak a}^{(\infty)}$ near ${\cal O}_1$. Therefore, there exists a $CR$-isomorphism from $O(p)$ onto ${\cal O}_0^{(\infty)}$ that transforms $G(M)|_{O(p)}$ into ${\cal R}^{(\infty)}|_{{\cal O}_0^{(\infty)}}$. 

Suppose that $G(M)={\frak V}_n$ for $n\in\NN$. Fix $q_0\in O(p)$ for which $H_{q_0}=H_0$, where, as before, $H_0$ is the subgroup of $G(M)$ with Lie algebra generated by $X_3$, $Y_3$, and let $f:M_{q_0}\ra {\cal P}$ be a holomorphic equivalence that transforms $H_{q_0}|_{M_{q_0}}$ into the group ${\cal T}$ and such that $\hat X_3|_{M_{q_0}}$ and $\hat Y_3|_{M_{q_0}}$ are transformed into the vector fields $\pm X_1|_{\cal P}$ and $Y_1|_{\cal P}$, respectively. For every $q\in O(p)$ we now find the unique map $g\in{\cal W}_n$ such that $gM_{q_0}=M_q$ and define $F(q):=\Bigl(f(g^{-1}(q)),e^{it/n}\Bigr)\in\CC^2$, where $t$ is the parameter value corresponding to $g$ (see (\ref{wn})). Clearly, $F$ is a real-analytic $CR$-isomorphism between $O(p)$ and ${\cal O}_1'$ that transforms $\hat X_3-1/2\,\hat Z_3$ into $1/n\, Z_1'|_{{\cal O}_1'}$ and transforms $\hat X_3$, $\hat Y_3$ on a neighborhood of $O(p)$ into holomorphic vector fields $\tilde X$, $\tilde Y$ of the forms appearing in the right-hand sides of formulas (\ref{vectorspecform1}), (\ref{vectorspecform2}), respectively, where $\rho,\sigma,\mu,\tau$ are functions holomorphic near ${\cal O}_1'$ and satisfying (\ref{z1}). Arguing as before, we obtain
\begin{equation}
\displaystyle\left[\frac{1}{n}\, Z_1',\tilde X\right]=-\displaystyle\frac{1}{2}\tilde Y,\,
\displaystyle\left[\frac{1}{n}\, Z_1',\tilde Y\right]=\displaystyle 2\tilde X-\frac{1}{n}\, Z_1'.\label{commut4}
\end{equation} 
Conditions (\ref{z1}), (\ref{commut4}) uniquely determine the functions $\rho,\sigma,\mu,\tau$ as follows:
$$
\begin{array}{ll}
\rho&=\pm\displaystyle\frac{i}{4}\left((z+2)w^{\pm n}-(z-2)w^{\mp n}\right)\mp i,\quad 
\sigma=-\displaystyle\frac{i}{4n}\left(w^{n+1}+w^{1-n}-2w\right),\\
\vspace{0cm}&\\
\mu&=\displaystyle\frac{1}{2}\left((z+2)w^{\pm n}+(z-2)w^{\mp n}\right)-z,\quad
\tau=-\displaystyle\frac{1}{2n}\left(w^{n+1}-w^{1-n}\right).
\end{array}
$$

Thus, we have shown that if $M$ is a (2,3)-manifold such that $G(M)$ is an $n$-sheeted cover of $SO_{2,1}(\RR)^0$ and $O(p)$ is a Levi-flat $G(M)$-orbit in $M$, then there exists a $CR$-isomorphism from $O(p)$ onto ${\cal O}_1'$ that transforms near $O(p)$ vector fields from the Lie algebra ${\frak g}(M)$ into vector fields near ${\cal O}_1'$ from the Lie algebra ${\frak a}^{(n)}$ generated by $Z_1'$ and 
$$
\begin{array}{ll}
&\displaystyle i\left((z+2)w^n-(z-2)w^{-n}\right)\partial/\partial z-\displaystyle\frac{i}{n}\left(w^{n+1}+w^{1-n}\right)\partial/\partial w,\\
\vspace{0.1cm}&\\
&\displaystyle \left((z+2)w^n+(z-2)w^{-n}\right)\partial/\partial z-\displaystyle\frac{1}{n}\left(w^{n+1}-w^{1-n}\right)\partial/\partial w.
\end{array}
$$
The $CR$-isomorphism is either the map $F$ constructed above or the map $S'\circ F$, where $S'$ is given by
$$
\begin{array}{lll}
z & \mapsto & z,\\
w & \mapsto & 1/w.
\end{array}
$$

Let $N$ be any of the (2,3)-manifolds ${\frak D}_s^{(n)}$, ${\frak D}_{s,t}^{(n)}$, $\hat{\frak D}_t^{(1)}$ (here $n=1$) introduced in {\bf (11)(d)}. The group $G(N)$ coincides with ${\cal R}^{(n)}$ (see {\bf (10)}, (\ref{autnun1})) and hence is an $n$-sheeted cover of $SO_{2,1}(\RR)^0$. Furthermore, ${\cal O}_0^{(n)}$ (see (\ref{calo61}), (\ref{calo6n})) is a Levi-flat $G(N)$-orbit in $N$. Hence, as we have shown above, there exists a $CR$-isomorphism from ${\cal O}_0^{(n)}$ onto ${\cal O}_1'$ that transforms ${\frak g}(N)$ near ${\cal O}_0^{(n)}$ into ${\frak a}^{(n)}$ near ${\cal O}_1'$. Therefore, there exists a $CR$-isomorphism from $O(p)$ onto ${\cal O}_0^{(n)}$ that transforms $G(M)|_{O(p)}$ into ${\cal R}^{(n)}|_{{\cal O}_0^{(n)}}$. 

The proof of the proposition is complete. \qed

\begin{remark}\label{explicit1}\rm It is in fact possible to write down a suitable $CR$-equivalence between ${\cal O}_0^{(j)}$ for $j\in\{1,2,\dots,\infty\}$ and either ${\cal O}_1$ or ${\cal O}_1'$ explicitly. For example, let us realize ${\cal O}_0^{(1)}$ as $\Delta\times\partial\Delta\subset\CC\PP^1\times\CC\PP^1$ (see {\bf (11)(c)}). Then the map given by  
$$
\begin{array}{lll}
z&=&-2\displaystyle\frac{ZW+1}{ZW-1},\\
\vspace{0.1cm}&&\\
w&=&W
\end{array}
$$
takes $\Delta\times\partial\Delta$ onto ${\cal O}_1'$ and transforms near $\Delta\times\partial\Delta$ the Lie algebra of vector fields arising from the action of $SU_{1,1}/\{\pm\hbox{id}\}\simeq SO_{2,1}(\RR)^0$ on $\CC\PP^1\times\CC\PP^1$ into ${\frak a}^{(1)}$ near ${\cal O}_1'$  (here we set $Z_0=W_0=1$ on $\Delta\times\partial\Delta$ and denote $Z:=Z_1$, $W:=W_1$).
\end{remark}

We will now prove the following theorem that finalizes our classification of (2,3)-manifolds in the case when every orbit is a real hypersurface. In the formulation below we use the notation introduced in Section \ref{list}.  

\begin{theorem}\label{completerealhypersurface} \sl Let $M$ be a (2,3)-manifold. Assume that the $G(M)$-orbit of every point in $M$ is of codimension 1 and that at least one orbit is Levi-flat. Then $M$ is holomorphically equivalent to one of the following pairwise non-equivalent manifolds: 
$$
\begin{array}{ll}
\hbox{(i)} & \hbox{$R_{b,s,t}$, where $b\in\RR$, $b\ne 0, 1$, and either $s=-\infty$, $t=1$ or}\\
&\hbox{$s=-1$, $0<t\le\infty$, and in the latter case $t\ne 1$, if $b= 1/2$};\\
\hbox{(ii)} & \hbox{$\hat R_{b,-1,t}$, where $b>0$, $b\ne 1$, $0<t<\infty$};\\
\hbox{(iii)} & \hbox{$\hat U_{1,t}$, where $-\infty<t<0$};\\
\hbox{(iv)} & \hbox{${\frak D}_{s,t}^{(j)}$, where $j\in\{1,2,\dots,\infty\}$, $-1\le s<1<t\le\infty$,}\\
&\hbox{and $s=-1$ and $t=\infty$ do not hold simultaneously.}
\end{array}
$$
\end{theorem}

\noindent{\bf Proof:} The proof is based on Proposition \ref{propleviflat} and the orbit gluing procedure introduced in the proof of Theorem \ref{classrealhypersurfaces}.

Observe that the set ${\frak L}:=\left\{p\in M:\,\hbox{$O(p)$ is Levi-flat}\right\}$ is closed in $M$. Hence, if ${\frak L}$ is also open, then every orbit in $M$ is Levi-flat. Let $p\in{\frak L}$ and suppose first that there exists a $CR$-isomorphism $f:O(p)\ra {\cal O}_1$ that transforms $G(M)|_{O(p)}$ into the group $G_0|_{{\cal O}_1}$. The group $G_0$ acts on ${\cal C}=\left\{(z,w)\in\CC^2:\hbox{Re}\,z>0\right\}$; every orbit of this action 
has the form
$$
{\frak b}_r:=\left\{(z,w)\in{\cal C}:\hbox{Re}\,w=r\right\},
$$
for $r\in\RR$, and hence is Levi-flat. Arguing as at step (II) of the orbit gluing procedure, we extend $f$ to a biholomorphic map between a $G(M)$-invariant neighborhood $U$ of $O(p)$ and a $G_0$-invariant neighborhood of ${\cal O}_1$ in ${\cal C}$ that satisfies (\ref{equivar}) for all $g\in G(M)$ and $q\in U$, where $\varphi:G(M)\ra G_0$ is an isomorphism. Since every $G_0$-orbit in ${\cal C}$ is Levi-flat, the set ${\frak L}$ is open. The group $G_0$ is not isomorphic to any of the groups $G_b$ for $b\in\RR^*$, $G_0'$, ${\frak G}$, ${\cal R}^{(j)}$, and it follows that every orbit $O(q)$ in $M$ is $CR$-equivalent to ${\cal O}_1$ by means of a map that transforms $G(M)|_{O(q)}$ into $G_0|_{{\cal O}_1}$.

We will now further utilize the orbit gluing procedure from the proof of Theorem \ref{classrealhypersurfaces}. Our aim is to show that $M$ is holomorphically equivalent to a $G_0$-invariant domain in ${\cal C}$. First of all, we need to prove that the map $F$ arising at step (III) extends to a holomorphic automorphism of ${\cal C}$. This map establishes a $CR$-isomorphism between ${\frak b}_{r_1}$ and ${\frak b}_{r_2}$ for some $r_1,r_2\in\RR$. Clearly, $F$ has the form $F=\nu\circ g$, where $\nu$ is a real translation in $w$, and $g\in \hbox{Aut}_{CR}({\frak b}_{r_1})$. Since $F=f_1\circ f^{-1}$ and the maps $f$ and $f_1$ transform the group $G(M)|_{O(s)}$ for some $s\in M$ into the groups $G_0|_{{\frak b}_{r_1}}$ and $G_0|_{{\frak b}_{r_2}}$, respectively, the element $g$ lies in the normalizer of $G_0|_{{\frak b}_{r_1}}$ in $\hbox{Aut}_{CR}({\frak b}_{r_1})$. 

The general form of an element of $\hbox{Aut}_{CR}({\frak b}_{r_1})$ is 
\begin{equation}
(z,r_1+iv)\mapsto (a_v(z),r_1+i\mu(v)),\label{genform}
\end{equation}
where $v\in\RR$, $a_v\in\hbox{Aut}({\cal P})$ for every $v$, and $\mu$ is a diffeomorphism of $\RR$. Considering $g$ in this form, we obtain that $a_v$ for every $v\in\RR$ lies in the normalizer of ${\cal T}$ in $\hbox{Aut}({\cal P})$ (see (\ref{groupt})), and hence $a_v\in {\cal T}$ for all $v$. Moreover, we obtain: $a_{v_1} a a_{v_1}^{-1}=a_{v_2} a a_{v_2}^{-1}$ for all $a\in {\cal T}$ and all $v_1, v_2\in\RR$. Therefore, $a_{v_1}^{-1}a_{v_2}$ lies in the center of ${\cal T}$, which is trivial. Hence we obtain that $a_{v_1}=a_{v_2}$ for all $v_1,v_2$. In addition, there exists $d\in\RR^*$ such that $\mu^{-1}(v)+\gamma\equiv\mu^{-1}(v+d\gamma)$, for all $\gamma\in\RR$. Differentiating this identity with respect to $\gamma$ at $0$ gives
\begin{equation}
\mu^{-1}(v)=v/d+t_0\label{formmu}
\end{equation}
for some $t_0\in\RR$. Therefore, $F$ extends to a holomorphic automorphism of ${\cal C}$ as the following map:
\begin{equation}
\begin{array}{lll}
z & \mapsto & \lambda z+i\beta,\\
w & \mapsto & d w+\sigma-idt_0,
\end{array}\label{specialform}
\end{equation}
where $\lambda>0$, $\beta,\sigma\in\RR$.   

Any $G_0$-invariant domain in ${\cal C}$ is given by
$$
\left\{(z,w)\in\CC^2:\hbox{Re}\,z>0,\, s<\hbox{Re}\,w<t\right\},
$$
for some $-\infty\le s<t\le\infty$. At step (IV) we observe that, since ${\cal O}_1$ splits ${\cal C}$, for $V$ sufficiently small we have $V=V_1\cup V_2\cup O(x)$, where $V_j$ are open connected non-intersecting sets. If $V_j\subset D$ for $j=1,2$, then for the domain $D$ we have $s>-\infty$, $t<\infty$, and the argument applied above to the map $F$ shows that $\hat F$ has the form (\ref{specialform}). Further, using the fact that $\hat F$ is $G_0$-equivariant, we obtain that $\hat F$ is a translation in $w$ and that ${\cal C}$ covers $M$, contradicting with the hyperbolicity of $M$. It then follows that $M$ is equivalent to $\Delta^2$ which is impossible, since $d(\Delta^2)=6$.

Next, if for $p\in{\frak L}$ there exists a $CR$-isomorphism between $O(p)$ and ${\cal O}_1'$ that transforms $G(M)|_{O(p)}$ into the group $G_0'|{{\cal O}_1'}$, a similar argument gives that $M$ is holomorphically equivalent to the product of $\Delta\times A$, where $A$ is either an annulus or a punctured disk. This is impossible either since $d(\Delta\times A)=4$.

Let now $p\in{\frak L}$ and suppose that there exists a $CR$-isomorphism $f:O(p)\ra {\cal O}_1$ that transforms $G(M)|_{O(p)}$ into the group $G_1|_{{\cal O}_1}$. The group $G_1$ acts on 
$$
{\cal D}:=\CC^2\setminus\left\{(z,w)\in\CC^2:\hbox{Re}\,z=\hbox{Re}\,w=0\right\}, 
$$
with codimension 1 orbits, and, as as before, we can extend $f$ to a biholomorphic map between a $G(M)$-invariant neighborhood $U$ of $O(p)$ and a $G_1$-invariant neighborhood of ${\cal O}_1$ in ${\cal D}$ that satisfies (\ref{equivar}) for all $g\in G(M)$ and $q\in U$, where $\varphi:G(M)\ra G_1$ is an isomorphism. A $G_1$-orbit in ${\cal D}$ is either of the form
$$
\left\{(z,w)\in\CC^2:\hbox{Re}\,w=r\hbox{Re}\,z,\,\hbox{Re}\,z>0\right\},
$$
or of the form
$$
\left\{(z,w)\in\CC^2:\hbox{Re}\,w=r\hbox{Re}\,z,\,\hbox{Re}\,z<0\right\},
$$
for $r\in\RR$, or coincides with either $\hat{\cal O}_1$ (see (\ref{calo1prime})), or
\begin{equation}
\hat{\cal O}_1^{-}:=\left\{(z,w)\in\CC^2: \hbox{Re}\,z=0,\,\hbox{Re}\,w<0\right\},\label{flat2}
\end{equation}
and hence is Levi-flat. Therefore every orbit in $M$ is Levi-flat, and it  follows as before that every orbit $O(q)$ in $M$ is $CR$-equivalent to ${\cal O}_1$ by means of a map that transforms $G(M)|_{O(q)}$ into $G_1|_{{\cal O}_1}$.

In order to show that $M$ is holomorphically equivalent to a $G_1$-invariant domain in ${\cal D}$, we need to deal with steps (III) and (IV) of the orbit gluing procedure. In this case we have $F=\nu\circ g$, where $\nu$ is a map of the form (\ref{autmu41}) with $A\in GL_2(\RR)$, and $g\in \hbox{Aut}_{CR}(o)$ for some $G_1$-orbit $o$. As before, $g$ lies in the normalizer of $G_1|_o$ in $\hbox{Aut}_{CR}(o)$. Let ${\frak X}$ be a map of the form (\ref{autmu41}) with $A\in GL_2(\RR)$ that transforms $o$ into ${\cal O}_1$ and $g_{\frak X}:={\frak X}\circ g\circ{\frak X}^{-1}$. Considering $g_{\frak X}$ in the general form (\ref{genform}) with $r_1=0$ we see that for every $\lambda>0$, $\beta,\gamma\in\RR$, the composition $a_{\lambda v+\gamma}\circ a^{\lambda,\beta}\circ a^{-1}_v$, where $a^{\lambda,\beta}(z):=\lambda z+i\beta$, belongs to ${\cal T}$ and is independent of $v$. This implies that $a_v(z)=\lambda_0 z+i(C_1 v+ C_2)$ for some $\lambda_0>0$, $C_1,C_2\in\RR$. Also, for every $\lambda>0$, $\gamma\in\RR$ there exist $\lambda_1>0$, $\gamma_1\in\RR$ such that $\mu\left(\lambda\mu^{-1}(v)+\gamma\right)=\lambda_1 v+\gamma_1$. It then follows, in particular, that either there exist $c\in\RR^*$, $d\in\RR$ such that $\mu^{-1}(v)+\gamma\equiv\mu^{-1}\left(e^{c\gamma}v+d(1-e^{c\gamma})\right)$, or there exists $d\in\RR^*$ such that $\mu^{-1}(v)+\gamma\equiv\mu^{-1}\left(v+d\gamma\right)$ for all $\gamma\in\RR$. Differentiating these identities with respect to $\gamma$ at $0$, we see that the first identity cannot hold and that $\mu^{-1}$, as before, has the form (\ref{formmu}) for some $t_0\in\RR$. It then follows that $g_{\frak X}$ extends to a holomorphic automorphism of ${\cal D}$ as the map
\begin{equation}
\begin{array}{lll}
z & \mapsto & \lambda_0 z+C_1w+iC_2,\\
w & \mapsto & dw-idt_0,
\end{array}\label{formmmmggg}
\end{equation}
and thus $F$ extends to an automorphism of ${\cal D}$ as well.

Any hyperbolic $G_1$-invariant domain in ${\cal D}$ has the form $S+i\RR^2$, where $S$ is an angle of size less than $\pi$ with vertex at the origin in the $\left(\hbox{Re}\,z,\hbox{Re}\,w\right)$-plane. If at step (IV) we have $V_j\subset D$ for $j=1,2$, then the argument applied above to the map $F$ shows that $\hat F$ has the form (\ref{autchi}) with $A\in GL_2(\RR)$. Further, using the fact that $\hat F$ is $G_1$-equivariant, we obtain that $\hat F=\hbox{id}$, which is impossible. This shows that $M$ is holomorphically equivalent to a hyperbolic $G_1$-invariant domain in ${\cal D}$. By means of a suitable linear transformation every such domain is equivalent to the tube domain whose base is the first quadrant, and thus $M$ is holomorphically equivalent to $\Delta^2$, which is impossible.

Let $p\in{\frak L}$ and suppose that there exists a $CR$-isomorphism $f:O(p)\ra {\cal O}_1$ that transforms $G(M)|_{O(p)}$ into the group $G_b|_{{\cal O}_1}$ for some $b\in\RR^*$, $b\ne 1$. The group $G_b$ acts on ${\cal D}$, and every $G_b$-orbit in ${\cal D}$ is either strongly pseudoconvex and has one of the forms
$$
\left\{(z,w)\in\CC^2:\hbox{Re}\,w=r\left(\hbox{Re}\,z\right)^b,\,\hbox{Re}\,z>0\right\},
$$
$$
\left\{(z,w)\in\CC^2:\hbox{Re}\,w=r\left(-\hbox{Re}\,z\right)^b,\,\hbox{Re}\,z<0\right\},
$$
for $r\in\RR^*$, or coincides with one of ${\cal O}_1$, $\hat{\cal O}_1$ (see (\ref{calo1prime})), $\hat{\cal O}_1^{-}$ (see (\ref{flat2})), and
\begin{equation}
{\cal O}_1^{-}:=\left\{(z,w)\in\CC^2: \hbox{Re}\,z<0,\,\hbox{Re}\,w=0\right\}.\label{othreeprimes}
\end{equation}
It then follows that every Levi-flat orbit in $M$ has a $G(M)$-invariant neighborhood in which every other orbit is strongly pseudoconvex.
Among the groups $G_c$ (with $c\in\RR^*$, $c\ne 1$, $c\ne b$), ${\frak G}$, ${\cal R}^{(j)}$ the only group isomorphic to $G_b$ is $G_{1/b}$. Thus, it follows that every Levi-flat orbit $O(q)$ in $M$ is $CR$-equivalent to ${\cal O}_1$ by means of a map that transforms $G(M)|_{O(q)}$ into either $G_b|_{{\cal O}_1}$ or $G_{1/b}|_{{\cal O}_1}$. In the latter case interchanging the variables we obtain a map that takes $O(q)$ into $\hat{\cal O}_1$ and transforms $G(M)|_{O(q)}$ into $G_b|_{\hat{\cal O}_1}$. Next, by Lemma \ref{groupsofmodels}, every strongly pseudoconvex orbit $O(q')$ is $CR$-equivalent to $\tau_b$ (see (\ref{tau})) by means of a $CR$-map that transforms $G(M)|_{O(q')}$ into $G_b|_{\tau_b}$.

We now turn to step (III) of the orbit gluing procedure. For the point $x\in\partial D$ there exists a real-analytic $CR$-isomorphism $f_1$ between $O(x)$ and one of ${\cal O}_1$, $\hat{\cal O}_1$, $\tau_b$ that transforms $G(M)|_{O(x)}$ into one of $G_b|_{{\cal O}_1}$, $G_b|_{\hat{\cal O}_1}$, $G_b|{\tau_b}$, respectively. In each of these three cases the corresponding point $s$ can be chosen so that $O(s)$ is strongly pseudoconvex. Then $F$ is a $CR$-isomorphism between strongly pseudoconvex $G_b$-orbits, and thus has the form $F=\nu\circ g$, where $\nu$ is a map of the form
\begin{equation}
\begin{array}{lll}
z & \mapsto & \pm z,\\
w & \mapsto & dw
\end{array}\label{speccial}
\end{equation}
with $d\in\RR^*$, and $g\in G_b$. Therefore, $F$ extends to an automorphism of ${\cal D}$.

Suppose that at step (IV) we have $V_j\subset D$ for $j=1,2$. Assume first $O(x)$ is strongly pseudoconvex. Then $\hat F=\nu\circ g$, where $\nu$ is a non-trivial map of the form (\ref{speccial}), and $g\in G_b$. Now using the fact that $\hat F$ is $G_b$-equivariant, we obtain that $\hat F=\hbox{id}$, which is impossible. Suppose now that $O(x)$ is Levi-flat. Then $\hat F=\nu\circ g$, where $\nu$ is one of the maps
\begin{equation}
\begin{array}{lll}
z & \mapsto & -z,\\
w & \mapsto & w,
\end{array}\quad
\begin{array}{lll}
z & \mapsto & z,\\
w & \mapsto & -w,
\end{array}\quad
\begin{array}{lll}
z & \mapsto & \pm w,\\
w & \mapsto & z,
\end{array}\quad
\begin{array}{lll}
z & \mapsto & w,\\
w & \mapsto & \pm z,
\end{array}\label{fourspecforms}
\end{equation}
and $g\in\hbox{Aut}_{CR}(o)$, where $o$ is a Levi-flat $G_b$-orbit. In this case, either $g$ lies in the normalizer of $G_b|_o$ in $\hbox{Aut}_{CR}(o)$, or $gG_b|_og^{-1}=G_{1/b}|_o$ and $gG_{1/b}|_og^{-1}=G_b|_o$. Transforming $o$ into ${\cal O}_1$ by a map ${\frak X}$ from list (\ref{fourspecforms}) and arguing as in the case of the group $G_1$ for the map $F$, we obtain that $g_{\frak X}:={\frak X}\circ g\circ{\frak X}^{-1}$ extends to a holomorphic automorphism of ${\cal D}$ as a map of the form (\ref{formmmmggg}) with $C_1=0$. It then follows that $\hat F$ has the form (\ref{autchi}) with $A\in GL_2(\RR)$. Now, using the $G_b$-equivariance of $\hat F$ we again see that $\hat F=\hbox{id}$ which is impossible. Hence $M$ is holomorphically equivalent to a $G_b$-invariant domain in ${\cal D}$, and we obtain (i) and (ii) of the theorem.

Let $p\in{\frak L}$ and suppose that there exists a $CR$-isomorphism $f:O(p)\ra {\cal O}_1$ that transforms $G(M)_{O(p)}$ into the group ${\frak G}|_{{\cal O}_1}$. The group ${\frak G}$ acts on ${\cal D}$, and every ${\frak G}$-orbit in ${\cal D}$ is either strongly pseudoconvex and has one of the forms
$$
\left\{(z,w)\in\CC^2:\hbox{Re}\,z=\hbox{Re}\,w\ln\left(r\hbox{Re}\,w\right),\,\hbox{Re}\,w>0\right\},
$$
$$
\left\{(z,w)\in\CC^2:\hbox{Re}\,z=\hbox{Re}\,w\ln\left(-r\hbox{Re}\,w\right),\,\hbox{Re}\,w<0\right\},
$$
for $r>0$, or coincides with one of ${\cal O}_1$, ${\cal O}_1^{-}$ (see (\ref{othreeprimes})).

It then follows that every Levi-flat orbit in $M$ has a $G(M)$-invariant neighborhood in which every other orbit is strongly pseudoconvex, that every Levi-flat orbit in $M$ is $CR$-equivalent to ${\cal O}_1$ by means of a map that transforms $G(M)|_{O(p)}$ into ${\frak G}|_{{\cal O}_1}$ and that every strongly pseudoconvex orbit is $CR$-equivalent to $\xi$ (see (\ref{xi})) by means of a map that transforms $G(M)|_{O(p)}$ into ${\frak G}|_{\xi}$.

At step (III), as in the case of the groups $G_b$ above, we can choose $s$ so that $O(s)$ is strongly pseudoconvex. It then follows that $F=\nu\circ g$, where $\nu$ is a map of the form
\begin{equation}
\begin{array}{lll}
z & \mapsto & d z,\\
w & \mapsto & d w
\end{array}\label{specccial}
\end{equation}
with $d\in\RR^*$, and $g\in {\frak G}$. Therefore, $F$ extends to an automorphism of ${\cal D}$.  

Suppose that at step (IV) we have $V_j\subset D$ for $j=1,2$. Assume first $O(x)$ is strongly pseudoconvex. Then $\hat F=\nu\circ g$, where $\nu$ is a non-trivial map of the form (\ref{specccial}), and $g\in {\frak G}$. Now, using the fact that $\hat F$ is ${\frak G}$-equivariant, we obtain that $\hat F=\hbox{id}$, which is impossible. Suppose now that $O(x)$ is Levi-flat. Then $\hat F=\nu\circ g$, where $\nu$ is map (\ref{specccial}) with $d=-1$, and $g\in\hbox{Aut}_{CR}(o)$, where $o$ is a Levi-flat ${\frak G}$-orbit.
The element $g$ lies in the normalizer of ${\frak G}$ in $\hbox{Aut}_{CR}(o)$. Transforming $o$ into ${\cal O}_1$ by a map ${\frak X}$ of the form (\ref{specccial}) with $d=\pm 1$ and considering $g_{\frak X}:={\frak X}\circ g\circ {\frak X}^{-1}$ in the general form (\ref{genform}) with $r_1=0$, we obtain, as before, that $\mu^{-1}$ has the form (\ref{formmu}) for some $d\in\RR^*$, $t_0\in\RR$, and that $a_v(z)=\lambda_0 z+i\beta(v)$, where $\lambda_0>0$ and $\beta(v)$ is a function satisfying for every $\lambda>0$ and $\gamma\in\RR$ the following condition:
$$
\partial/\partial v\left[\beta\left(\lambda\mu^{-1}(v)+\gamma\right)-\lambda\beta\left(\mu^{-1}(v)\right)+\lambda\ln\lambda\left(\lambda_0\mu^{-1}(v)-v\right)\right]\equiv 0.
$$
Setting $\lambda=1$ in the above identity gives that $g_{\frak X}$ extends to an automorphism of ${\cal D}$ as a map of the form (\ref{formmmmggg}). Therefore, $\hat F$ has the form (\ref{autchi}) with $A\in GL_2(\RR)$, and using the ${\frak G}$-equivariance of $\hat F$ we again see that $\hat F=\hbox{id}$ which is impossible. Hence $M$ is holomorphically equivalent to a ${\frak G}$-invariant domain in ${\cal D}$, and we have obtained (iii) of the theorem.

Let $p\in{\frak L}$ and suppose that there exists a $CR$-isomorphism $f:O(p)\ra {\cal O}_0^{(j)}$ that transforms $G(M)|_{O(p)}$ into the group ${\cal R}^{(j)}|_{{\cal O}_0^{(j)}}$ for some\linebreak $j\in\{1,2,\dots,\infty\}$. The group ${\cal R}^{(j)}$ acts on ${\frak D}^{(j)}$, where ${\frak D}^{(1)}:={\frak D}_{-1,\infty}^{(1)}$ (see (\ref{frakdst1spec})), ${\frak D}^{(j)}:=M^{(j)}\setminus{\cal O}^{(2j)}$ for $1<j<\infty$, and ${\frak D}^{(\infty)}:=M^{(\infty)}$ (see {\bf (11)(a), (d)}). Apart from  ${\cal O}_0^{(j)}$, every ${\cal R}^{(j)}$-orbit in ${\frak D}^{(j)}$ is strongly pseudoconvex and is one of $\nu_{\alpha}^{(j)}$, for $-1<\alpha<1$, or $\eta_{\alpha}^{(2j)}$, for $\alpha>1$. It then follows that Levi-flat orbits in $M$ are isolated, and every such orbit $O(q)$ is $CR$-equivalent to ${\cal O}_0^{(j)}$ by means of a $CR$-isomorphism that transforms $G(M)|_{O(q)}$ into the group ${\cal R}^{(j)}|_{{\cal O}_0^{(j)}}$.

At step (III) we again choose $s$ so that $O(s)$ is strongly pseudoconvex which gives that $F$ extends to ${\frak D}^{(j)}$ as an element of ${\cal R}^{(j)}$. At step (IV), suppose that $V_j\subset D$ for $j=1,2$. Observe that $O(x)$ cannot be strongly pseudoconvex since otherwise $\hat F$ would be a $CR$-isomorphism between two distinct strongly pseudoconvex ${\cal R}^{(j)}$-orbits in ${\frak D}^{(j)}$, while in fact ${\cal R}^{(j)}$-orbits are pairwise $CR$ non-equivalent. On the other hand, $O(x)$ cannot be Levi-flat either, since otherwise $\hat F$ would be a $CR$-isomorphism between two distinct Levi-flat ${\cal R}^{(j)}$-orbits in ${\frak D}^{(j)}$, while $O_0^{(j)}$ is the only Levi-flat orbit in ${\frak D}^{(j)}$. This implies that $M$ is holomorphically equivalent to a ${\cal R}^{(j)}$-invariant domain in ${\frak D}^{(j)}$ which leads to (iv) of the theorem.

The proof of the theorem is complete. \qed

\section{Codimension 2 Orbits}\label{codimtwoorbits}
\setcounter{equation}{0}

In this section we finalize our classification by allowing codimension 2 orbits to be present in the manifold. We will prove the following theorem (as before, we use the notation introduced in Section \ref{list}).

\begin{theorem}\label{theoremcodimtwo}\sl Let $M$ be a (2,3)-manifold. Assume that a $G(M)$-orbit of codimension 1 and a $G(M)$-orbit of codimension 2 are present in $M$. Then $M$ is holomorphically equivalent to one of the following manifolds: 
$$
\begin{array}{ll}
\hbox{(i)} & \hbox{${\frak S}_1$};\\
\hbox{(ii)} & \hbox{$E_t$ with $1<t<\infty$};\\
\hbox{(iii)} & \hbox{$E_t^{(2)}$ with $1<t<\infty$};\\
\hbox{(iv)} & \hbox{$\Omega_t$ with $-1<t<1$};\\
\hbox{(v)} & \hbox{$D_s$ with $1\le s<\infty$};\\
\hbox{(vi)} & \hbox{$D_s^{(2)}$ with $1<s<\infty$};\\
\hbox{(vii)} & \hbox{$D_s^{(n)}$ with $n\ge 3$, $1\le s<\infty$};\\
\hbox{(viii)} & \hbox{${\frak D}_s^{(n)}$ with $n\ge 1$, $-1<s<1$};\\
\hbox{(ix)} & \hbox{$\hat{\frak D}_t^{(1)}$ with $1<t<\infty$}.
\end{array}
$$
\end{theorem}

\noindent {\bf Proof:} Since a codimension 1 orbit is present in $M$, it follows that there are at most two codimension 2 orbits (see \cite{AA}). Let $O$ be one such orbit. Parts (iii) and (iv) of Proposition \ref{dim} yield that for every $p\in O$ the group $I_p^0$ is isomorphic to $U_1$ (in particular, $G(M)$ has a subgroup isomorphic to $U_1$), and there exists an $I_p^0$-invariant connected complex curve $C_p$ in $M$ that intersects $O$ transversally at $p$. If $O$ is a complex curve, one such curve $C_p$ corresponds -- upon local linearization of the $I_p$-action -- to the $L_p^0$-invariant subspace $\{w=0\}$ of $T_p(M)$, where the coordinates $(z,w)$ in $T_p(M)$ are chosen with respect to the decomposition of $T_p(M)$ specified in (iii) of Proposition \ref{dim}, with $\{z=0\}$ corresponding to $V_p=T_p(O(p))$; if, in addition, the isotropy linearization is given by (\ref{stab2}), then the maximal extension of this curve is the only maximally extended complex curve in $M$ with these properties. If $O$ is a totally real orbit, $C_p$ can be constructed from any of the two $L_p^0$-invariant subspaces $\left\{z=\pm iw\right\}$ of $T_p(M)$, where the coordinates $(z,w)$ in $T_p(M)$ are chosen so that $V_p=\left\{\hbox{Im}\,z=0,\,\hbox{Im}\,w=0\right\}$ (see (iv) of Proposition \ref{dim}); locally near $p$ there are no other such curves. Clearly, there exists a neighborhood $U$ of $p$ such that $U\cap\left(C_p\setminus\{p\}\right)$ is equivalent to a punctured disk.

Since there is a codimension 1 orbit in $M$, the group $G(M)$ is either isomorphic to one of the groups listed in Lemma \ref{groupsofmodels} (if a strongly pseudoconvex orbit is present in $M$) or to one of $G_1$, $G_0$, $G_0'$ (if all codimension 1 orbits are Levi-flat -- see (\ref{auttau}) and (\ref{groupcircle})). Since $G_0$ and $G_1$ do not contain subgroups isomorphic to $U_1$, the group $G(M)$ in fact cannot be isomorphic to either of these groups. Let $M'$ be the manifold obtained from $M$ by removing all codimension 2 orbits, and suppose that $G(M)$ is isomorphic to $G_0'$. The subgroup of $G_0'$ isomorphic to $U_1$ is unique and consists of all rotations in $w$, it is normal and maximal compact in $G_0'$; we denote it by $J$. It follows from the proof of Theorem \ref{completerealhypersurface} that $M'$ is holomorphically equivalent to 
$$
{\cal V}_{s,t}:=\left\{(z,w)\in\CC^2:\hbox{Re}\,z>0,\, s<|w|<t\right\},
$$
where $0\le s<t\le\infty$, and either $s>0$ or $t<\infty$, by means of a map $f$ that satisfies (\ref{equivar}) for all $g\in G(M)$, $q\in M'$ and an isomorphism $\varphi:G(M)\ra G_0'$. Clearly, $I_p=I:=\varphi^{-1}(J)$ for every $p\in O$. In particular, $I_p$ acts trivially on $O$ for every $p\in O$; hence $O$ is a complex curve with isotropy linearization given by (\ref{stab2}), and there are no totally real orbits in $M$. The group $G_0'$ acts on $\tilde {\cal C}:={\cal P}\times\CC\PP^1$ (we set $g(z,\infty):=(\lambda z+i\beta,\infty)$ for every $g$ of the form (\ref{groupcircle})). This action has two complex curve orbits
\begin{equation}
\begin{array}{lll}
{\cal O}_7&:=&{\cal P}\times\{0\},\\
{\cal O}_8&:=&{\cal P}\times\{\infty\}.
\end{array}\label{o67}
\end{equation}
It is straightforward to observe that every connected $J$-invariant complex curve in ${\cal V}_{s,t}$ extends to a curve of the form
$$
N_{z_0}:=\{z=z_0\}\cap {\cal V}_{s,t},
$$
for some $z_0\in {\cal P}$, which is either an annulus (possibly with infinite outer radius) or a punctured disk. Fix $p_0\in O$, let $C_{p_0}$ be the unique maximally extended $I$-invariant complex curve in $M$ that intersects $O$ at $p_0$ transversally, and let $z_0\in {\cal P}$ be such that $f(C_{p_0}\setminus\{p_0\})=N_{z_0}$. Since for a sequence $\{p_j\}$ in $C_{p_0}$ converging to $p_0$ the sequence $\{f(p_j)\}$ approaches either $\{z=z_0, |w|=s\}$ or $\{z=z_0,|w|=t\}$ and $C_{p_0}\setminus\{p_0\}$ is equivalent to a punctured disk near $p_0$, we have either $s=0$ or $t=\infty$, respectively. 

Assume first that $s=0$. We extend $f$ to a map from $\hat M:=M'\cup O$ onto the domain
$$
{\cal V}_t:=\left\{(z,w)\in\CC^2:\hbox{Re}\,z>0,\, |w|<t\right\}={\cal V}_{0,t}\cup{\cal O}_7,
$$
by setting $f(p_0):=q_0:=(z_0,0)\in {\cal O}_7$, with $z_0$ constructed as above. The extended map is 1-to-1 and satisfies (\ref{equivar}) for all $g\in G(M)$, $q\in\hat M$. To prove that $f$ is holomorphic on all of $\hat M$, it suffices to show that $f$ is continuous on $O$. We will prove that every sequence $\{p_j\}$ in $\hat M$ converging to $p_0$ has a subsequence along which the values of $f$ converge to $q_0$. Let first $\{p_j\}$ be a sequence in $O$. Clearly, there exists a sequence $\{g_j\}$ in $G(M)$ such that $p_j=g_jp_0$ for all $j$. Since $G(M)$ acts properly on $M$, there exists a converging subsequence $\{g_{j_k}\}$ of $\{g_j\}$, and we denote by $g_0$ its limit. It then follows that $g_0\in I$ and, since $f$ satisfies (\ref{equivar}), we obtain that $\{f(p_{j_k})\}$ converges to $q_0$. Next, if $\{p_j\}$ is a sequence in $M'$, then there exists a sequence $\{g_j\}$ in $G(M)$ such that $g_jp_j\in C_{p_0}$. Clearly, the sequence $\{g_jp_j\}$ converges to $p_0$ and hence $\{f(g_jp_j)\}$ converges to $q_0$. Again, the properness of the $G(M)$-action on $M$ yields that there exists a converging subsequence $\{g_{j_k}\}$ of $\{g_j\}$. Let $g_0$ be its limit; as before, we have $g_0\in I$.  Property (\ref{equivar}) now implies $f(p_{j_k})=\left[\varphi(g_{j_k})\right]^{-1}f(g_{j_k}p_{j_k})$, and therefore $\{f(p_{j_k})\}$ converges to $q_0$. Thus, $f$ is holomorphic on $\hat M$. If $O'$ was another complex curve orbit, then, since $t<\infty$, arguing as above we could extend $f$ biholomorphically to a map from $M'\cup O'$ onto ${\cal V}_t$ that takes $O'$ onto ${\cal O}_7$. Then there exist non-intersecting tubular neighborhoods $U$ and $U'$ of $O$ and $O'$, respectively, such that $f(U\setminus O)=f(U'\setminus O')$, which contradicts the fact that $f$ is biholomorphic on $M'$. Hence, $O$ is the only codimension 2 orbit, and $M$ is holomorphically equivalent to ${\cal V}_t$. This is, however, impossible since $d({\cal V}_t)=6$.

Assume now that $t=\infty$. Arguing as in the case $s=0$ and mapping $O$ onto ${\cal O}_8$, we can extend $f$ to a biholomorphic map between $M$ and the domain in $\tilde {\cal C}$ given by
$$
\left\{(z,w)\in\CC^2: \hbox{Re}\,z>0,\, |w|>s\right\}\cup {\cal O}_8={\cal V}_{s,\infty}\cup{\cal O}_8,
$$
which is holomorphically equivalent to ${\cal V}_1$. This is again impossible, and we have ruled out the case when $G(M)$ is isomorphic to $G_0'$. 

It then follows that there is always a strongly pseudoconvex orbit in $M$ and hence $G(M)$ is isomorphic to one of the groups listed in Lemma \ref{groupsofmodels}. Observe that the groups that arise in subcases (g), (h), (j'), (k), (m'), (n') of case (A) as well as in cases (D) and (F) do not have non-trivial compact subgroups; thus these situations do not in fact occur. In addition, arguing as in the proof of Theorem \ref{classrealhypersurfaces}, we rule out case (B). 

We now assume that a complex curve orbit is present in $M$. Let $O$ be such an orbit. Then (iii) of Proposition \ref{dim} gives that $O$ is equivalent to ${\cal P}$. Furthermore, if for $p\in O$ the group $I_p^0$ acts on $O$ non-trivially (see (\ref{stab1})), then there exists a finite normal subgroup $H\subset I_p$ such that $G(M)/H$ is isomorphic to $\hbox{Aut}({\cal P})\simeq SO_{2,1}(\RR)^0$; if $I_p^0$ acts on $O$ trivially (see (\ref{stab2})), then there is a 1-dimensional normal compact subgroup $H\subset I_p$ such that $G(M)/H$ is isomorphic to the subgroup ${\cal T}\subset\hbox{Aut}({\cal P})$ (see (\ref{groupt})). In particular, every maximal compact subgroup of $G(M)$ is 1-dimensional and therefore is isomorphic to $U_1$. It then follows that for every $p\in O$ the group $I_p^0$ is a maximal compact subgroup of $G(M)$ and hence $I_p$ is connected. Observe now that in subcases (l), (l'), (l'') of case (A) as well as in case (C) the group $G(M)$ is compact. In case (G) the group $G(M)$ is isomorphic to $U_1\times\RR^2$; thus no factor of $G(M)$ by a finite subgroup is isomorphic to $SO_{2,1}(\RR)^0$ and the factor-group of $G(M)$ by its maximal compact subgroup is not isomorphic to ${\cal T}$.  Furthermore, in subcases (j), (j'') the group $G(M)$ is isomorphic to $U_1\ltimes\RR^2$. This group has no 1-dimensional compact normal subgroups and cannot be factored by a finite subgroup to obtain a group isomorphic to $SO_{2,1}(\RR)^0$. Therefore, if a complex curve orbit is present in $M$, we only need to consider subcases (m), (m''), (n), (n'') of case (A), and case (E).  

We start with case (E) and assume first that for some point $p\in M$ there exists a $CR$-isomorphism between $O(p)$ and the hypersurface $\varepsilon_b$ for some $b>0$ (see (d) in the proof of Theorem \ref{classrealhypersurfaces}), that transforms $G(M)|_{O(p)}$ into $G_{\varepsilon_b}|_{\varepsilon_b}$, where $G_{\varepsilon_b}$ is the group of all maps
\begin{equation}
\begin{array}{lll}
z &\mapsto& \lambda z+i\beta,\\
w &\mapsto& e^{i\psi}\lambda^{1/b}w,
\end{array}\label{gvarepsilonb}
\end{equation}
with $\lambda>0$, $\psi,\beta\in\RR$. 

We proceed as in the case of the group $G_0'$ considered above. In this case Levi-flat orbits are not present in $M$, and it follows from the proof of Theorem \ref{classrealhypersurfaces} that $M'$ is holomorphically equivalent to ${\frak R}_{b,s,t}$ (see (\ref{frakrbst})) for some $0\le s<t\le\infty$, with either $s>0$ or $t<\infty$, by means of a map $f$ that satisfies (\ref{equivar}) for all $g\in G(M)$, $q\in M'$ and an isomorphism $\varphi:G(M)\ra G_{\varepsilon_b}$. The only 1-dimensional compact subgroup of $G_{\varepsilon_b}$ is the maximal compact normal subgroup $J^{\varepsilon_b}$ given by the conditions $\lambda=1$, $\beta=0$ in (\ref{gvarepsilonb}). Clearly, $I_p=I:=\varphi^{-1}(J^{\varepsilon_b})$ for all $p\in O$, which implies, as before, that $O$ is a complex curve with isotropy linearization given by (\ref{stab2}), and there are no totally real orbits. The group $G_{\varepsilon_b}$ acts on $\tilde {\cal C}$, and, as before, this action has two complex curve orbits ${\cal O}_7$ and ${\cal O}_8$ (see (\ref{o67})).

Further, every connected $J^{\varepsilon_b}$-invariant complex curve in ${\frak R}_{b,s,t}$ extends to a curve of the form
$$
\begin{array}{ll}
\{z=z_0\}\cap {\frak R}_{b,s,t}=&\Biggl\{(z,w)\in\CC^2: z=z_0,\\
&\left(\hbox{Re}\,z_0/t\right)^{1/b}<|w|<\left(\hbox{Re}\,z_0/s\right)^{1/b}\Biggr\},
\end{array}
$$ 
for some $z_0\in {\cal P}$, which is either an annulus or a punctured disk. As before, we obtain that either $s=0$, or $t=\infty$.

If $t=\infty$, we extend $f$ to a biholomorphic map from $\hat M=M'\cup O$ onto the domain
\begin{equation}
{\frak E}_{b,s}:=\left\{(z,w)\in\CC^2: \hbox{Re}\,z>s|w|^b\right\}={\frak R}_{b,s,\infty}\cup{\cal O}_7.\label{frakebs}
\end{equation}
Since $s>0$, the orbit $O$ is the only codimension 2 orbit, and hence $M$ is holomorphically equivalent to ${\frak E}_{b,s}$.
Similarly, if $s=0$, then $M$ is holomorphically equivalent to the domain
\begin{equation}
\left\{(z,w)\in\CC^2: 0<\hbox{Re}\,z<t|w|^b\right\}\cup{\cal O}_8={\frak R}_{b,0,t}\cup{\cal O}_8,\label{prelim}
\end{equation}
which is equivalent to the domain  
$$
{\cal E}_{b,t}:=\left\{(z,w)\in\CC^2:\hbox{Re}\,z>0,\,|w|<\left(t/\hbox{Re}\,z\right)^{1/b}\right\}.
$$
This is, however, impossible since $d({\frak E}_{b,s})\ge 4$ and $d({\cal E}_{b,t})=4$.

Assume now that in case (E) for some point $p\in M$ there exists a $CR$-isomorphism $f$ between $O(p)$ and the hypersurface
$\varepsilon_b$ for some $b\in\QQ$, $b>0$, that transforms $G(M)|_{O(p)}$ into ${\frak V}_b|_{\varepsilon_b}$ (see Lemma \ref{groupsofmodels} for the definition of ${\frak V}_b$). Let $b=k_1/k_2$, for $k_1,k_2\in\NN$, with $(k_1,k_2)=1$. 

As before, $f$ extends to a biholomorphic map between $M'$ and  
${\frak R}_{b,s,t}$, where $0\le s<t\le\infty$, and either $s>0$ or $t<\infty$. The map $f$ satisfies (\ref{equivar}) for all $g\in G(M)$, $q\in M'$ and an isomorphism $\varphi:G(M)\ra {\frak V}_b$. The group ${\frak V}_b$ acts on $\tilde {\cal C}$, and, as before, this action has two complex curve orbits ${\cal O}_7$ and ${\cal O}_8$.
Every 1-dimensional compact subgroup of ${\frak V}_b$ is the isotropy subgroup of the points $(z_0,0)\in{\cal O}_7$ and $(z_0,\infty)\in{\cal O}_8$ for a uniquely chosen $z_0\in{\cal P}$. For $z_0\in{\cal P}$ denote by $J_{z_0}^{\varepsilon_b}$ the corresponding maximal compact subgroup of ${\frak V}_b$.

For every $z_0\in {\cal P}$ there is a family ${\cal F}_{z_0}^{\frak R}$ of connected closed complex curves in ${\frak R}_{b,s,t}$ invariant under the $J_{z_0}^{\varepsilon_b}$-action, such that every $J_{z_0}^{\varepsilon_b}$-invariant connected complex curve in ${\frak R}_{b,s,t}$ extends to a curve from ${\cal F}_{z_0}^{\frak R}$. We will now describe ${\cal F}_1^{\frak R}$ (here $z_0=1$); for arbitrary $z_0\in{\cal P}$ we have ${\cal F}_{z_0}^{\frak R}=g\left({\cal F}_1^{\frak R}\right)$, where $g\in {\frak V}_b$ is constructed from an element $\tilde g\in\hbox{Aut}({\cal P})$ such that $z_0=\tilde g(1)$. The family ${\cal F}_1^{\frak R}$ consists of the curves
$$
\left\{(z,w)\in\CC^2:(z^2-1)^{k_2}=\rho\,w^{k_1}\right\}\cap {\frak R}_{b,s,t},
$$
where $\rho\in\CC$. Each of these curves is equivalent to either an annulus or a punctured disk. The latter occurs only for $\rho=0$ if either $s=0$ or $t=\infty$, and for $\rho\ne 0$ if $t=\infty$. If either
$s=0$ or $t=\infty$, the corresponding curves accumulate to either the point $(1,\infty)\in{\cal O}_8$ or the point $(1,0)\in{\cal O}_7$, respectively.

Fix $p_0\in O$. Since $\varphi(I_{p_0})$ is a 1-dimensional compact subgroup of ${\frak V}_b$, there is a unique $z_0\in{\cal P}$ such that $\varphi(I_{p_0})=J_{z_0}^{\varepsilon_b}$. Consider any connected $I_{p_0}$-invariant complex curve $C_{p_0}$ in $M$ intersecting $O$ transversally at $p_0$. Since $f(C_{p_0}\setminus\{p_0\})$ is $J_{z_0}^{\varepsilon_b}$-invariant, it extends to a complex curve $C\in{\cal F}_{z_0}^{\frak R}$. If a sequence $\{p_j\}$ in $C_{p_0}\setminus\{p_0\}$ accumulates to $p_0$, the sequence $\{f(p_j)\}$ accumulates to one of the two ends of $C$, and therefore we have either $s=0$ or $t=\infty$.

Assume first that $t=\infty$. In this case, arguing as earlier, we can extend $f$ to a biholomorphic map between $\hat M$ and ${\frak E}_{b,s}$ by setting $f(p_0):=q_0:=(z_0,0)\in {\cal O}_7$, where $p_0$ and $z_0$ are related as specified above. As before, it is straightforward to show that $O$ is the only codimension 2 orbit in $M$, and hence it follows that $M$ is holomorphically equivalent to ${\frak E}_{b,s}$. Similarly, it can be proved that for $s=0$ the manifold $M$ is holomorphically equivalent to ${\cal E}_{b,t}$. As before, this is impossible and thus in case (E) no orbit is a complex curve.       

We now consider the remaining subcases of case (A). Suppose first that there is an orbit in $M$ whose model is either some $\nu_{\alpha}$, or some $\eta_{\alpha}$, or some $\eta_{\alpha}^{(2)}$. It then follows from the poofs of Theorems \ref{classrealhypersurfaces} and \ref{completerealhypersurface} that $M'$ is holomorphically equivalent to one of the following: $\Omega_{s,t}$ with $-1\le s<t\le1$ (see (\ref{omegast})); $D_{s,t}$ with $1\le s<t\le\infty$ (see (\ref{dst})); $D_{s,t}^{(2)}$ with $1\le s<t\le\infty$ (see (\ref{dstinftyn})); ${\frak D}_{s,t}^{(1)}$ with $-1\le s<1<t\le \infty$, where $s=-1$ and $t=\infty$ do not hold simultaneously (see (\ref{frakdst1})).

Suppose first that $M'$ is equivalent to $\Omega_{s,t}$, and let $f$ be an equivalence map. The group ${\cal R}_{\nu}$ (see (\ref{autnu})) acts on the domain $\Omega_1$ (see (\ref{omegat})) with the totally real codimension 2 orbit ${\cal O}_5$ (see (\ref{calo5})). Every 1-dimensional compact subgroup of ${\cal R}_{\nu}$ is the isotropy subgroup of a unique point in ${\cal O}_5$. For $q_0\in{\cal O}_5$ denote by $J_{q_0}^{\nu}$ its isotropy subgroup under the action of ${\cal R}_{\nu}$. There is a family ${\cal F}_{q_0}^{\Omega}$ of connected closed complex curves in $\Omega_{s,t}$ invariant under the $J_{q_0}^{\nu}$-action, such that every connected $J_{q_0}^{\nu}$-invariant complex curve in $\Omega_{s,t}$ extends to a curve from ${\cal F}_{q_0}^{\Omega}$. As before, it is sufficient to describe this family only for a particular choice of $q_0$. The family ${\cal F}_{(0,0)}^{\Omega}$ consists of the connected components of non-empty sets of the form
\begin{equation}
\begin{array}{l}
\left\{(z,w)\in\CC^2:z^2+w^2=\rho\right\}\cap \Omega_{s,t},\\
\left\{(z,w)\in\CC^2: z=iw\right\}\cap \Omega_{s,t},\\
\left\{(z,w)\in\CC^2: z=-iw\right\}\cap \Omega_{s,t},
\end{array}\label{enomega}
\end{equation}
where $\rho\in\CC^*$. Each of the curves from ${\cal F}_{(0,0)}^{\Omega}$ is equivalent to either an annulus or a punctured disk.  The latter is possible only for the last two curves and only for $s=-1$, in which case they accumulate to $(0,0)\in{\cal O}_5$.

Now, arguing as in the second part of case (E) above, we obtain that $s=-1$ and extend $f$ to a map from $\hat M$ onto $\Omega_t$ such that $f(O)={\cal O}_5\subset\Omega_t$. It can be shown, as before, that $f$ is holomorphic on $\hat M$. However, $O$ is a complex curve in $\hat M$ whereas ${\cal O}_5$ is totally real in $\Omega_t$. Hence $M'$ cannot be equivalent to $\Omega_{s,t}$.
 
Assume next that $M'$ is equivalent to $D_{s,t}$ by means of a map $f$. The group $R_{\eta}$ (see (\ref{auteta})) acts on the domain $D_1$ (see (\ref{ds})) with the complex curve orbit ${\cal O}$ (see (\ref{calo})). We again argue as in the second part of case (E) above. Every 1-dimensional compact subgroup of $R_{\eta}$ is the isotropy subgroup of a unique point in ${\cal O}$. For $q_0\in{\cal O}$ denote by $J_{q_0}^{\eta}$ its isotropy subgroup under the action of $R_{\eta}$. There is a family ${\cal F}_{q_0}^{D}$ of connected closed complex curves in $D_{s,t}$ invariant under the $J_{q_0}^{\eta}$-action, such that every connected $J_{q_0}^{\eta}$-invariant complex curve in $D_{s,t}$ extends to a curve from ${\cal F}_{q_0}^{D}$. The family ${\cal F}_{(i,0)}^{D}$ consists of the sets
$$
\begin{array}{l}
\left\{(z,w)\in\CC^2: 1+z^2+\rho w^2=0\right\}\cap D_{s,t},\\
\left\{w=0\right\}\cap D_{s,t},
\end{array}
$$
where $\rho\in\CC$. Each of the curves from ${\cal F}_{(i,0)}^{D}$ is equivalent to an annulus for $t<\infty$ and to a punctured disk if $t=\infty$, in which case it accumulates to $(i,0)\in{\cal O}$.

As before, we now obtain that $t=\infty$ and extend $f$ to a biholomorphic map from $\hat M$ onto $D_s$ such that $f(O)={\cal O}$. It is straightforward to see that $O$ is the only codimension 2 orbit; hence $M$ is holomorphically equivalent to $D_s$, and we have obtained (v) of the theorem.

Suppose now that $M'$ is equivalent to ${\frak D}_{s,t}^{(1)}$, and let $f$ be an equivalence map. The group ${\cal R}^{(1)}$ (see {\bf (10)}) acts on $M^{(1)}$ (see (\ref{m1})) with the complex curve orbit ${\cal O}^{(2)}$ (see (\ref{calo(2)})) and the totally real orbit ${\cal O}_6$ (see (\ref{calo6})). Every 1-dimensional compact subgroup of ${\cal R}^{(1)}$ is the isotropy subgroup of a unique point in each of ${\cal O}^{(2)}$, ${\cal O}_6$. For $q_1\in{\cal O}^{(2)}$ and $q_2\in{\cal O}_6$ that have the same isotropy subgroup under the ${\cal R}^{(1)}$-action, denote this subgroup by $J_{q_1,q_2}^{{\frak D}}$. As before, there is a family ${\cal F}_{q_1,q_2}^{{\frak D}}$ of connected complex closed curves in ${\frak D}_{s,t}^{(1)}$ invariant under the $J_{q_1,q_2}^{{\frak D}}$-action, such that every connected $J_{q_1,q_2}^{{\frak D}}$-invariant complex curve in ${\frak D}_{s,t}^{(1)}$ extends to a curve from ${\cal F}_{q_1,q_2}^{{\frak D}}$. The family ${\cal F}_{(0:1:i:0),(0,0,-i)}^{{\frak D}}$ consists of the connected components of the sets
\begin{equation}
\begin{array}{l}
\left\{(z_1,z_2,z_3)\in\CC^3: z_1^2+z_2^2+\rho z_3^2=0\right\}\cap {\frak D}_{s,t}^{(1)},\\
\left\{z_3=0\right\}\cap {\frak D}_{s,t}^{(1)},\\
\left\{(z_1,z_2,z_3)\in\CC^3: z_1=iz_2\right\}\cap {\frak D}_{s,t}^{(1)},\\
\left\{(z_1,z_2,z_3)\in\CC^3: z_1=-iz_2\right\}\cap {\frak D}_{s,t}^{(1)},\\
\end{array}\label{curvesfrakd}
\end{equation}
where $\rho\in\CC^*$. Each of the sets from ${\cal F}_{(0:1:i:0),(0,0,-i)}^{{\frak D}}$ is equivalent to either an annulus or a punctured disk. If $s>-1$, the latter can only occur for $t=\infty$, in which case the corresponding curves accumulate to $(0:1:i:0)\in{\cal O}^{(2)}$; if $t<\infty$, it occurs only for the last two curves provided $s=-1$, and in this case they accumulate to $(0,0,-i)\in{\cal O}_6$.

It now follows, as before, that either $s=-1$ or $t=\infty$. If $s=-1$ we can extend $f$ to a biholomorphic map between $M'$ and $\hat{\frak D}_t^{(1)}$ (see (\ref{frakdst1})) that takes $O$ onto ${\cal O}_6$. This is impossible since $O$ is a complex curve in $M$ and ${\cal O}_6$ is totally real in $\hat{\frak D}_t^{(1)}$. Hence $t=\infty$, and we can extend $f$ to a biholomorphic map between $\hat M$ and ${\frak D}_s^{(1)}$ (see (\ref{frakdst1})). It is straightforward to see that $O$ is the only codimension 2 orbit in $M$, and thus $M$ is holomorphically equivalent to ${\frak D}_s^{(1)}$, which is a manifold listed in (viii) of the theorem.

Next, the case when $M'$ is equivalent to $D_{s,t}^{(2)}$ is treated as the preceding one. Here we parametrize maximal compact subgroups of ${\cal R}^{(1)}$ by points in ${\cal O}^{(2)}$, and for the point $(0:1:i:0)\in{\cal O}^{(2)}$ the corresponding family of complex curves consists of sets constructed as family (\ref{curvesfrakd}), where the curves appearing on the left must be intersected with $D_{s,t}^{(2)}$ rather than ${\frak D}_{s,t}^{(1)}$ (note, however, that the second last intersection is empty). As above, we obtain that $t=\infty$ and that $M$ is holomorphically equivalent to $D_s^{(2)}$ (see (\ref{ds2n})). We now recall that $D_1^{(2)}$ is equivalent to $\Delta^2$ (see {\bf (11)(c)}), and, excluding the value $s=1$, obtain (vi) of the theorem. 

We now assume that $M'$ is holomorphically equivalent to one of the $n$-sheeted covers, for $n\ge 2$, of the previously considered possibilities: $\Omega_{s,t}^{(n)}$ (the cover of $\Omega_{s,t}$) with $-1\le s<t\le1$ -- see (\ref{omegastinftyn}); $D_{s,t}^{(n)}$ (the cover of $D_{s,t}$) with $1\le s<t\le\infty$, where $n$ is odd -- see {\bf (10)}; $D_{s,t}^{(2n)}$ (the cover of $D_{s,t}^{(2)}$) with $1\le s<t\le\infty$ -- see (\ref{dstinftyn}); ${\frak D}_{s,t}^{(n)}$ (the cover of ${\frak D}_{s,t}^{(1)}$) with $-1\le s<1<t\le \infty$, where $s=-1$ and $t=\infty$ do not hold simultaneously -- see {\bf (11)(d)}. We will now formulate a number of useful properties that hold for the covers. These properties (that we hereafter refer to  as Properties (P)) follow from the explicit construction of the covers in {\bf (10)}, {\bf (11)}. 

Let $S$ be one of $\Omega_{s,t}$, $D_{s,t}$, $D_{s,t}^{(2)}$, ${\frak D}_{s,t}^{(1)}$ and let $S^{(n)}$ be the corresponding $n$-sheeted cover of $S$ (for $S=D_{s,t}$ we assume that $n$ is odd). Let $H:=G(S)$ and $H^{(n)}:=G\left(S^{(n)}\right)$. Then we have:
\vspace{0.3cm}

(a) the group $H^{(n)}$ consists of all lifts from $S$ to $S^{(n)}$ of all elements of $H$, and the natural projection $\pi: H^{(n)}\ra H$ is a Lie group homomorphism and realizes $H^{(n)}$ as an $n$-sheeted cover of $H$;

(b) it follows from (a) that for every maximal compact subgroup $K_0\subset H$ (all such subgroups are isomorphic to $U_1$) the subgroup $\pi^{-1}(K_0)$ is maximal compact in $H^{(n)}$, and all maximal compact subgroups of $H^{(n)}$ are obtained in this way;

(c) for every maximal compact subgroup $K\subset H^{(n)}$ the family of all $K$-invariant complex curves in $S^{(n)}$ consists of the lifts from $S$ to $S^{(n)}$ of all $\pi(K)$-invariant complex curves in $S$, where every connected $\pi(K)$-invariant curve $C$ is lifted to a unique connected $K$-invariant curve $C^{(n)}$ (in particular, $C^{(n)}$ covers $C$ in an $n$-to-1 fashion);

(d) if $S$ is one of $D_{s,t}$, $D_{s,t}^{(2)}$, ${\frak D}_{s,t}^{(1)}$, then every maximal compact subgroup $K\subset H^{(n)}$ is the isotropy subgroup -- with respect to the $H^{(n)}$-action -- of a unique point in ${\cal O}^{(n)}$ (see {\bf (11)(b)}) in the first case, and a unique point in ${\cal O}^{(2n)}$ (see (\ref{cal02n})) in each of the other two cases; every $K$-invariant closed complex curve in $S^{(n)}$ equivalent to a punctured disk accumulates to this point (provided, for $S={\frak D}_{s,t}^{(1)}$, we assume that $s>-1$).  
\vspace{0.3cm}

Properties (P) yield that if $M'$ is equivalent to either $D_{s,t}^{(n)}$ for odd $n$ or $D_{s,t}^{(2n)}$ for $n\ge 2$, then $t=\infty$ and $M$ is holomorphically equivalent to either $D_s^{(n)}$ (see (\ref{dsn})) or $D_s^{(2n)}$ (see (\ref{ds2n})), respectively; this gives (vii) of the theorem.

Suppose now that $M'$ is equivalent to $\Omega_{s,t}^{(n)}$ by means of a map $f$. Then Properties (P) imply that $s=-1$. Recall that $\Psi_{\nu}\circ\Phi_{\nu}^{(n)}:\Omega_{-1,t}^{(n)}\ra\Omega_{-1,t}$ is an $n$-to-1 covering map (see {\bf (10)}). Consider the composition $\tilde f:=\Psi_{\nu}\circ\Phi_{\nu}^{(n)}\circ f$. This is an $n$-to-1 covering map from $M'$ onto $\Omega_{-1,t}$ satisfying (\ref{equivar}) for all $g\in G(M)$, $q\in M'$, where $\varphi: G(M)\ra {\cal R}_{\nu}$ is an $n$-to-1 covering homomorphism. Fix $p_0\in O$. Since $K_0:=\varphi(I_{p_0})$ is a maximal compact subgroup of ${\cal R}_{\nu}$, there is a unique $q_0\in{\cal O}_5$ such that $K_0$ is the isotropy subgroup of $q_0$ under the ${\cal R}_{\nu}$-action on $\Omega_1$. We now define $\tilde f(p_0):=q_0$. Thus, we have extended $\tilde f$ to an equivariant map from $\hat M$ onto $\Omega_t$ that takes $O$ onto ${\cal O}_5$. As before, it can be shown that $\tilde f$ is holomorphic on $\hat M$. However, ${\cal O}_5$ is totally real in $\Omega_{-1,t}$ and therefore $M'$ cannot in fact be equivalent to $\Omega_{s,t}^{(n)}$.

Let $M'$ be equivalent to ${\frak D}_{s,t}^{(n)}$, and let $f$ be an equivalence map. In this case Properties (P) imply that either $s=-1$ or $t=\infty$. If $s=-1$, arguing as in the preceding paragraph, we extend the map $\tilde f:=\Phi^{(n)}\circ f$ to a holomorphic map from $\hat M$ onto $\hat{\frak D}_t^{(1)}$ that takes $O$ onto ${\cal O}_6$. As before, this is impossible since ${\cal O}_6$ is totally real in $\hat{\frak D}_t^{(1)}$, and therefore we in fact have $t=\infty$. In this case Properties (P) yield that $M$ is holomorphically equivalent to ${\frak D}_s^{(n)}$, and we have obtained (viii) of the theorem.  

We now assume that every codimension 2 orbit in $M$ is totally real. We will go again through all the possibilities for the group $G(M)$ listed in Lemma \ref{groupsofmodels}, paying attention to constraints imposed on $G(M)$ by this condition. In what follows $O$ denotes a totally real orbit in $M$. In case (E) with $G(M)$ isomorphic to $G_{\varepsilon_b}$ (see (\ref{gvarepsilonb})) we obtain, as before, that $I_p=I:=\varphi^{-1}(J^{\varepsilon_b})$ for every $p\in O$, and thus $I_p$ acts trivially on $O(p)$ for every $p\in O$ which contradicts (iv) of Proposition \ref{dim}. A similar argument gives a contradiction in case (G). In case (E) with $G(M)$ isomorphic to ${\frak V}_b$ the argument given above for the case of complex curve orbits shows that $f$ extends to a biholomorphic map between $\hat M$ and either ${\frak E}_{b,s}$ (see (\ref{frakebs})) or domain (\ref{prelim}), with either $f(O)={\cal O}_7$ or $f(O)={\cal O}_8$, respectively, which is impossible, since $O$ is totally real, whereas ${\cal O}_6$, ${\cal O}_7$ are complex curves. 
Next, in subcase (l') of case (A) the group $G(M)$ is isomorphic to $SU_2$, which implies that $M$ is holomorphically equivalent to one of the manifolds listed in \cite{IKru2}. However, none of the manifolds on the list has a totally real orbit. Therefore, it remains to consider subcases (j), (j''), (l), (l''), (m), (m''), (n), (n'') of case (A), and case (C).

We start with case (C). In this situation $G(M)$ is isomorphic to $SU_2$, if $m$ is odd and to $SU_2/\left\{\pm\hbox{id}\right\}$, if $m$ is even. To rule out the case of odd $m$ we again use the result of \cite{IKru2}. We now assume that $m$ is even. In this case $M'$ is holomorphically equivalent to $S_{m,s,t}$ (see (\ref{smst})), with $0\le s<t<\infty$, by means of a map $f$ that satisfies (\ref{equivar}) for all $g\in G(M)$, $q\in M'$ and some isomorphism $\varphi:G(M)\ra SU_2/\left\{\pm\hbox{id}\right\}$.  

Fix $p_0$ in $O$ and consider the connected compact 1-dimensional subgroup $\varphi(I_{p_0}^0)\subset SU_2/\left\{\pm\hbox{id}\right\}$. 
It then follows that $\varphi(I_{p_0}^0)$ is conjugate in $SU_2/\left\{\pm\hbox{id}\right\}$ to the subgroup $J^{{\cal L}}$ that consists of all elements of the form
$$
\left(
\begin{array}{ll}
e^{i\psi} & 0\\
0 & e^{-i\psi}
\end{array}
\right) \left\{\pm\hbox{id}\right\},
$$
where $\psi\in\RR$ (see e.g. Lemma 2.1 of \cite{IKru1}). Suppose that $p_0$ is chosen so that $\varphi(I_{p_0}^0)=J^{{\cal L}}$. Let $C_{p_0}$ be a connected $I_{p_0}^0$-invariant complex curve in $M$ that intersects $O$ transversally at $p_0$. Then $f(C_{p_0}\setminus\{p_0\})$ is a connected $J^{{\cal L}}$-invariant complex curve in $S_{m,s,t}$. It is straightforward to see that every connected $J^{{\cal L}}$-invariant complex curve in $S_{m,s,t}$ extends to a closed curve equivalent to either an annulus or a punctured disk. The only closed connected $J^{{\cal L}}$-invariant curves in $S_{m,s,t}$ that can be equivalent to a punctured disk (which only occurs for $s=0$) are
$$
\Bigl(\{z=0\}/\ZZ_m\Bigr)\cap S_{m,s,t}
$$
and
$$
\Bigl(\{w=0\}/\ZZ_m\Bigr)\cap S_{m,s,t}.
$$
Therefore, the curve $f(C_{p_0}\setminus\{p_0\})$ extends to one of these curves, and we have $s=0$.

Let $B_t$ be the ball of radius $t$ in $\CC^2$ and   $\widehat{B_t}$  its blow-up at the origin, i.e.
$$
\widehat
{B_t}:=\left\{\Bigl[(z,w),(\xi:\zeta)\Bigr]\in B_t\times\CC\PP^1:z\zeta=w\xi\right\},
$$
where $(\xi:\zeta)$ are the homogeneous coordinates in
$\CC\PP^1$. We define an action of $U_2$ on $\widehat{B_t}$ as
follows: for $g\in U_2$ and $\Bigl[(z,w),(\xi:\zeta)\Bigr]\in\widehat{B_t}$ set
$$
g\Bigl[(z,w),(\xi:\zeta)\Bigr]:=\Bigl[g(z,w),g(\xi:\zeta)\Bigr],
$$
where in the right-hand side we use the standard actions of $U_2$ on $\CC^2$ and $\CC\PP^1$. Next, we denote by $\widehat{B_t}/\ZZ_m$ the quotient of $\widehat {B_t}$ by the equivalence relation $\Bigl[(z,w),(\xi:\zeta)\Bigr]\sim e^{\frac{2\pi
    i}{m}}\Bigl[(z,w),(\xi:\zeta)\Bigr]$. Let $\left\{\Bigl[(z,w),(\xi:\zeta)\Bigr]\right\}\in\widehat{B_t}/\ZZ_m$ be  the
equivalence class of $\Bigl[(z,w),(\xi:\zeta)\Bigr]\in\widehat{B_t}$. We now define in a natural way an action of
$SU_2/\left\{\pm\hbox{id}\right\}$ on $\widehat{B_t}/\ZZ_m$: for
$\left\{\Bigl[(z,w),(\xi:\zeta)\Bigr]\right\}\in\widehat{B_t}/\ZZ_m$ and $g\left\{\pm\hbox{id}\right\}\in SU_2/\left\{\pm\hbox{id}\right\}$ we set
$$
g\left\{\pm\hbox{id}\right\}\,\left\{\Bigl[(z,w),(\xi:\zeta)\Bigr]\right\}:=\left\{g\Bigl[(z,w),(\xi:\zeta)\Bigr]\right\}.
$$
The points $\left\{\Bigl[(0,0),(\xi:\zeta)\Bigr]\right\}$ form an $SU_2/\left\{\pm\hbox{id}\right\}$-orbit that we denote by ${\cal O}_9$; this is a complex curve equivalent to $\CC\PP^1$. Everywhere below we identify $S_{m,0,t}$ with $\widehat{B_t}/\ZZ_m\setminus{\cal O}_9$.

For $q_0\in{\cal O}_9$ let $J_{q_0}^{{\cal L}}$ be the isotropy subgroup under the action of\linebreak $SU_2/\left\{\pm\hbox{id}\right\}$. It is straightforward to see that every subgroup $J_{q_0}^{{\cal L}}$ is conjugate to $J^{{\cal L}}$ in $SU_2/\left\{\pm\hbox{id}\right\}$ and that for every $q_0$ there is exactly one $q_0'\in{\cal O}_9$, $q_0'\ne q_0$, such that $J_{q_0}^{{\cal L}}=J_{q_0'}^{{\cal L}}$ (for example, $J^{{\cal L}}$ is the isotropy subgroup of each of $\left\{\Bigl[(0,0),(1:0)\Bigr]\right\}$ and $\left\{\Bigl[(0,0),(0:1)\Bigr]\right\}$). Fix $q_0\in{\cal O}_9$ and let $p_0\in O$ be such that $\varphi(I_{p_0}^0)=J_{q_0}^{{\cal L}}$. As we noted at the beginning of the proof of the theorem, there are exactly two connected $I_{p_0}^0$-invariant complex curves $C_{p_0}$ and $\tilde C_{p_0}$ in a neighborhood of $p_0$ that intersect $O$ at $p_0$ transversally. The curves $f(C_{p_0}\setminus\{p_0\})$ and $f(\tilde C_{p_0}\setminus\{p_0\})$ extend to the two distinct closed $J_{q_0}^{{\cal L}}$-invariant complex curves in $\widehat{B_t}/\ZZ_m\setminus{\cal O}_9$ that are equivalent to a punctured disk. Since there are no other closed $J_{q_0}^{{\cal L}}$-invariant complex curves in $\widehat{B_t}/\ZZ_m\setminus{\cal O}_9$ equivalent to a punctured disk, it follows that $I_{p_0'}^0\ne I_{p_0}^0$ for every $p_0'\in O$, $p_0'\ne p_0$.

Observe that if $q,q'\in{\cal O}_9$, $q\ne q'$, are such that $J_{q}^{{\cal L}}=J_{q'}^{{\cal L}}=:J$, then one of the $J$-invariant complex curves equivalent to a punctured disk accumulates to $q$ and the other to $q'$. Therefore, we can extend ${\frak F}:=f^{-1}$ to a map from $\widehat{B_t}/\ZZ_m$ onto $\hat M$ by setting ${\frak F}(q_0):=p_0$, where $q_0\in {\cal O}_9$ and $p_0\in O$ are related as indicated above (hence ${\frak F}$ is 2-to-1 on ${\cal O}_9$). As before, it can be shown that ${\frak F}$ is continuous on $\widehat{B_t}/\ZZ_m$ and thus is holomorphic there. However, ${\frak F}$ maps the complex curve ${\cal O}_9\subset\widehat{B_t}/\ZZ_m$ onto the totally real submanifold $O\subset\hat M$, which is impossible. Hence, $M'$ cannot be equivalent to $S_{m,s,t}$.

We now consider the remaining subcases of case (A). In subcase (j) the manifold $M'$ is holomorphically equivalent to ${\frak S}_{s,t}$ for $0\le s<t<\infty$ (see (\ref{frakst})). The group ${\cal R}_{\chi}$ (see (\ref{autchi})) acts on $\CC^2$ with the only codimension 2 orbit ${\cal O}_2$ (see (\ref{calo2})). The isotropy subgroup of a point $(i y_0,i v_0)\in{\cal O}_2$ is the group $J_{(y_0,v_0)}^{\chi}$ of all transformations of the form (\ref{autchi}) with $\beta=y_0-\cos\psi\cdot y_0-\sin\psi\cdot v_0$, $\gamma=v_0+\sin\psi\cdot y_0-\cos\psi\cdot v_0$, where
$$
A=\left(
\begin{array}{rr}
\cos\psi & \sin\psi\\
-\sin\psi & \cos\psi
\end{array}
\right),
$$
with $\psi\in\RR$. Note that these subgroups are maximal compact in ${\cal R}_{\chi}$ (which implies that $I_p$ is connected for every $p\in O$), and the isotropy subgroups of distinct points in ${\cal O}_2$ do not coincide.

We now argue as in the second part of case (E) for complex curve orbits. There is a family ${\cal F}_{(y_0,v_0)}^{\frak S}$ of connected closed complex curves in ${\frak S}_{s,t}$ invariant under the $J_{(y_0,v_0)}^{\chi}$-action, such that every connected $J_{(y_0,v_0)}^{\chi}$-invariant complex curve in ${\frak S}_{s,t}$ extends to a curve from ${\cal F}_{(y_0,v_0)}^{\frak S}$. The family ${\cal F}_{(0,0)}^{\frak S}$ consists of the connected components of non-empty sets analogous to (\ref{enomega}), where the sets on the left must be intersected with ${\frak S}_{s,t}$ rather than $\Omega_{s,t}$. Among the curves from ${\cal F}_{(0,0)}^{\frak S}$, only the last two can be equivalent to a punctured disk. This occurs only for $s=0$, in which case the curves accumulate to $(0,0)\in{\cal O}_2$. Arguing as before, we can now construct a biholomorphic map between $M$ and ${\frak S}_t$ (see (\ref{frakt})). Clearly, ${\frak S}_t$ is equivalent to ${\frak S}_1$, and we have obtained (i) of the theorem.

Consider subcase (j''). In this situation $M'$ is holomorphically equivalent to the $n$-sheeted cover ${\frak S}_{s,t}^{(n)}$ of ${\frak S}_{s,t}$ for $0\le s<t<\infty$, $n\ge 2$ (see (\ref{frakstinftyn})), by means of a map $f$. From the explicit construction of the covers in {\bf (4)} it follows that Properties (P) hold for $S={\frak S}_{s,t}$. Let $\tilde f:=\Phi_{\chi}^{(n)}\circ f$, where $\Phi_{\chi}^{(n)}:{\frak S}_{s,t}^{(n)}\ra {\frak S}_{s,t}$ is the $n$-to-1 covering map defined in (\ref{phichin}). Arguing as in the case of complex curve orbits when $M'$ was assumed to be equivalent to $\Omega_{s,t}^{(n)}$, we extend $\tilde f$ to a holomorphic map from $\hat M=M'\cup O$ onto ${\frak S}_t$, that takes $O$ onto ${\cal O}_2$. 

Suppose that the differential of $\tilde f$ is degenerate at a point in $O$. Then, since $\tilde f$ satisfies (\ref{equivar}), its differential degenerates everywhere on $O$. Since $O$ is totally real, it follows that the differential of $\tilde f$ is degenerate everywhere in $\hat M$. This is impossible since $\tilde f$ is a covering map on $M'$, and thus $\tilde f$ is non-degenerate at every point of $O$. Hence, for every $p\in O$ there exists a neighborhood of $p$ in which $\tilde f$ is biholomorphic. Fix $p_0\in O$ and let $C_{p_0}$ be a connected $I_{p_0}$-invariant complex curve intersecting $O$ at $p_0$ transversally (observe that $I_{p_0}$ is connected). Then it follows from (c) of Properties (P) that $f(C_{p_0}\setminus\{p_0\})$ covers $\tilde f(C_{p_0}\setminus\{p_0\})$ in an $n$-to-1 fashion, and hence $\tilde f$ cannot be biholomorphic in any neighborhood of $p_0$. This contradiction yields that $M'$ cannot be equivalent to ${\frak S}_{s,t}^{(n)}$.

Consider subcase (l). In this situation $M'$ is holomorphically equivalent to $E_{s,t}$ for $1\le s<t<\infty$ (see (\ref{est})). The group $R_{\mu}$ (see (\ref{autmu})) acts on $\CC\PP^2$ with the totally real orbit ${\cal O}_3$ (see (\ref{calo3})). Every connected 1-dimensional compact subgroup of ${\cal R}_{\mu}$ is conjugate in ${\cal R}_{\mu}$ to the subgroup $J^{\mu}$ that consists of all matrices of the form
\begin{equation}
\left(
\begin{array}{ccc}
1 & 0 & 0\\
0 & \cos\psi & \sin\psi\\
0 & -\sin\psi & \cos\psi\\
 \end{array}
 \right),\label{stabilizer}
\end{equation}
where $\psi\in\RR$ (this follows, for instance, from Lemma 2.1 of \cite{IKru1}). It is straightforward to see that the isotropy subgroup $J^{\mu}_q$  of a point $q\in{\cal O}_3$ under the ${\cal R}_{\mu}$-action is conjugate to $J^{\mu}$ (note that $J^{\mu}=J^{\mu}_{(1:0:0)}$) and that the isotropy subgroups of distinct points do not coincide.

There is a family ${\cal F}_q^E$ of connected closed complex curves in $E_{s,t}$ invariant under the $J^{\mu}_q$-action, such that every connected $J^{\mu}_q$-invariant complex curve in $E_{s,t}$ extends to a curve from ${\cal F}_q^E$. The family ${\cal F}_{(1:0:0)}^E$ consists of the connected components of non-empty sets of the form 
$$
\begin{array}{l}
\left\{(\zeta:z:w)\in\CC\PP^2:z^2+w^2=\rho\zeta^2\right\}\cap E_{s,t},\\
\left\{(\zeta:z:w)\in\CC\PP^2: z=iw\right\}\cap E_{s,t},\\
\left\{(\zeta:z:w)\in\CC\PP^2: z=-iw\right\}\cap E_{s,t},
\end{array}
$$
where $\rho\in\CC^*$. Among the curves from ${\cal F}_{(1:0:0)}^E$, only the last two can be equivalent to a punctured disk. This occurs only for $s=1$, in which case the curves accumulate to $(1:0:0)\in{\cal O}_3$. Arguing as before, we can now construct a biholomorphic map between $M$ and $E_t$ (see (\ref{et})), which gives (ii) of the theorem.

Further, in subcase (l'') $M'$ is holomorphically equivalent to $E_{s,t}^{(2)}$ for some $1\le s<t<\infty$ (see (\ref{est24})). Let $f$ be an equivalence map that satisfies (\ref{equivar}) for all $g\in G(M)$, $q\in M'$ and some isomorphism $\varphi: G(M)\ra {\cal R}_{\mu}^{(2)}$ (see (\ref{autmu2})). The group ${\cal R}_{\mu}^{(2)}$ acts on ${\cal Q}_{+}$ (see (\ref{quadricplus})) with the totally real orbit ${\cal O}_4$ (see (\ref{calo4})). All 1-dimensional compact subgroups are described as in subcase (l) -- see (\ref{stabilizer}). The isotropy subgroup $J^{\mu^{(2)}}_q$ of a point $q\in{\cal O}_4$ under the ${\cal R}_{\mu}^{(2)}$-action is conjugate to $J^{\mu}$, and for every $q\in{\cal O}_4$ there exists exactly one $q'\in{\cal O}_4$, $q'\ne q$, for which $J^{\mu^{(2)}}_q=J^{\mu^{(2)}}_{q'}$ (note that $q'=-q$ and $J^{\mu}=J^{\mu^{(2)}}_{(\pm1,0,0)}$).

Again, there is a family ${\cal F}_q^{E^{(2)}}$ of connected closed complex curves in $E_{s,t}^{(2)}$ invariant under the $J^{\mu^{(2)}}_q$-action, such that every connected $J^{\mu^{(2)}}_q$-invariant complex curve in $E_{s,t}^{(2)}$ extends to a curve from ${\cal F}_q^{E^{(2)}}$. The family ${\cal F}_{(\pm1,0,0)}^{E^{(2)}}$ consists of the connected components of non-empty sets of the form 
$$
\begin{array}{lll}
&&\left\{(z_1,z_2,z_3)\in\CC^3:z_2^2+z_3^2=\rho z_1^2\right\}\cap E_{s,t}^{(2)},\\
{\frak C}_1&:=&\left\{(z_1,z_2,z_3)\in\CC^3: z_1=1,\,z_2=iz_3\right\}\cap E_{s,t}^{(2)},\\
{\frak C}_2&:=&\left\{(z_1,z_2,z_3)\in\CC^3: z_1=1,\,z_2=-iz_3\right\}\cap E_{s,t}^{(2)},\\
{\frak C}_3&:=&\left\{(z_1,z_2,z_3)\in\CC^3: z_1=-1,\,z_2=iz_3\right\}\cap E_{s,t}^{(2)},\\
{\frak C}_4&:=&\left\{(z_1,z_2,z_3)\in\CC^3: z_1=-1,\,z_2=-iz_3\right\}\cap E_{s,t}^{(2)},
\end{array}
$$
where $\rho\in\CC^*$. Among the curves from ${\cal F}_{(\pm1,0,0)}^{E^{(2)}}$, only ${\frak C}_j$ can be equivalent to a punctured disk, which occurs only for $s=1$. It then follows that $s=1$, and in this case ${\frak C}_1$, ${\frak C}_2$ accumulate to $(1,0,0)\in{\cal O}_4$, while ${\frak C}_3$, ${\frak C}_4$ accumulate to $(-1,0,0)\in{\cal O}_4$.

Fix $p_0\in O$ and let $q_0\in{\cal O}_4$ be such that $\varphi(I_{p_0}^0)=J^{\mu^{(2)}}_{q_0}$ and such that, for a $I_{p_0}^0$-invariant complex curve $C_{p_0}$ intersecting $O$ at $p_0$ transversally, the curve $f\left(C_{p_0}\setminus\{p_0\}\right)$ extends to a complex curve from ${\cal F}_{q_0}^{E^{(2)}}$ that accumulates to $q_0$. We extend ${\frak F}:=f^{-1}$ to a map from $E_t^{(2)}$ (see (\ref{et2})) onto $\hat M=M'\cup O$ that takes ${\cal O}_4$ onto $O$. Define ${\frak F}(q_0):=p_0$ and for any $h\in{\cal R}_{\mu}^{(2)}$ set ${\frak F}(hq_0):=\varphi^{-1}(h)p_0$. Since $\varphi^{-1}\left(J^{\mu^{(2)}}_{q_0}\right)\subset I_{p_0}$, this map is well-defined. Furthermore, the extended map satisfies (\ref{equivar1}) for all $h\in{\cal R}_{\mu}^{(2)}$, $q\in E_t^{(2)}$, and for every $q\in{\cal O}_4$ there exists a $J^{\mu^{(2)}}_q$-invariant complex curve ${\frak C}$ in $E_t^{(2)}$ that intersects ${\cal O}_4$ at $q$ transversally and such that ${\frak F}\left({\frak C}\setminus\{q\}\right)$ is an $I_{{\frak F}(q)}^0$-invariant complex curve that accumulates to ${\frak F}(q)$. Arguing as in the second part of case (E) for complex curve orbits, we now obtain that ${\frak F}$ is holomorphic on $E_t^{(2)}$. Further, as in subcase (j'') above, we see that ${\frak F}$ is locally biholomorphic in a neighborhood of every point in ${\cal O}_4$.

We will now show that ${\frak F}$ is 1-to-1 on ${\cal O}_4$. Suppose that for some $q,q'\in{\cal O}_4$, $q\ne q'$, we have ${\frak F}(q)={\frak F}(q')=p$ for some $p\in O$. Since ${\frak F}$ satisfies (\ref{equivar1}), we have $J^{\mu^{(2)}}_q=J^{\mu^{(2)}}_{q'}=\varphi(I_p^0)$, and therefore $q'=-q$. Consider the four $J^{\mu^{(2)}}_q$-invariant connected complex curves in $E_{1,t}^{(2)}$ equivalent to a punctured disk; a pair of these curves accumulates to $q$, while the other pair accumulates to $-q$. The curves are mapped by ${\frak F}$ into four distinct $I_p^0$-invariant complex curves in $M'$ whose extensions in $\hat M$ intersect $O$ transversally at $p$. However, as we noted at the beginning of the proof of the theorem, there are exactly two $I_p^0$-invariant complex curves near $p$ that intersect $O$ transversally at $p$. This contradiction yields that ${\frak F}$ is a biholomorphic map from $E_t^{(2)}$ onto $\hat M$. It can be now shown, as before, that $O$ is the only codimension 2 orbit in $M$, which gives that $M$ is holomorphically equivalent to $E_t^{(2)}$, and we have obtained (iii) of the theorem.

It now remains to consider subcases (m), (m''), (n), (n''). We will proceed as in the situation when a complex curve orbit was assumed to be present in $M$. If $M'$ is equivalent to one of $D_{s,t}$, $D_{s,t}^{(n)}$ for $n\ge 2$, ${\frak D}_{s,t}^{(n)}$ for $n\ge 1$ (where in the last case we assume that $s>-1$), we obtain a contradiction since $O$ is totally real in $M$ whereas ${\cal O}$ and ${\cal O}^{(n)}$ for $n\ge 2$ are complex curves in the corresponding manifolds. Further, if $M'$ is equivalent to $\Omega_{s,t}$, we obtain that $s=-1$ and $M$ is holomorphically equivalent to $\Omega_t$. Recalling that $\Omega_1$ is equivalent to $\Delta^2$ (see {\bf (11)(c)}) and excluding the value $t=1$, we obtain (iv) of the theorem. Next, If $M'$ is equivalent to ${\frak D}_{-1,t}^{(1)}$, then $M$ is equivalent to $\hat{\frak D}_t^{(1)}$, which are the manifolds in (ix) of the theorem.

Suppose now that $M'$ is equivalent to $\Omega_{s,t}^{(n)}$ for some $-1\le s<t\le\infty$. In this case we obtain a holomorphic map $\tilde f$ from $\hat M$ onto $\Omega_t$ that takes $O$ onto ${\cal O}_5$ and such that $\tilde f|_{M'}$ is an $n$-to-1 covering map from $M'$ onto $\Omega_{-1,t}$. Now, arguing as in subcase (j''), we obtain that the differential of $\tilde f$ is non-degenerate at every point of $O$ which leads to a contradiction. Finally, a similar argument leads to a contradiction if $M'$ is equivalent to ${\frak D}_{-1,t}^{(n)}$ for $n\ge 2$.   

The proof of the theorem is complete.\qed

{\obeylines
Department of Mathematics
The Australian National University
Canberra, ACT 0200
AUSTRALIA
E-mail: alexander.isaev@maths.anu.edu.au
}

\end{document}